\documentclass[eulercal, sepfigandtablenums, tikz]{nmd/nmd-article}
\pdfoutput=1
\ifnmd 
\else 
\fi

\usepackage{accents}
\usepackage{enumitem}
\usepackage{placeins}

\usepackage{array}
\usepackage{multirow}
\usepackage{makecell}
\newcommand{\tablefont}{\small}

\newcommand{\interior}[1]{{#1}^\circ}

\renewcommand{\bdry}{$\partial$\hyp}

\def\R{\mathbb{R}}
\def\Q{\mathbb{Q}}

\def\BQ{\mathbb Q}

\def\calG{\cG}

\def\calT{\cT}

\def\S{\Sigma}

\def\t{\tau}

\def\tight{T}

\def\LA{LA}

\def\be{  \begin{equation} }
\def\ee{  \end{equation} }

\newcommand{\wt}{\widetilde}

\renewcommand{\T}{\mathcal{T}}
\def\Z{\mathbb Z}

\def\C{\mathbb C}

\DeclareMathOperator{\weight}{wt}

\def\ML{\mathcal{M}\mathcal{L}}

\newcommand{\vE}{\vec{E}}
\newcommand{\vF}{\vec{F}}
\newcommand{\vG}{\vec{G}}
\newcommand{\vH}{\vec{H}}
\newcommand{\vK}{\vec{K}}

\newcommand{\va}{\vec{a}}
\newcommand{\vb}{\vec{b}}
\newcommand{\vc}{\vec{c}}
\newcommand{\vv}{\vec{v}}
\newcommand{\ww}{\vec{w}}
\newcommand{\xx}{\vec{x}}
\newcommand{\yy}{\vec{y}}

\newcommand{\evec}{\vec{e}}
\newcommand{\bt}{\mathbf{t}}

\newcommand{\LW}[1][\T]{\mathcal{LW}_{#1}}
\newcommand{\CLW}[1][\T]{\mathcal{CLW}_{#1}}
\newcommand{\PT}[1][\T]{\mathcal{P}_{#1}}
\newcommand{\ST}[1][\T]{\mathcal{S}_{#1}}
\newcommand{\AT}[1][\T]{\mathcal{A}_{#1}}
\newcommand{\NT}[1][\T]{\mathcal{N}_{#1}}
\newcommand{\intC}{C^{\hspace{0.03em} \circ}}
\newcommand{\intCtil}{\Ctil^{\hspace{0.03em} \circ}}
\newcommand{\intD}{\interior{D}}
\newcommand{\intX}{\interior{X}}
\newcommand{\intE}{E^{\hspace{0.05em} \circ}}
\newcommand{\dep}[1]{{\mathrm{dep}(#1)}}
\newcommand{\Rp}{\R_{\geq 0}}

\begin{document}


\title{Counting essential surfaces \\ in 3-manifolds}

\author{Nathan M. Dunfield}
\address{Department of Mathematics\\
         University of Illinois\\
         1409 W. Green Street \\
         Urbana, IL 61801, USA}
\email{nathan@dunfield.info}
\urladdr{http://dunfield.info}

\author{Stavros Garoufalidis}
\address{
  International Center for Mathematics, Department of Mathematics \\
  Southern University of Science and Technology \\
  Shenzhen, CHINA}
\email{stavros@mpim-bonn.mpg.de}
\urladdr{http://people.mpim-bonn.mpg.de/stavros}

\author{J. Hyam Rubinstein}
\address{School of Mathematics and Statistics \\
         The University of Melbourne \\
         Parkville, VIC, 3010, AUSTRALIA}
\email{joachim@unimelb.edu.au}

\begin{abstract}
  We consider the natural problem of counting isotopy classes of
  essential surfaces in 3-manifolds, focusing on closed essential
  surfaces in a broad class of hyperbolic \3-manifolds.  Our main
  result is that the count of (possibly disconnected) essential
  surfaces in terms of their Euler characteristic always has a short
  generating function and hence has quasi-polynomial behavior.  This
  gives remarkably concise formulae for the number of such surfaces,
  as well as detailed asymptotics.  We give algorithms that allow us
  to compute these generating functions and the underlying surfaces,
  and apply these to almost 60,000 manifolds, providing a wealth of
  data about them. We use this data to explore the delicate question
  of counting only connected essential surfaces and propose some
  conjectures.  Our methods involve normal and almost normal surfaces,
  especially the work of Tollefson and Oertel, combined with
  techniques pioneered by Ehrhart for counting lattice points in
  polyhedra with rational vertices. We also introduce a new way of
  testing if a normal surface in an ideal triangulation is essential
  that avoids cutting the manifold open along the surface; rather, we
  use almost normal surfaces in the original triangulation.
\end{abstract}
%

%

\renewcommand{\andify}[1]{\relax}
\renewcommand{\authors}{%
\noindent
Nathan M. Dunfield, 
Stavros Garoufalidis,

\vspace{0.4em}

\noindent
and J. Hyam Rubinstein 
}

\maketitle

\tableofcontents


\section{Introduction}
\label{sec.introduction}

Essential surfaces have played a central role in 3-manifold topology
for at least the last 70 years, being both a key tool and a
fundamental object of study.  Roughly, these are compact embedded
surfaces $F \subset M^3$ where $\pi_1 F \to \pi_1M$ is injective;
throughout this introduction, see Section~\ref{sec.background} for
precise definitions and conventions.  While some compact \3-manifolds
contain no essential surfaces at all (the \3-sphere, lens spaces),
others contain infinitely many isotopy classes of essential surfaces
of the same topological type (the \3-torus contains infinitely many
essential \2-tori).  However, for $M$ that are irreducible and
atoroidal (i.e.~contain no essential spheres or tori), the number of
essential $F$ of a fixed topological type is always finite
\cite[Corollary~2.3]{Jaco-Oertel}.  For example, any hyperbolic
\3-manifold is irreducible and atoroidal, and these form the main class
of interest here. 

A natural problem is thus to describe in a structured way the set of
essential surfaces in a given \3-manifold $M$, in particular to list
and to count them.  Focusing on those $F$ that are closed, connected,
and orientable, define $a_M(g)$ to be the number of isotopy classes of
essential surfaces in $M$ of genus $g$.  There are plenty of
hyperbolic 3-manifolds where $a_M(g) = 0$ for all $g$, including all
those that are exteriors of 2-bridge knots \cite{HatcherThurston1985}.
In contrast, for the exterior of $X$ of the Conway knot $K11n34$, we
can use Theorem~\ref{thm.algol} below to compute the values of
$a_X(g)$ shown in Table~\ref{tab:conway}, as well as further values
such as $a_X(50) = 56,892$ and $a_X(100) = 444,038$.


\begin{table}[hbt]
\centering
\begin{tikzpicture}
  \begin{scope}[scale=0.35, line width=1.4pt, shift={(-13, 3)}]
    \draw (3.36, 1.42) .. controls (3.36, 0.78) and (3.11, 0.13) .. 
          (2.55, 0.13) .. controls (1.89, 0.13) and (1.75, 0.97) .. (1.75, 1.75);
    \draw (1.75, 1.75) .. controls (1.75, 2.58) and (2.25, 3.36) .. (3.03, 3.36);
    \draw (3.69, 3.36) .. controls (4.12, 3.36) and (4.55, 3.36) .. (4.97, 3.36);
    \draw (4.97, 3.36) .. controls (5.75, 3.36) and (6.59, 3.22) .. 
          (6.59, 2.55) .. controls (6.59, 1.99) and (5.94, 1.75) .. (5.30, 1.75);
    \draw (4.64, 1.75) .. controls (4.22, 1.75) and (3.79, 1.75) .. (3.36, 1.75);
    \draw (3.36, 1.75) .. controls (2.93, 1.75) and (2.50, 1.75) .. (2.08, 1.75);
    \draw (1.42, 1.75) .. controls (0.13, 1.75) and (0.13, 3.95) .. 
          (0.13, 5.78) .. controls (0.13, 7.77) and (0.82, 9.82) .. 
          (2.55, 9.82) .. controls (3.72, 9.82) and (4.97, 9.54) .. (4.97, 8.53);
    \draw (4.97, 7.87) .. controls (4.97, 7.44) and (4.97, 7.02) .. (4.97, 6.59);
    \draw (4.97, 6.59) .. controls (4.97, 6.16) and (4.97, 5.73) .. (4.97, 5.30);
    \draw (4.97, 4.64) .. controls (4.97, 4.33) and (4.97, 4.01) .. (4.97, 3.69);
    \draw (4.97, 3.03) .. controls (4.97, 2.60) and (4.97, 2.17) .. (4.97, 1.75);
    \draw (4.97, 1.75) .. controls (4.97, 0.85) and (5.70, 0.13) .. 
          (6.59, 0.13) .. controls (8.20, 0.13) and (8.20, 2.77) .. (8.20, 4.97);
    \draw (8.20, 4.97) .. controls (8.20, 5.40) and (8.20, 5.83) .. (8.20, 6.26);
    \draw (8.20, 6.92) .. controls (8.20, 8.20) and (6.48, 8.20) .. (4.97, 8.20);
    \draw (4.97, 8.20) .. controls (3.43, 8.20) and (1.75, 7.92) .. 
          (1.75, 6.59) .. controls (1.75, 5.75) and (2.25, 4.97) .. (3.03, 4.97);
    \draw (3.69, 4.97) .. controls (4.12, 4.97) and (4.55, 4.97) .. (4.97, 4.97);
    \draw (4.97, 4.97) .. controls (5.94, 4.97) and (6.91, 4.97) .. (7.87, 4.97);
    \draw (8.53, 4.97) .. controls (9.17, 4.97) and (9.82, 5.22) .. 
          (9.82, 5.78) .. controls (9.82, 6.44) and (8.98, 6.59) .. (8.20, 6.59);
    \draw (8.20, 6.59) .. controls (7.24, 6.59) and (6.27, 6.59) .. (5.30, 6.59);
    \draw (4.64, 6.59) .. controls (3.87, 6.59) and (3.36, 5.81) .. (3.36, 4.97);
    \draw (3.36, 4.97) .. controls (3.36, 4.44) and (3.36, 3.90) .. (3.36, 3.36);
    \draw (3.36, 3.36) .. controls (3.36, 2.93) and (3.36, 2.50) .. (3.36, 2.08);

    \node at (4.5, -1) {\small $K11n34$};
  \end{scope}

  \node at (4, 2.75) {\tablefont
    \newcommand{\zp}{\phantom{0}}
    \newcommand{\zpc}{\phantom{0,}}
    \newcolumntype{i}{@{\hspace{2.5em}}c}
    \begin{tabular}{ccicicic}
      \toprule
      $g$ & $a_X(g)$ & $g$ & $a_X(g)$ & $g$ & $a_X(g)$ & $g$ & $a_X(g)$\\ \midrule
      1 & \zp 0 & \zp 7 & \zp 87 &  13 & \zpc 602  & 19  & 1,993 \\
      2 & \zp 6 & \zp 8 &    208 &  14 &   1,168  & 20  & 3,484 \\
      3 & \zp 9 & \zp 9 &    220 &  15 &   1,039  & 21  & 2,924 \\
      4 &    24 &    10 &    366 &  16 &   1,498  & 22  & 4,126 \\
      5 &    37 &    11 &    386 &  17 &   1,564  & 23  & 3,989 \\
      6 &    86 &    12 &    722 &  18 &   2,514  & 24  & 6,086 \\
      \bottomrule
    \end{tabular}
  };
\end{tikzpicture}
\caption{The first few values of $a_X(g)$ where $X$ is the exterior of
the Conway knot shown at left.}
\label{tab:conway}
\end{table}

While the sequence in Table~\ref{tab:conway} is a complete mystery to
us, if we broaden our perspective to include disconnected surfaces, we
get a relatively simple pattern that we can describe completely.
Specifically, for any $M$ define $b_M(n)$ to be the number of isotopy
classes of closed orientable essential surfaces $F$ in $M$ with
$\chi(F) = n$.  For the Conway exterior $X$, we show (see
Figure~\ref{fig:LWconway}):
\begin{equation}\label{eq:bX}
  b_X(-2n) = \frac{2}{3}n^3 + \frac{9}{4} n^2 + \frac{7}{3} n +
  \frac{7 + (-1)^n}{8}  \mtext{for all $n \geq 1$.} 
\end{equation}
The formula for $b_X$ would be a polynomial in $n$ were it not for the
final term which oscillates mod 2.  The first main result of this
paper, Theorem~\ref{thm.main} below, shows that the count $b_M$ always
has this kind of almost polynomial structure for a broad class of
\3-manifolds $M$.

\subsection{Main results}

We can encode a function $s \maps \N \to \Q$ by its \emph{generating
  function} $S(x) = \sum_{n=0}^\infty s(n) x^n$ in the formal power
series ring $\BQ[[x]]$.  We say this generating function is
\emph{short} when $S(x) = P(x)/Q(x)$ for polynomials $P$ and $Q$ in
$\Q[x]$ where $Q(x)$ is a product of cyclotomic polynomials. For
example, the function $s(n) = b_X(-2n)$ from (\ref{eq:bX}) above has a
short generating function, namely
\[
  S(x) = \frac{-x^5 + 3x^4 - 2x^3 + 2x^2 + 6x}{(x + 1)(x - 1)^4}
\]
Having a short generating function is equivalent to $s(n)$ being a
quasi-polynomial for all but finitely many values of $n$, see
Section~\ref{sec:genfunc}. Quasi-polynomials first arose in Ehrhart's
work on counting lattice points in polyhedra with rational vertices
\cite{Ehrhart} and have many applications to enumerative
combinatorics~\cite[Chapter 4]{Stanley}; curiously, they also appear
in quantum topology \cite{Ga:qdeg,GL:qholo,Ga:slope}.  We can now
state:
\begin{theorem}
  \label{thm.main}
  Suppose $M$ is a compact orientable irreducible \bdry irreducible
  atoroidal acylindrical \3-manifold that does not contain a closed
  nonorientable essential surface.  Let $b_M(n)$ be the number of
  isotopy classes of closed essential surfaces $F$ in $M$ with
  $\chi(F) = n$, and $B_M(x) = \sum_{n=1}^{\infty} b_M(-2n) x^n$ be
  the corresponding generating function. Then $B_M(x)$ is short.
\end{theorem}
Here, we can ensure that $M$ has no closed nonorientable essential
surfaces by requiring that $H_1(\partial M; \F_2) \to H_1(M; \F_2)$ is
onto, see Proposition~\ref{prop:nonorient}.  Thus,
Theorem~\ref{thm.main} applies to the exterior of any hyperbolic knot
in $S^3$.  We discuss possible extensions of Theorem~\ref{thm.main} to
nonorientable surfaces, as well as to surfaces with boundary, in
Section~\ref{sec:nonorient}.

All aspects of Theorem~\ref{thm.main} can be made algorithmic, both in
theory and in practice. The theoretical part is:
\begin{theorem}
  \label{thm.algol}
  There exists an algorithm that takes as input a triangulation $\cT$ of a
  manifold $M$ as in Theorem~\ref{thm.main} and computes
  $P(x), Q(x) \in \Q[x]$ such that $B_M(x) = P(x)/Q(x)$.  Moreover,
  there is an algorithm that given $n \in \N$ outputs a list of normal
  surfaces in $\cT$ uniquely representing all isotopy classes of
  essential surfaces with $\chi = -2n$.  Finally, there is an
  algorithm that given an essential normal surface $F$ with
  $\chi(F) = -2n$ finds the isotopic surface in the preceding
  list.
\end{theorem}
In Section~\ref{sec.ideal}, we refine Theorem~\ref{thm.algol} into a
practical algorithm that uses ideal triangulations and their special
properties.  Then in Section~\ref{sec.computations}, we compute
$B_M(x)$ for almost 60,000 examples.  It is natural to ask whether one
could permit nonorientable essential surfaces in
Theorem~\ref{thm.main}, as well as essential surfaces with boundary;
we outline some of the difficulties inherent in such extensions in
Section~\ref{sec:nonorient} below.
 
\subsection{Motivation and broader context}
\label{sec:motivation}
From Theorem~\ref{thm.main} and the discussion in
Section~\ref{sec:genfunc}, the sequence $b_M(-2n)$ grows at most
polynomially in $n$.  It is not always the case that $b_M(-2n)$ is
asymptotic to $c n^d$: we found an example where $b_M(-2n)$ is
$n/2 + 1$ for $n$ even and $0$ for $n$ odd.  However, by
Lemma~\ref{lem:quasismooth}, we get precise asymptotics if we smooth
the sequence by setting $\bbar_M(-2n) = \sum_{k = 1}^{n} b_M(-2k)$:
\begin{corollary}\label{cor:asymp}
  For each $M$ as in Theorem~\ref{thm.main}, either $b_M(-2n) = 0$ for
  all $n$ or there exists $d \in \N$ and $c > 0$ in $\Q$ such that 
  $\lim_{n \to \infty} \bbar_M(-2n)/n^d = c$.
\end{corollary}
We conjecture in Section~\ref{sec:ML} below that $d$ is the dimension
of the space $\ML_0(M)$ of measured laminations without boundary in
$M$, and $c$ is the volume of a certain subset of $\ML_0(M)$.

As $a_M(g) \leq b_M(-2g + 2)$ for each $g$, we have that
$a_M(g)$ also grows at most polynomially in $g$.  In stark contrast, if we
allow \emph{immersed} surfaces, then Kahn-Markovic
\cite{Kahn-Markovic:counting} showed that, for any closed hyperbolic
3-manifold $M$, the number of essential immersed surfaces of genus $g$
grows like $g^{2 g}$.

This distinction between counts of embedded versus immersed surfaces
parallels the following story a dimension down.  For a closed
hyperbolic surface $Y$ of genus $g$, Mirzakhani \cite{Mirzakhani2008}
showed that the number $s_Y(L)$ of embedded essential multicurves in
$Y$ of geodesic length at most $L$ satisfies $s_Y(L) \sim n(Y) L^{6g}$
for some $n(Y) > 0$; in contrast, the number $c_Y(L)$ of primitive
closed geodesics of length at most $L$ satisfies $c_Y(L) \sim e^L/L$,
see e.g.~\cite{Buser1992}.  In fact, Mirzakhani proved much more:
given an essential multicurve $\gamma$, the count $s_Y(L, \gamma)$ of
multicurves in the mapping class group orbit of $\gamma$ also
satisfies $s_Y(L, \gamma) \sim n_\gamma(Y) L^{6g}$ with
$n_\gamma(Y) > 0$. In particular, this gives asymptotics for the
counts of all \emph{connected} essential curves, analogous in our
setting to $a_M$ as opposed to $b_M$; we hint at how this connection
might be further developed in Section~\ref{sec:approaching_ag}. There
are also similarities between the setting of \cite{Mirzakhani2008} and
the measured lamination perspective on Theorem~\ref{thm.main} outlined
in Section \ref{sec:ML}.  The fact that we count surfaces by Euler
characteristic, which is discrete, rather than by a continuous notion
such as length or area, is what allows us get precise formulas for
$b_M$ as well as asymptotics.  (More directly analogous to the surface
case, one could try to count embedded essential surfaces in a closed
hyperbolic 3-manifold $M$ in terms of the area of their least area
representatives.  As such representatives satisfy
$\pi \abs{\chi(F)} \leq \Area(F) \leq 2 \pi \abs{\chi(F)}$ by
\cite[Lemma 6]{Hass1995}, it is not inconceivable that there are good
asymptotics here as well given Corollary~\ref{cor:asymp}.)

The algorithm of Section~\ref{sec.ideal} relies heavily on ideal
triangulations and their normal and almost-normal surfaces. Curiously,
normal surfaces are also used to construct recent topological quantum
invariants of 3-manifolds, specifically the 3D-index of Dimofte,
Gaiotto and Gukov~\cite{DGG1,DGG2}. The latter is a collection of
Laurent series with integer coefficients which are defined using an
ideal triangulation and depend only on the number of tetrahedra around
each edge of the triangulation, as encoded in the Neumann-Zagier
matrices. The 3D-index is a topological invariant of cusped hyperbolic
3-manifolds \cite{GHRS} that can be expressed as a generating series
of generalized normal surfaces in a 1-efficient
triangulation~\cite{GHHR}, a class of surfaces that includes both
normal and almost normal surfaces. It would be very interesting to
connect the topological invariants of Theorem~\ref{thm.main} with the
3D-index.

\subsection{The key ideas behind Theorem~\ref{thm.main}}
\label{sec:keyideas}

We first explain how the perspective of branched surfaces, especially
the work of Oertel \cite{Oertel:branched}, naturally relates the
sequence $b_M(-2n)$ to counting lattice points in an expanding family
of rational polyhedra; combined with Ehrhart's work \cite{Ehrhart} on
the latter topic, this discussion will make Theorem~\ref{thm.main}
very plausible. We then sketch how Tollefson \cite{To:isotopy}
reinterpreted and extended Oertel's branched surface picture in the
context of normal surface theory, and how this viewpoint allows us to
actually prove Theorem~\ref{thm.main}.  For ease of exposition, we
assume throughout that $M$ is closed and contains only
orientable surfaces by Proposition~\ref{prop:nonorient}.

A branched surface $\cB$ in a \3-manifold $M$ is the analog, one
dimension up, of a train track on a surface; see \cite{Floyd-Oertel,
  Oertel:branched} for definitions and general background.  A surface
$F$ is carried by $\cB$ if it is isotopic into a fibered neighborhood
$N(\cB)$ of $\cB$ so that it is transverse to the vertical interval
fibers.  Such an $F$ is determined by the nonnegative integer weights
it associates to the sectors of $\cB$, which are the components of
$\cB$ minus its singular locus.  Such weights correspond to a surface
if and only if they satisfy a system of homogenous linear equations
that are analogous to the switch conditions for a train track.  The
set of all nonnegative \emph{real} weights satisfying these equations
gives a finite-sided polyhedral cone $\ML(\cB)$, which corresponds to
measured laminations carried by $\cB$.  Here, each integer lattice
point in $\ML(\cB)$ corresponds to a surface carried by $\cB$.  As the
equations defining $\ML(\cB)$ have integer coefficients, each edge ray
of the cone $\ML(\cB)$ contains a lattice point.

For $M$ as in Theorem~\ref{thm.main}, by Theorem 4 of
\cite{Oertel:branched} there is a finite set $\cB_1,\ldots,\cB_n$ of
branched surfaces that together carry all essential surfaces in $M$
and also carry \emph{only} essential surfaces.  Moreover, two surfaces
carried by one $\cB_i$ are isotopic if and only if they correspond to
the same lattice point in $\ML(\cB_i)$.  Putting aside the important
issue of surfaces being carried by several of these branched surfaces,
here is how to count essential surfaces carried by a fixed $\cB_i$.
First, there is a linear function $\chibar \maps \ML(\cB_i) \to \R$
which on lattice points gives the Euler characteristic of the
corresponding surface.  Because $M$ is irreducible and atoroidal,
every essential surface has $\chi < 0$; as each edge ray of the
cone $\ML(\cB_i)$ contains a lattice point corresponding to an
essential surface, we conclude that $\chibar$ is proper and
nonpositive on $\ML(\cB_i)$.  Hence $P = \chibar^{\ -1}(-1)$ is a
compact polytope with, it turns out, rational vertices.  Thus, the
contribution to $b_M(-2n)$ of surfaces carried by $\cB_i$ is exactly
the number of lattice points in $2n \cdot P$, where the latter denotes
the dilation of $P$ by a factor of $2n$.  The foundational work of
Ehrhart \cite{Ehrhart} shows that this count of lattice points is
quasi-polynomial.

If no surface is carried by multiple $\cB_i$, the sketch just given
would essentially prove Theorem~\ref{thm.main} as sums of
quasi-polynomials are again quasi-polynomial.  However, there is no
avoiding this issue in general, and we deal with it by using the work
of Tollefson \cite{To:isotopy}, who built on \cite{Oertel:branched} to
provide a concrete description of isotopy classes of essential
surfaces in the context of normal surface theory.  If we fix a
triangulation $\cT$ of $M$, then every essential surface in $M$ can be
isotoped to be normal with respect to $\cT$; throughout, see
Section~\ref{sec:normal} for definitions and general background.
There can be many normal representatives of the same essential
surface, so to reduce this redundancy, Tollefson focuses on those that
are least weight in that they meet the 1-skeleton of $\cT$ in as few
points as possible.  We define a \emph{lw-surface} to be a normal
surface that is essential and least weight.  To prove
Theorem~\ref{thm.main}, we need to count such lw-surfaces modulo
isotopy in $M$.

Let $\ST$ be the normal surface solution space, which is a finite
rational polyhedral cone whose admissible integer points correspond to
normal surfaces in $\cT$, and let $\PT$ be its projectivization.  A
normal surface $F$ is \emph{carried} by a face $C$ of $\PT$ if the
projectivization of the lattice point corresponding to $F$ is in $C$.
An admissible face $C$ of $\PT$ is a \emph{lw-face} if every normal
surface it carries is a lw-surface.  While it is not obvious that any
lw-faces exist, Tollefson showed that every lw-surface is carried by a
lw-face. To make the parallel with the previous discussion explicit,
each lw-face $C$ has a corresponding branched surface $\cB_C$ which
carries, in the prior sense, exactly the surfaces carried by $C$ in
the current sense.  The collection of all lw-faces is a complex we
denote $\LW$; see Figure~\ref{fig:LWconway} for an example in the case
of a triangulation of the Conway knot exterior.

Tollefson shows moreover that every lw-surface is carried by a
lw-face that is \emph{complete}: if $F$ and $G$ are isotopic
lw-surfaces and $C$ carries $F$ then it also carries $G$. The
isotopies between lw-surfaces carried by the same complete lw-face can
be understood using a foliation of $C$ by affine subspaces parallel to
some fixed linear subspace $W_C$; roughly, surfaces $F$ and $G$
carried by $C$ are isotopic if their lattice points differ by an
element of $W_C$.  See Section~\ref{sec.thm.main} and especially
Theorem~\ref{thm.compiso} for details, including the key notion of
$\dep C$.  This translates the problem of counting essential surfaces
carried by a complete face to one of counting \emph{projections} of
lattice points in the cone over $C$ after we quotient out by $W_C$.
This is exactly the setting of recent work of Nguyen and Pak
\cite{N-Pak}, which we use in Section~\ref{sec:counts} to complete the
proof of Theorem~\ref{thm.main}.

\subsection{Making Theorem~\ref{thm.main} algorithmic}

Since Haken, normal surfaces have played a key role in the study of
algorithmic questions about \3-manifolds.  Despite this, Tollefson in
\cite{To:isotopy} did not give an algorithm for finding the lw-faces
of $\PT$ nor determining their properties such as completeness.
Section~\ref{sec:decision} here focuses on establishing
Theorem~\ref{thm.lwfaces}, which gives an algorithm for computing all
complete lw-faces.  One important tool for this is
Theorem~\ref{thm.path}, which shows that if $F$ and $G$ are isotopic
lw-surfaces then there is a sequence of isotopic lw-surfaces
$F = F_1, F_2, \ldots, F_{n-1}, F_n = G$ with each pair $(F_i, F_{i+1})$
disjoint and cobounding a product region.  Combined with results from
Section~\ref{sec.thm.main}, especially Theorem~\ref{thm.interior},
we can strengthen the arguments behind Theorem~\ref{thm.main} to
prove Theorem~\ref{thm.algol}.

\subsection{Ideal triangulations and almost normal surfaces}

When the \3-manifold $M$ has nonempty boundary, the proofs of
Theorems~\ref{thm.main} and \ref{thm.algol} use ideal triangulations
rather than finite ones (see Section~\ref{sec:ideal}). Our
computations were with $M$ where $\partial M$ is a single torus whose
interior admits a complete hyperbolic metric of finite-volume, and we
used ideal triangulations there as well, especially as they have
several advantages.  For example, they typically have fewer tetrahedra
than finite triangulations, which speeds up normal surface
computations.  Most importantly, when the ideal triangulation admits a
strict angle structure, Lackenby showed \cite{Lackenby:heeg} that the
number of connected normal surfaces of a fixed genus is finite and
described how they can be enumerated. In Section~\ref{sec.ideal}, we
explain how to exploit this to give a \emph{practical} version of the
algorithms in Theorem~\ref{thm.algol}.  Unlike the proof of
Theorem~\ref{thm.algol}, we make heavy use of almost normal surfaces,
including those with tubes, and in particular the process of
tightening (also called normalizing) an almost normal surface.

The usual method for testing if a normal surface $F$ in $M$ is
essential is to cut $M$ open along $F$, triangulate the result, and
then use normal surfaces to search for a compressing disk; a key
difficulty with this is that the triangulation of $M \setminus F$ is
usually much more complicated than the original one.  Here, we
introduce a completely new method for determining when $F$ is
essential that does not require cutting $M$ open but rather uses
almost normal surfaces in the original triangulation
(Section~\ref{sec:NTgraph}).

Our implementation of the algorithm in Section~\ref{sec.ideal} can be
found at \cite{PaperData} and makes heavy use of Regina
\cite{Regina}, SageMath \cite{SageMath}, and Normaliz \cite{Normaliz}.
It includes code for tightening almost normal surfaces, as well as
dealing with general normal surfaces with tubes, both of which have
explored extensively in theory but never before in practice.

\subsection{Computations and patterns}
\label{sec:comp_intro}

Sections~\ref{sec.computations} and \ref{sec.genuspat} detail our
experiments using the algorithm of Section~\ref{sec.ideal}.
In particular, we applied it to more than 59,000 manifolds, including
more than 4,300 where $\dim \LW > 0$.  We include overall statistics
about the complexes $\LW$, the generating functions $B_M(x)$, and the
sequences $a_M(g)$ in Tables~\ref{tab:overall}--\ref{tab:genmessy} and
\ref{tab:ag_common}--\ref{tab:lambert}, as well as detailed examples
of $\LW$ in Figures~\ref{fig:K15n51747}--\ref{fig:K13n1019}.  In
Section~\ref{sec:indep}, we give examples showing that, perhaps
surprisingly, neither of $B_M(x)$ and $a_M(g)$ determines the other.


For the more mysterious $a_M(g)$, while we are unable to find a
pattern in these sequences in many cases, there are some $M$ where we
conjecture relatively simple formulae for $a_M(g)$; see
Conjecture~\ref{conj:K13n586} and Table~\ref{tab:lambert}.  In
Conjecture~\ref{conj:asymp}, we posit the existence of general
asymptotics for (a smoothed version of) $a_M(g)$ based on the striking
plots in Figures~\ref{fig:asymp} and \ref{fig:asymp2345}, where we
computed $a_M(g)$ out to $g=200$ in many cases.
  
\subsection{The view from measured laminations}
\label{sec:ML}

For surfaces, a central tool for studying their topology, geometry,
and dynamics is measured laminations; for example, the space $\ML(F)$
of all measured laminations on a surface $F$ plays a key role in
\cite{Mirzakhani2008}.  In 3-dimensions, building on work of Morgan
and Shalen \cite{MorganShalen1984, MorganShalen1988}, independently
Hatcher \cite{Hatcher:laminations} and Oertel
\cite{Oertel:laminations} studied measured laminations on \3-manifolds
in detail, organizing them into a topological space $\ML(M)$.  Note
here an essential surface, with or without boundary, can be viewed as
a measured lamination, and the set of all essential surfaces nearly
injects into $\ML(M)$ (see page 6 of \cite{Hatcher:laminations} for
the caveat which involves the two nonorientable surfaces in a
semifibration) with its image being a discrete set of points.  While
for a surface $F$ of genus $g$ the space $\ML(F)$ is just homeomorphic
to $\R^{6g-6}$, for a general 3-manifold $M$ the space $\ML(M)$ can be
singular, being built from open strata each of which is a PL
manifold. The charts on the individual strata come from branched
surfaces; specifically, one uses the polyhedral cones $\ML(\cB_i)$
associated with certain essential branched surfaces $\cB_i$ as
sketched in Section~\ref{sec:keyideas}; see \cite{Hatcher:laminations,
  Oertel:laminations} for details.

Let $\ML_0(M)$ denote the subset of measured laminations that are
disjoint from $\partial M$.  The topological dimension of $\ML_0(M)$
is the maximum of $\dim(\ML(\cB))$ for the appropriate class of
essential branched surfaces $\cB$ without boundary.  Because of the
theory of Oertel \cite{Oertel:branched} that underlies
\cite{To:isotopy}, we are highly confident that:
\begin{conjecture}
  \label{conj:dimML}
  The dimension of $\ML_0$ is the maximum of $\dim C - \dim W_C + 1$
  where $C$ is an essential lw-face of $\PT$ and $W_C$ is defined in
  Theorem~\ref{thm.compiso}.
\end{conjecture}
If Conjecture~\ref{conj:dimML} holds, then in
Corollary~\ref{cor:asymp} where $\bbar_M(-2n) \sim c n^d$ one has
$d = \dim(\ML_0)$, thus giving an intrinsic characterization of that
exponent.  We further posit that the coefficient $c$ in these
asymptotics has the following natural interpretation.  As mentioned,
the PL structure on the strata of $\ML_0(M)$ comes from charts to
$\ML(\cB_i)$ for certain branched surfaces $\cB_i$; in particular, one
gets PL coordinate change maps between (possibly empty) subsets of
each pair $\ML(\cB_i)$ and $\ML(\cB_j)$, see Proposition~4.1 of
\cite{Hatcher:laminations}. These coordinate change maps must take
lattice points to lattice points, since these correspond to the
special measured laminations that come from essential surfaces.  Hence
the coordinate change maps should have derivatives that are in
$\GL{n}{\Z}$ and so are (unsigned) volume preserving. This would give
a well-defined measure (in the Lebesgue class) on each strata of
$\ML_0(M)$; this is a direct analog of Thurston's notion of volume on
$\ML(F)$ where $F$ is a surface, which is defined in terms of the
integral PL structure on $\ML(F)$ coming from train track charts.

Recall for any branched surface $\cB_i$, there is a linear map
$\chibar \maps \cB_i \to \R$ which gives the Euler characteristic of
the corresponding surface at each lattice point.  These should piece
together to give a PL map $\chibar \maps \ML_0(M) \to \R$.  In the
setting of Theorem~\ref{thm.main}, the subset
$P = \chibar^{\ -1}\big([-1, 0]\big)$ in $\ML_0(M)$ will be compact.  We
conjecture that the coefficient $c$ is precisely $\vol(P)$.

\subsection{Understanding counts by genus}
\label{sec:approaching_ag}

The key problem to overcome in understanding $a_M(g)$ is to
determine, for a complete lw-face $C$, which lattice points in
$\Ctil = \R_{\geq 0} \cdot C$ correspond to connected surfaces.  Agol,
Hass, and Thurston showed in \cite[\sec 4]{AgolHassThurston2006} how
counting the number of connected components of a normal surface can be
reframed as counting the number of orbits of a family of interval
isometries acting on $\{1, 2, \ldots, N\}$.  Such families of interval
isometries include both classical and non-classical interval exchange
transformations on surfaces \cite{Gadre2012}, but are considerably
more general.  Geometrically, a family of interval isometries can be
thought of as an interval $I$ of some length $L$ to which finitely
many bands of specified widths are attached, without any restriction
on how many bands are glued to any subinterval of $I$.  For normal
surfaces, the interval $I$ is basically an arbitrary concatenation of
the edges of the ambient triangulation $\cT$, and the bands correspond
to families of normal arcs in the corners of each face of $\cT^2$; see
Corollary 13 of \cite{AgolHassThurston2006} for details.  (For each
admissible face $C$ of $\PT$, one can also think about this in terms
of the associated branched surface $\cB_C$.) Thus, a general theory of
the number of orbits of the integer points of such a family of
isometries should allow one to develop a detailed picture for $a_M(g)$.

Currently, the best understood case is for a suitable train track
$\tau$ on a surface $F$, where Mirzakhani \cite{Mirzakhani2008} gives
asymptotics on the portion of integer points in $\ML(\tau)$ that
correspond to connected curves, see also \cite{Bell2019} for a
detailed discussion.  (Here, one uses the total weight of a point in
$\ML(\tau)$ as the ``length'' of the associated multicurve,
rather than Euler characteristic in the 3-dimensional setting.)  Even
for simple train tracks, it seems that the counts of connected curves
can be irregular in the sense of Section~\ref{sec:reg}, so there is
work to be done even in that setting.

\subsection{Acknowledgements} 
We thank Alexander Barvinok and Josephine Yu for discussions about
counting lattice points that were crucial to the proof of
Theorem~\ref{thm.main}, as well as Craig Hodgson for helpful
discussions on several related projects.  We also thank the referee
for their very careful reading of this paper and detailed
comments. Dunfield was partially supported by U.S. National Science
Foundation grant DMS-1811156, and Rubinstein partially supported by
Australian Research Council grant DP160104502.


\section{Background and conventions}
\label{sec.background}

\subsection{Numbers} We use $\N$ to denote the nonnegative integers,
i.e.~$\N = \{0, 1, 2, \dots\}$.

\subsection{Surfaces in 3-manifolds}
\label{sec:surfaces}

Throughout the rest of this paper, every \3-manifold $M$ will be
compact, orientable, irreducible (every embedded sphere bounds a ball)
and \bdry irreducible (every properly embedded disk bounds a ball with
some disk in $\partial M$).  Surfaces need not be orientable, but will
always be embedded in any ambient \3-manifold, and in particular be
compact. Moreover, a surface $F$ in a \3-manifold $M$ will be assumed
to be \emph{properly} embedded with $F \cap \partial M = \partial F$,
except for compressing disks and \bdry compressing disks which we
define next.  A \emph{compressing disk} for a surface $F$ in a
\3-manifold $M$ is a disk $D \subset M$ where $D \cap F = \partial D$
and $\partial D$ does not bound a disk in $F$.  An orientable surface
$F$ in $M$ is \emph{incompressible} when it has no compressing disks
and is neither a sphere nor a disk.  (A more general notion of
incompressibility allows certain spheres and disks, but none such
exist in an irreducible and \bdry irreducible manifold.)  Since $M$ is
\bdry irreducible, any parallel copy of a component of $\partial M$ is
incompressible.

A \emph{\bdry compressing disk} $D$ for a surface $F$ in $M$ is one
where $\partial D$ consists of an arc $\alpha$ in $F$ and an arc
$\beta$ in $\partial M$, the interior of $D$ is disjoint from
$F \cup \partial M$, and $\alpha$ does not bound a disk in $F$ with a
segment of $\partial F$.  An orientable surface $F$ in $M$ is
\emph{\bdry incompressible} when it has no \bdry compressing disks and
is not itself a disk. A surface $F$ in $M$ is \emph{\bdry parallel}
when every connected component of $F$ is isotopic, keeping
$\partial F$ fixed, into $\partial M$; when $\partial F = \emptyset$,
this is equivalent to $F$ being ambient isotopic to a union of
parallel copies of components of $\partial M$.

An orientable surface $F$ in $M$ is \emph{essential} when it is
incompressible, \bdry incompressible, and no connected component is
\bdry parallel.
A \3-manifold is \emph{atoroidal} when it does not contain an
essential torus (this is sometimes called \emph{geometrically
  atoroidal}).  Similarly, it is \emph{acylindrical} when it does not
contain an essential annulus (also called \emph{anannular}).

\subsection{Nonorientable surfaces}

For a nonorientable surface $F$ in $M$, we define it to be
incompressible, \bdry incompressible, or essential when the boundary
of a regular neighborhood of $F$ has that same property.  One
could instead apply the above definitions directly to nonorientable
surfaces, which give significantly weaker conditions in general.
Sources such as \cite{Floyd-Oertel, To:isotopy} use the terms
injective and \bdry injective for what we here call incompressible and
\bdry incompressible to distinguish the possible definitions in the
nonorientable case.  Some corner cases are worth mentioning.  First,
with our conventions, a connected surface $F$ in $M$ is incompressible
if and only if $\pi_1 F \to \pi_1 M$ is injective and $F$ is not a
sphere, a disk, or $\RP^2$.  Also, if $M$ is the twisted interval
bundle over a nonorientable closed surface $F$, then $F$ is
incompressible but not essential.

In our main results, we require that $M$ contain no closed
nonorientable essential surfaces, and the following proposition
provides an easily checkable sufficient condition for this to be the
case:

\begin{proposition}\label{prop:nonorient}
  Suppose $M$ is a compact orientable \3-manifold. Every closed
  embedded surface in $M$ is orientable if and only if $H_2(\partial M;
  \F_2) \to H_2(M; \F_2)$ is onto.
\end{proposition}
Thus a closed $M$ contains only orientable surfaces
if and only if $H_2(M; \F_2) = 0$.  Using
the long exact sequence of the pair, you can check that the homological
condition in Proposition~\ref{prop:nonorient} is equivalent to
$\dim H_1(M; \F_2) = \frac{1}{2} \dim H_1(\partial M; \F_2)$.
\begin{proof}[Proof of Proposition~\ref{prop:nonorient}]
  It suffices to consider the case when $M$ is connected.  First, note
  that any closed surface $F$ in $M$ gives a class in $H_2(M;\F_2)$.
  Moreover, any $c$ in $H_2(M;\F_2)$ can be represented by a closed
  surface $F$ that is connected (by adding tubes between components if
  needed) and nonempty (by adding a sphere bounding a ball if
  $c = 0$).  In the rest of this proof, all surfaces will be
  connected, nonempty, and embedded in $M$.

  As $M$ is orientable, any nonorientable surface $F$ is
  nonseparating.  Also, given a nonseparating orientable surface $F$
  we can build a nonorientable surface as follows: take an embedded
  arc $\alpha$ in $M$ that meets $F$ only at its endpoints and goes
  from one side of $F$ to the other; attaching a tube to $F$ along
  $\alpha$ gives the desired nonorientable surface.  Thus every closed
  surface $F$ in $M$ is orientable if and only if every closed surface
  is separating.  So we will prove that the homological hypotheses of
  the proposition are equivalent to every closed surface in $M$ being
  separating.  When $M$ is closed, the proposition is now immediate
  since a closed surface $F$ is $0$ in $H_2(M; \F_2)$ if and only if
  it is separating.

  To prove the proposition when $M$ has boundary, it suffices to show
  that the class $[F]$ of a closed surface $F$ is in the image of
  $H_2(\partial M; \F_2)$ if and only if $F$ is separating.  If $F$ is
  separating, then $F$ divides $M$ into two pieces $A$ and $B$ and we
  have $[F] = [A \cap \partial M] = [B \cap \partial M]$, so $[F]$
  comes from $H_2(\partial M; \F_2)$ as claimed.  If instead $F$ is
  nonseparating, let $\gamma$ be a loop disjoint from $\partial M$
  that meets $F$ in a single point; hence the homology intersection
  pairing $H_2(M; \F_2) \times H_1(M; \F_2) \to \F_2$ has
  $[F] \cdot [\gamma] = 1$.  As any $c \in H_2(\partial M; \F_2)$ has
  $c \cdot [\gamma] = 0$, it follows that $[F]$ does not come from
  $H_2(\partial M; \F_2)$.  So we have characterized which $F$ give
  classes coming from $H_2(\partial M; \F_2)$, completing the proof.
\end{proof}

\subsection{Triangulations}
\label{sec:tri}

A \emph{triangulation} of a compact 3-manifold is a cell complex made
from finitely many tetrahedra by gluing some of their \2-dimensional
faces in pairs via orientation-reversing affine maps so that the link
of every vertex is either a sphere or a disc.  (For such face gluings,
the link condition is equivalent to the complex being a 3-manifold,
see e.g.~\cite[Prop.~3.2.7]{Thurston1997}.) In particular, a
triangulation is not necessarily a simplicial complex, but rather what
is sometimes called a semi\hyp simplicial, pseudo\hyp simplicial, or
singular triangulation.

An \emph{ideal triangulation} of a compact 3-manifold with nonempty
boundary is a cell complex $\cT$ made out of finitely many tetrahedra
by gluing \emph{all} of their \2-dimensional faces in pairs as above
with no conditions on the vertex links.  Here, the manifold $M$ being
triangulated is not the underlying space of $\cT$ but rather the
subset of it gotten by removing a small regular neighborhood of each
vertex.  Put another way, the manifold $M$ is what you get by gluing
together \emph{truncated} tetrahedra in the corresponding pattern.
Hence $M$ will be a compact 3-manifold with nonempty boundary, and
$\cT \setminus \cT^0$ is homeomorphic to the interior of $M$, where
$\cT^i$ denotes the $i$-skeleton of $\cT$. See
e.g.~\cite{Tillmann2008} for more background on ideal triangulations.

We will work with both kinds of triangulations in this paper and will
sometimes refer to the first kind as \emph{finite triangulations} for
clarity.

\subsection{Normal surfaces}
\label{sec:normal}

Our conventions and notation for normal surfaces closely follow
\cite[\sec 2]{To:isotopy}, which the reader should consult for
additional details beyond the sketch we give here.  Throughout, we
consider a fixed triangulation $\T$ of a compact \3-manifold $M$,
which can be either finite or ideal.  However, in the ideal case, we
only consider closed normal surfaces, not the spun-normal ones of
\cite{Tillmann2008, Walsh2011}.  An \emph{elementary disk} $E$ in a
tetrahedron $\Delta$ is a disk meeting each face of $\partial \Delta$
in either a straight line or the empty set; note $\partial E$ is
determined by $E \cap \Delta^1$, and \cite[pg.~1089]{To:isotopy} gives
a convention so that the interior of $E$ is determined by
$E \cap \Delta^1$ as well.  A surface $F$ in $M$ is \emph{normal} when
it is in general position with the skeleta of $\T$ and meets each
tetrahedron of $\T$ in elementary disks.  A normal surface $F$ is
completely determined by $F \cap \T^1$.  A \emph{normal isotopy} of
$M$ is one that leaves every simplex in $\T$ invariant.  The normal
isotopy classes of elementary disks in a tetrahedron $\Delta$ are
called the \emph{disk types}, of which there are seven: three kinds of
triangles and four kinds of quadrilaterals (or quads for short).
Fixing an ordering of the $t$ tetrahedra in $\T$ and the seven disk
types, a normal surface $F$ gives a tuple $\vF \in \N^{7 t}$ by
counting the number of occurrences of each disk type; these are called
the \emph{normal coordinates} of $F$, or more precisely the
triangle-quad normal coordinates.  Note that the vector $\vF$
determines $F$ up to normal isotopy.

The coordinates of $\vF$ satisfy a system of homogenous linear
equations, called the \emph{matching equations} in \cite{To:isotopy},
one for each arc type in a face of $\T^2$.  In the vector space
$\R^{7t}$, the intersection of the solutions to the matching equations
with the positive orthant gives a polyhedral cone $\ST$
called the \emph{normal solution space}.  A vector $\xx \in \ST$ is
\emph{admissible} when for every tetrahedron of $\T$ there is at most one
quad coordinate of $\xx$ that is nonzero. The points in $\ST$
corresponding to normal surfaces are precisely the admissible integral
points.

A key property of a normal surface $F$ is its \emph{weight}
$\weight (F)$ which is the number of times it intersects $\cT^1$ and
can be viewed as its combinatorial area.  This notion of weight
extends to a linear function $\weight \maps \R^{7t} \to \R$ as
follows.  For an elementary disk $E_i$ corresponding to coordinate
$i$, each vertex of $E_i$ is incident on an edge of $\cT^1$; take
$c_i$ to be the sum of the reciprocals of the valences of those edges.
Defining  $\weight (\xx) = \sum_{i} c_i x_i$, we have $\weight(F) =
\weight(\vF)$ for every normal surface $F$.

The \emph{projective solution space} $\PT$ for $\T$ is abstractly the
quotient of $\ST \setminus \{0\}$ modulo positive scaling. It is
useful to concretely identify $\PT$ with a subset of $\ST$, and here
\cite{To:isotopy} uses the points of $\ST$ whose coordinates sum
to 1.  However, we instead use the convention that
$\PT = \setdef{\xx \in \ST}{\weight(\xx) = 1}$ as this simplifies the
statement of a key result of \cite{To:isotopy}.  We will use
$\vF^* = (1/\weight(\vF)) \vF$ to denote the projectivization of $\vF$ and
call it the \emph{projective normal class} of $F$.

The \emph{carrier} $C_F$ of a normal surface $F$ is the unique minimal
face of $\PT$ containing $\vF^*$.  The faces of $\ST$ and hence
$\PT$ correspond to having some of the defining inequalities $x_i \geq 0$
become equalities. Thus the carrier $C_F$ is the face of $\PT$ cut out
by the requirement that $x_i = 0$ whenever $F_i = 0$.

If normal surfaces $F$ and $G$ are \emph{compatible} in the sense that
they never have distinct quad types in a single tetrahedron, then they
have a natural ``cut and paste'' geometric sum that is also a normal
surface.  This new surface is called their \emph{normal sum} and
denoted $F + G$. Its normal coordinates are $\vF + \vG$ and in
particular the normal sum is determined up to normal isotopy by the
normal isotopy classes of $F$ and $G$, even though $F \cap G$ can
change under normal isotopy of the surfaces individually.

\subsection{Short generating functions and quasi-polynomials}
\label{sec:genfunc}

Throughout this subsection, see Chapter 4 of \cite{Stanley} for details
and further background.  We can encode a function $s \maps \N \to \Q$
by its \emph{generating function} $S(x) = \sum_{n=0}^\infty s(n) x^n$
in the ring $\BQ[[x]]$ of formal power series.  This generating
function is \emph{short} when $S(x) = P(x)/Q(x)$ for polynomials $P$
and $Q$ in $\Q[x]$ where $Q$ is a product of cyclotomic
polynomials.  Equivalently, the generating function is short if and
only if
\[
  S(x) = \sum_{i=1}^k \frac{c_ix^{a_i}}{(1-x^{b_i})^{d_i}}
  \mtext{for some $c_i \in \Q$ and $a_i, b_i, d_i \in \N$.}
\]
If $s$ has a short generating function $S = P/Q$ where further
$\deg P < \deg Q$, then we say that $s$ is a
\emph{quasi-polynomial}.  Equivalently, a function $s \maps \N \to \Q$
is quasi-polynomial if and only if there exists $L \in \N$ and
polynomials $f_0, f_1,\ldots,f_{L - 1} \in \Q[x]$ such that
$s(n) = f_k(n)$ if $n \equiv k \bmod L$, see Proposition~4.4.1 of
\cite{Stanley}.  When $s$ has a short generating function, it
is equal to a fixed quasi-polynomial except for finitely many inputs
\cite[Proposition~4.2.2]{Stanley}.

We will be interested exclusively in $s$ where $s(n) \in \N$ for all
$n$. When such an $s$ has a short generating function
$S(x) = P(x)/Q(x)$, where $Q \in \Z[x]$ is a product of cyclotomic
polynomials, then $P$ must also be in $\Z[x]$; this is because
$P(x) = S(x) Q(x)$ as elements of $\Q[[x]]$ and $S(x) Q(x)$ is a
product of elements in $\Z[[x]]$.

We end this section with the lemma that gives
Corollary~\ref{cor:asymp} from Theorem~\ref{thm.main}:
\begin{lemma}
  \label{lem:quasismooth}
  Suppose $s \maps \N \to \Q$ with all $s(n) \geq 0$ has a short
  generating function and consider $\sbar(n) = \sum_{k = 0}^n s(k)$.
  Then either $s(n) = 0$ for all large $n$ or there exists $d \in \N$
  and $c > 0$ in $\Q$ such that $\sbar(n) \sim c n^d$.
\end{lemma}

\begin{proof} 
  Since we only care about asymptotics, assume that $s$ is a
  quasi-polynomial with $f_0, f_1, \ldots, f_{L-1} \in \Q[x]$ where
  $s(n) = f_\ell(n)$ if $n \equiv \ell \bmod L$.  Assume some
  $f_\ell \neq 0$ as otherwise we are done.  Set
  $e = \max(\deg f_\ell)$, which is at least 0, and let $c_\ell$ be
  the coefficient on $x^e$ in $f_\ell$, so that
  $f_\ell(n) = c_\ell n^e + O(n^{e - 1})$.  Then as $s(n) \geq 0$ for
  all $n$ it follows that $c_\ell > 0$ if $\deg f_\ell = e$ and
  otherwise $c_\ell = 0$; in particular, all $c_\ell \geq 0$ and
  $\sum_\ell c_\ell > 0$. Separating the sum in $\sbar(n)$ into
  congruence classes modulo $L$, we write
  \begin{equation}
    \label{eq.modL}
    \sbar(n)=\sum_{\ell=0}^{L-1} \sbar^{(\ell)}(n) \mtext{where}
    \sbar^{(\ell)}(n)=\sum_{\substack{j=0\\ j \equiv \ell \bmod L}}^n s(j) =
    \sum_{k=0}^{\lfloor{(n-\ell)/L\rfloor}} f_\ell(\ell + L k)
  \end{equation}
  Using that $(\ell + L k)^e$ is a polynomial in $k$ with leading
  term $L^e k^e$, we get
  $
  f_\ell(\ell + L k) = c_\ell (\ell + L k)^e + O(n^{e - 1}) = c_\ell
  L^e k^e + O(n^{e - 1})$ where $n = \ell + L k$.  Thus
  \begin{align*}
    \sbar^{(\ell)}(n)
    &= 
      \sum_{k=0}^{\lfloor{(n-\ell)/L\rfloor}} \left( c_\ell L^e k^e + O(n^{e-1}) \right)
    = c_\ell L^e \left( \sum_{k=0}^{\lfloor{(n-\ell)/L\rfloor}} k^e \right)
      + O(n^e) 
    \\ &=  \frac{c_\ell L^e}{e + 1}
         \left\lfloor \frac{n - \ell}{L} \right\rfloor^{e + 1} + O(n^e)
       = \frac{c_\ell}{(e + 1)L} n^{e + 1} + O(n^e)
  \end{align*}
  where we have used
  $\sum_{k = 0}^m k^e = \frac{m^{e + 1}}{e + 1} + O(m^e)$. Set
  $d = e + 1$ and $c = \frac{1}{d L} \sum_\ell c_\ell > 0$, and it now
  follows from (\ref{eq.modL}) that
  $\sbar(n) \sim c n^d$ as required.
\end{proof}

\section{Isotopy classes of essential normal surfaces}
\label{sec.thm.main}

In this section, we discuss and refine Tollefson's work on isotopy
classes of incompressible surfaces from the point of view of normal
surface theory.  In particular, this allows us to build a bijection
between isotopy classes of such surfaces and certain equivalence
classes of lattice points in a collection of rational cones.  We will
use this framework to prove Theorem~\ref{thm.main} in
Section~\ref{sec:counts}.  We begin by explaining some key facts from
Tollefson \cite{To:isotopy} using the notation that we reviewed in
Section~\ref{sec:normal}. Throughout, we consider a compact orientable
irreducible $\partial$-irreducible 3-manifold $M$ equipped with a
fixed finite triangulation $\cT$.

A \emph{lw-surface} is a compact orientable incompressible \bdry
incompressible normal surface that is least weight among all normal
surfaces in its isotopy class; such surfaces play a key role in
\cite{To:isotopy}. (The term lw-surface is not actually used in
\cite{To:isotopy} but makes its results easier to state.)  A face $C$
of $\cP_\cT$ is a \emph{lw-face} when every orientable normal surface
carried by $C$ is a lw-surface.  We use $\LW$ to denote the set of all
lw-faces of $\PT$. Clearly, $\LW$ is a subcomplex of $\PT$. A lw-face
$C$ is \emph{complete} if whenever it carries an orientable normal
surface $F$ it also carries every lw-surface isotopic to $F$.  A key
fact for us is:
\begin{theorem}[{\cite[Theorem~4.5]{To:isotopy}}]\label{thm.toll}
  Every lw-surface is carried by a complete lw-face.  In
  particular, any lw-face is contained in some complete lw-face.
\end{theorem}
On a complete lw-face, Tollefson characterizes the various possible
forms for isotopy relations among the surfaces that it carries.  As we
will describe, these have to be relatively simple on the interior
$\intC$ of $C$, but proper faces of $C$ can have different
isotopy relations. Tollefson introduces the notion of a
\emph{PIC-partition} to encode all of these.  We will not work with
PIC-partitions directly, but reframe the underlying structure in a way
more suited for the proof of Theorem~\ref{thm.main}.  To give our structure
theorem, we first need some definitions and a useful characterization
of when a face of $\LW$ is complete.

If $F$ is an orientable surface in $M$ and $m$ a positive integer, a
disjoint union of $m$ parallel copies of $F$ is called a
\emph{multiple} of $F$ and denoted $m F$. When $F$ is normal, we
always take $m F$ to be a normal surface whose normal coordinates are
$m \vF$.  To mirror what happens algebraically in the normal case, for
a nonorientable surface $F$ one defines $2 F$ as the boundary $G$ of a
regular neighborhood of $F$ and then $m F$ as either $\frac{m}{2}G$ or
$F \cup \frac{m-1}{2}G$ depending on the parity of $m$. Surfaces $F$
and $G$ are \emph{projectively isotopic} when they have multiples that
are isotopic.  We say two normal surfaces are \emph{projectively
  normally isotopic} if they have multiples that are normally
isotopic.  Note here that the admissible rational points of $\PT$
correspond exactly to projective normal isotopy classes of normal
surfaces.

\begin{remark}
  Our definitions of least-weight and completeness for a face $C$
  differ from those in \cite{To:isotopy} in that we only look at
  \emph{orientable} normal surfaces $F$ carried by $C$ whereas
  \cite{To:isotopy} allows nonorientable $F$.  However, it is easy to
  see our definitions are equivalent to the originals.  For example,
  if $C$ is a lw-face with our definition and $F$ is a nonorientable
  surface carried by $C$, then $2 F$ is a lw-surface and hence $F$
  itself is incompressible and \bdry incompressible.  Moreover, if $G$
  is any normal surface isotopic to $F$ then $2 G$ is isotopic to
  $2 F$ and hence $\weight(2 G) \geq \weight(2 F)$ which implies
  $\weight(G) \geq \weight(F)$; thus $F$ is least weight among all
  normal surfaces in its isotopy class.  The equivalence of the two
  definitions of completeness is similar, using that if $C$ carries
  $2 G$ then it carries $G$.
\end{remark}

Important for us throughout this paper is that whether a face is
(complete) least-weight is determined by any one surface carried by
its interior:
\begin{theorem}\label{thm.interior}
  Suppose $F$ is an orientable normal surface carried by the interior
  of a face $C$ of $\PT$. Then the following are equivalent:
  \begin{enumerate}[nosep]
  \item \label{item:Clw}
    $C$ is a lw-face.
  \item \label{item:Flw}
    $F$ is a lw-surface.
  \item \label{item:compslw}
    Every connected component of $F$ is a lw-surface.
  \end{enumerate}
  If $C$ is a lw-face, the following are equivalent:
  \begin{enumerate}[nosep, resume]
  \item \label{item:Ccomplete}
    $C$ is complete.
  \item \label{item:CisoF}
    $C$ carries every lw-surface isotopic to $F$.
  \item \label{item:Cisocomps}
    $C$ carries every lw-surface isotopic to a
    connected component of $F$.
  \end{enumerate}
\end{theorem}
We will prove Theorem~\ref{thm.interior} below in Section~\ref{sec:lwdetails}.

\subsection{Dependent faces}\label{sec:dep}
A face $D$ of a lw-face $C$ is \emph{$C$-dependent} if there exists a
surface carried by $D$ that is projectively isotopic to one carried by
$\intC$; otherwise, the face $D$ is \emph{$C$-independent}.  The
collection of $C$-independent faces of $C$ clearly forms a subcomplex
$\cD$ of $\partial C$ and we define $\dep{C}$ to be
$C \setminus \bigcup_{D \in \cD} D$.  Note that if $D$ is a
$C$-dependent face of $C$, then $\intD \subset \dep{C}$ since if any
$\xx \in \intD$ was in a $C$-independent face $E$ then $D$ would be a
face of $E$, contradicting that $D$ is $C$-dependent. As any point of
$C$ is in the interior of some face, we see that $\dep{C}$ is also
the union of $\intD$ over all $C$-dependent faces $D$ of $C$.

Tollefson completely
characterized the isotopy relations among the surfaces carried by each
$\dep{C}$.  We rework this as:
\begin{theorem}\label{thm.compiso}
  For each face $C$ of $\LW$ there is a rational linear subspace $W_C$
  such that the following holds.  Any two surfaces $F$ and $G$ carried
  by $\dep{C}$ are projectively isotopic if and only if
  $\vF^* - \vG^*$ is in $W_C$.  Moreover, if $F$ and $G$ are
  orientable then they are isotopic if and only if $\vF - \vG$ is in
  $W_C$.  Also,
  \begin{equation}\label{eq:altdep}
    \dep{C} = \setdef{\xx \in C}{\mbox{$\xx + W_C$ meets $\intC$}}
  \end{equation}
  so that in particular any $F$ carried by $\dep{C}$ is projectively
  isotopic to one carried by $\intC$.  Finally, the subspace $W_C$ is
  contained in $\ker(\weight)$ and given any $F$ carried by
  $\intC$ there exist surfaces $F_1, \ldots, F_k$ projectively
  isotopic to $F$ and carried by $\intC$ such that the
  $\vF^* - \vF_i^*$ span $W_C$.
\end{theorem}
The example in Section~\ref{sec:isoexs} may help you understand the
statement of Theorem~\ref{thm.compiso}.

For the practical algorithms in Section~\ref{sec.ideal}, we will need
the following additional properties of $W_C$:
\begin{corollary}\label{cor:WC} 
  Suppose $C$ is a face of $\LW$.  If surfaces $F$ and $G$ carried by
  $C$ are projectively isotopic, then $\vF^* - \vG^* \in W_C$.  Also,
  if $D$ is a face of $C$, then $W_D \subset W_C$.  Finally, if $F$ is
  any orientable surface carried by $\intC$, then $W_C$ is spanned by
  all $\vG - \vH$ where $G$ is a connected component of $F$ and
  $H$ is isotopic to $G$ and carried by $C$.
\end{corollary}
Combined with Theorem~\ref{thm.compiso}, the next result will be key
to proving Theorem~\ref{thm.main}:
\begin{theorem}\label{thm:disjointdep}
  The complex $\LW$ is the disjoint union of the $\dep{C}$ as $C$
  ranges over the complete lw-faces of $\PT$.  Moreover, if $C$ and
  $C'$ are distinct complete lw-faces, then no surface carried by
  $\dep{C}$ is projectively isotopic to one carried by
  $\dep{C'}$. Consequently, for an orientable incompressible \bdry
  incompressible surface $F$, there is a unique complete lw-face $C$
  such that $\dep{C}$ carries a surface (non-projectively) isotopic to
  $F$.
\end{theorem}

\begin{remark}
  We defined $\PT = \setdef{\xx \in \ST}{\weight(\xx) = 1}$ rather
  than $\PT' = \setdef{\xx \in \ST}{\sum x_i= 1}$ in order to state
  Theorem~\ref{thm.compiso} in the above form.  Tollefson uses $\PT'$,
  and ends up with a partition of $\dep{C}$ along a family of
  typically \emph{nonparallel} affine subspaces whereas our partition
  is along \emph{parallel} affine subspaces.  While both $\PT$ and
  $\PT'$ are projectivizations of $\ST$, the map that identifies them
  is not affine but rather projective and so this is not a
  contradiction.
\end{remark}

\subsection{Complete lw-faces in detail}
\label{sec:lwdetails}

We begin with the proof of Theorem~\ref{thm.interior} as it is needed
to prove Theorem~\ref{thm:disjointdep}.

\begin{proof}[Proof of Theorem~\ref{thm.interior}]
  First, recall that faces of $\PT$ are defined by setting a subset of
  the normal coordinates to $0$.  Consequently, a normal surface $G$ is
  carried by $C$ if and only if every connected component of it is
  carried by $C$.  More generally, if $K$ and $L$ are compatible
  normal surfaces, then $C$ carries $K + L$ if and only if it carries
  $K$ and $L$ individually.
  
  A normal surface $F$ is carried by $\intC$ if and only if the
  carrier of $F$ is equal to $C$; the equivalence of (\ref{item:Clw})
  and (\ref{item:Flw}) is thus Theorem 4.2 of \cite{To:isotopy}. From
  the definition it is clear that (\ref{item:compslw}) implies
  (\ref{item:Flw}), so to complete the proof of the first part of the
  theorem we will show (\ref{item:Clw}) implies (\ref{item:compslw}).
  This holds because if $F'$ is a component of $F$ then, as noted above,
  $C$ carries $F'$ and thus $F'$ is a lw-surface as $C$ is a lw-face.

  For the second part, by definition (\ref{item:Ccomplete}) implies
  (\ref{item:CisoF}), and (\ref{item:Ccomplete}) implies
  (\ref{item:Cisocomps}) since every component of $F$ is also carried
  by $C$.  Since $C$ carries a surface if and only if it carries all
  of its components, we see that (\ref{item:Cisocomps}) easily gives
  (\ref{item:CisoF}).  So it remains to prove (\ref{item:CisoF})
  implies (\ref{item:Ccomplete}).  So suppose $C$ is a least-weight
  face of $\PT$ such that every lw-surface projectively isotopic to
  $F$ is carried by $C$.  We must show that $C$ is complete, so
  suppose $K$ is a lw-surface carried by $C$ and $L$ is a lw-surface
  isotopic to $K$.  As $\vF^* \in \intC$, we can pick a
  lw-surface $E$ with $\vE^* \in \intC$ and $\vF^*$ in the interior of
  the line segment joining $\vK^*$ to $\vE^*$.  Then there are
  positive integers $\{m, a, b\}$ with $m F = a K + b E$.  Applying
  Corollary~4.3 of \cite{To:isotopy} with $G = a K$, $G' = a L$, and
  $H = H' = b E$, we conclude that $a L$ and $b E$ are compatible and
  that $a L + b E$ is isotopic to $m F$.  By hypothesis, as
  $a L + b E$ is projectively isotopic to $F$, it is carried by $C$.
  Then $a L$ is carried by $C$ and hence $L$ is carried by $C$ as well
  as $a > 0$.  Hence $C$ is complete as claimed.
\end{proof}

We now turn to the proof of Theorem~\ref{thm.compiso} for which we
will need:
\begin{lemma}\label{lem:polygeom}
  Suppose $C$ is a compact convex polyhedron in $\R^n$ and $W$ a
  subspace of $\R^n$.  Set $\dep{C, W} = \setdef{\xx \in
    C}{\mbox{$\xx + W$ meets $\intC$}}$.  If $D$ is a face of $C$,
  then the intersection $D \cap \dep{C, W}$ is either empty or
  contains $\intD$.
\end{lemma}
\begin{proof}
  Passing to a subspace if necessary, we assume that $\dim C = n$ and
  hence $\intC$ is open in $\R^n$.  There are finitely many nonzero
  linear functionals $\ell_i$ on $\R^n$, say indexed by a set $I$, and
  $\alpha_i \in \R$, such that $C = \setdef{\xx \in \R^n}{%
    \mbox{$\ell_i(\xx) \geq \alpha_i$ for all $i \in I$}}$. Then
  $\intC = \setdef{\xx \in \R^n}{%
    \mbox{$\ell_i(\xx) > \alpha_i$ for all $i \in I$}}$.  For a face
  $D$ of $C$, define $I_D$ to be the indices in $I$ where
  $\ell_i(\xx) = \alpha_i$ on all of $D$.

  Now assume $D \cap \dep{C, W}$ is nonempty, and pick $\xx \in D$ and
  $\ww \in W$ with $\xx + \ww$ in $\intC$.  For any $i \in I$, we have
  $\ell_i(\xx + \ww) > \alpha_i$, which for those $i \in I_D$ implies
  $\ell_i(\ww) > 0$ since $\ell_i(\xx) = \alpha_i$.  Given $\yy $ in
  $\intD = \setdef{ \xx \in D}{\mbox{$\ell_i(\xx) > \alpha_i$ for all
      $i \notin I_D$}}$, we need to show that it is in $\dep{C,
    W}$. For $\epsilon > 0$, consider $\vv = \yy + \epsilon \ww$.  For
  $i \notin I_D$, we have
  $\ell_i(\vv) = \ell_i(\yy) + \epsilon \ell_i(\ww)$; since
  $\ell_i(\yy) > \alpha_i$, we can thus make $\ell_i(\vv) > \alpha_i$
  as well by choosing $\epsilon$ small enough. On the other hand, for
  $i \in I_D$ we have
  $\ell_i(\vv) = \alpha_i + \epsilon \ell_i(\ww) > \alpha_i$ for any
  positive $\epsilon$ as $\ell_i(\ww) > 0$ for such $i$.  Thus $\vv$
  is in $\intC$ for small $\epsilon$ and so $\yy \in \dep{C, W}$ as
  needed.
\end{proof}

\begin{proof}[Proof of Theorem~\ref{thm.compiso}]
  This result is essentially a reframing of Theorem 5.5 of
  \cite{To:isotopy} on the existence of a PIC-partition for $C$, but
  to see this one must use a number of details from the proof of that
  theorem.  Hence we will simply prove Theorem~\ref{thm.compiso}
  directly relying on results earlier in that paper.  Suppose $F$ is
  any lw-surface carried by $\intC$.  Let $V_F$ be the subspace of
  $\R^{7t}$ spanned by all $\vG$ where $G$ is carried by $C$ and
  projectively isotopic to $F$.  Using that $V_F$ is
  finite-dimensional, we can find mutually isotopic lw-surfaces
  $F_1, \ldots, F_k$ carried by $C$, each projectively isotopic to $F$
  with $\vF_1^* = \vF^*$, such that $\{\vF_1, \ldots, \vF_k\}$ is a
  basis for $V_F$.  Then Theorem~5.3 of \cite{To:isotopy} implies that
  every normal surface carried by $A_F = C \cap V_F$ is projectively
  isotopic to $F$.

  Consider the affine subspace
  $X = \setdef{\xx \in \R^{7t}}{\weight(\xx) = 1}$ and note
  $\PT = \ST \cap X$ where $\ST$ is the normal solution space.  Define
  $X_F = X \cap V_F$ which is also the smallest affine subspace
  containing $\{\vF_1^*, \ldots, \vF_k^*\}$.  Note here that $A_F$ is
  also $C \cap X_F$.  If $\vG$ is another orientable normal surface
  carried by $\intC$, then either $A_F = A_G$ or
  $A_F \cap A_G = \emptyset$ depending on whether or not $F$ and $G$
  are projectively isotopic.  When $A_F$ and $A_G$ are disjoint, we
  claim that $X_F$ and $X_G$ are still parallel; formally, the tangent
  space to an affine subspace $Y \subset \R^{7t}$ is
  $TY = \setdef{\yy_1 - \yy_2}{\yy_1, \yy_2 \in Y}$ and we will show
  $TX_F = TX_G$.

  As $\vG^* \in \intC$, we can find a lw-surface $H$ with $\vH^* \in
  \intC$ and $\vG^*$ on the interior of the line segment between
  $\vF_1^*$ and $\vH^*$.  Hence there are positive integers $\{m, a,
  b\}$ such that $m G = a F_1 + b H$.  Set $G_i = a F_i + b H$.  By
  Corollary~4.3 of \cite{To:isotopy}, all the $G_i$ are isotopic to
  $G_1 = m G$ and hence lie in $V_G$. The $F_i$ are isotopic
  lw-surfaces and so have the same weight, and consequently so do the 
  $G_i$ and hence
  \[
    \vF_i^* - \vF_j^*= \frac{\weight(G_1)}{a \cdot \weight(F_1)}
    \left(\vG_i^* - \vG_j^*\right) \mtext{for all $i, j$.}
  \]
  In particular, we have $T X_F \subset T X_G$.  Reversing the roles
  of $F$ and $G$ shows $T X_G = T X_F$ as claimed.

  Now set $W_C = T X_F$ for any orientable $F$ carried by $\intC$.
  Note that $W_C$ is spanned by the
  $\vF_i - \vF_j = \weight(F_1)\big(\vF_i^* - \vF_j^*\big)$ from
  above, and hence by the $\vF^* - \vF_i^*$ since $\vF_1^* = \vF^*$.
  As $F$ is arbitrary, this verifies the claims in the last
  sentence of the statement of the theorem since in addition each
  $\vF_i - \vF_j$ is in $\ker(\weight)$.

  We extend our notion of $A_F$ to all
  $\yy \in \intC$ by setting $A_{\yy} = (\yy + W_C) \cap C$.  Let
  $\Atil = \bigcup_{\yy \in \intC} A_{\yy}$.  We will show:
  \begin{claim}\label{claim:key}
    $\Atil = \dep{C}$.
  \end{claim}
  Before proving the claim, let us show that Theorem~\ref{thm.compiso}
  follows from it.  First, the equation (\ref{eq:altdep}) holds since
  $\Atil$ is also $\setdef{\xx \in C}{\mbox{$\xx + W_C$ meets
      $\intC$}}$.  Second, the claim that surfaces $F$ and $G$ carried
  by $\Atil$ are projectively isotopic if and only if
  $\vF^* - \vG^* \in W_C$ follows because $\Atil$ is partitioned by
  the $A_{\yy}$ which for rational $\yy$ correspond exactly to
  projective isotopy classes of surfaces carried by $C$.  Finally, if
  $F$ and $G$ are orientable surfaces carried by $\dep{C}$, we need to
  show they are isotopic if and only if $\vF - \vG \in W_C$.  If they
  are isotopic, they must have the same weight and so $\vF - \vG$ is a
  multiple of $\vF^* - \vG^*$, and the latter must be in $W_C$ as $F$
  and $G$ are projectively isotopic.  Conversely, if
  $\vF - \vG \in W_C$, the surfaces are projectively isotopic and as
  $W_C \subset \ker(\weight)$ it follows $\weight(F) = \weight(G)$.
  As $F$ and $G$ are orientable, least weight, of the same weight, and
  projectively isotopic, they are actually isotopic as needed.

  To prove Claim~\ref{claim:key}, note both sets contain $\intC$, so
  for each face $D$ of $\partial C$ we will check that both sets agree
  on $\intD$.  By Lemma~\ref{lem:polygeom}, either $\Atil \cap D$ is
  empty or it contains $\intD$.  If the former, then no surface
  carried by $D$ can be projectively isotopic to one carried by
  $\intC$, and so $D$ is $C$-independent and hence
  $\dep{C} \cap D = \emptyset$ as well.  If the latter, then any
  rational point in $\intD$ gives a surface projectively isotopic to
  one carried by $\intC$; hence, $D$ is $C$-dependent. As noted in
  Section~\ref{sec:dep}, this implies $\intD \subset \dep{C}$. This
  proves the claim and hence the theorem.
\end{proof}

\begin{proof}[Proof of Corollary~\ref{cor:WC}] First, suppose $F$ and
  $G$ are projectively isotopic and carried by $C$.  Since this does
  not change $\vF^* - \vG^*$, we will assume $F$ and $G$ are
  orientable and actually isotopic.  Let $H$ be an orientable surface
  carried by $\intC$.  Then by Corollary~4.3 of \cite{To:isotopy}, we
  have $H + F$ and $H + G$ are isotopic, and hence by
  Theorem~\ref{thm.compiso} above we have
  $\overrightarrow{(H + F)} - \overrightarrow{(H + G)} = \vF - \vG =
  \weight(F)(\vF^* - \vG^*)$ is in $W_C$ as needed.

  Second, if $D$ is a face of $C$, then by the last part of
  Theorem~\ref{thm.compiso} we have $W_D$ is spanned by certain
  $\vF^* - \vF_i^*$ where $F$ and $F_i$ are carried by $D$.  By what
  we just showed, all of these are in $W_C$ as well, proving
  $W_D \subset W_C$.

  Finally, fix an orientable surface $F$ carried by $\intC$ and define
  $Z$ to be the span of all $\vG - \vH$ where $G$ is a connected
  component of $F$ and $H$ is isotopic to $G$ and carried by $C$.  We
  need to show $W_C = Z$.  By the first part of this corollary, we
  know $Z \subset W_C$. From Theorem~\ref{thm.compiso}, there are
  surfaces $F_1, \ldots, F_k$ projectively isotopic to $F$ where the
  $\vF^* - \vF_i^*$ span $W_C$.  We can moreover arrange that each
  $F_i$ is isotopic to $F$ so that the $\vF - \vF_i$ span $W_C$.  To
  see $\vF - \vF_i$ is in $Z$, let $G_1,\ldots,G_n$ be the connected
  components of $F$.  Under an isotopy between $F$ and $F_i$, let
  $G_j'$ be the connected component of $F_i$ corresponding to $G_j$.
  Then $\vF - \vF_i = \sum (\vG_j - \vG'_j)$ which is in $Z$, giving
  $W_C = Z$ and completing the proof of the corollary.
\end{proof}

Turning now to the proof of Theorem~\ref{thm:disjointdep}, we begin
with a lemma:

\begin{lemma}\label{lem:maxCindep}
  Suppose $C$ is a complete lw-face. A maximal $C$-independent face
  $D$ of $C$ is also complete.
\end{lemma}
\begin{proof}
  Pick a lw-surface $F$ carried by $\intD$.  By
  Theorem~\ref{thm.interior}, it suffices to show that given a
  lw-surface $G$ isotopic to $F$ then $G$ is carried by $D$.  By
  completeness of $C$, we know $G$ is carried by $C$. By Theorem~5.3
  of \cite{To:isotopy}, every normal surface carried by the segment
  $L = [\vF^*, \vG^*]$ in $C$ is projectively isotopic to $F$.  Let
  $E$ be the minimal face of $C$ containing $L$; since $L$ is just a
  segment, it meets $\interior{E}$.  We cannot have $E$ be $C$ as then
  $F$ is projectively isotopic to some surface in $\intC$, violating
  that $D$ is $C$-independent.  For the same reason, the face $E$
  cannot be $C$-dependent as then $\intE \subset \dep{C}$ and hence by
  Theorem~\ref{thm.compiso} any surface carried by $\intE$ is
  projectively isotopic to one carried by $\intC$.  Thus $E$ must be
  $C$-independent and we know that it contains $\vF^*$ which is an
  interior point of the \emph{maximal} $C$-independent face $D$;
  consequently, we must have $E = D$ and so $G$ is carried by $D$.
  Thus $D$ is complete as claimed.
\end{proof}

\begin{proof}[Proof of Theorem~\ref{thm:disjointdep}]
We start with:
\begin{claim}\label{lem:depcarry}
  Any lw-surface $F$ is carried by $\dep{C}$ for some complete lw-face
  $C$.
\end{claim}
By Theorem~\ref{thm.toll}, the surface $F$ is carried by some complete
lw-face.  It is immediate from the definition that the intersection of
two complete lw-faces is again complete, so there exists a minimal
complete lw-face $C$ carrying $F$.  Let $D$ be the face of $C$
containing $F$ in its interior.  If $D$ is $C$-dependent, then
$\intD \subset \dep{C}$ and so $F \in \dep{C}$ as desired. So assume
$D$ is $C$-independent.  Let $E$ be a maximal $C$-independent face of
$C$ containing $D$, and note $E \neq C$ as $C$ is $C$-dependent. By
Lemma~\ref{lem:maxCindep}, the face $E$ is complete and so we have
found a smaller complete face containing $F$ than $C$, a
contradiction.  So we have proven Claim~\ref{lem:depcarry}.

To prove Theorem~\ref{thm:disjointdep} it remains to show:
\begin{claim}\label{lem:depdisjoint}
  If $C_1$ and $C_2$ are distinct (but perhaps not disjoint)
  complete lw-faces, then no surface carried by $\dep{C_1}$
  is projectively isotopic to one carried by $\dep{C_2}$.  In
  particular, the sets $\dep{C_1}$ and $\dep{C_2}$ are disjoint.
\end{claim}
Suppose not and that $F_1$ and $F_2$ are projectively isotopic normal
surfaces carried by $\dep{C_1}$ and $\dep{C_2}$ respectively.  By
Theorem~\ref{thm.compiso}, we can further assume each $F_i$ is carried
by $\intC_i$.  Replacing them with multiples if necessary, we can
assume that they are actually isotopic.  By Corollary~4.6 of
\cite{To:isotopy}, it follows that both $F_1$ and $F_2$ must be
carried by $C_1 \cap C_2$; as each $F_i$ is carried by $\intC_i$, we
must have $C_1 = C_2$, a contradiction.  This proves
Claim~\ref{lem:depdisjoint} and hence the theorem.
\end{proof}



\section{Surface counts are almost quasi-polynomial}
\label{sec:counts}

The first of this section's two main results is:
\begin{theorem}\label{thm.oneface}
  Suppose $M$ is a closed irreducible atoroidal 3-manifold that
  contains no nonorientable essential surfaces and $C$ is a complete
  lw-face of $\PT$.  Let $b_C(n)$ be the number of isotopy classes of
  closed essential surfaces $F$ carried by $\dep C$ with $\chi(F) = n$,
  and let $B_C(x) = \sum_{n=1}^{\infty} b_C(-2n) x^n$ be the corresponding
  generating function.  Then $B_C(x)$ is short.
\end{theorem}
The other main result of this section is
Theorem~\ref{thm.idealoneface}, which is the analog of
Theorem~\ref{thm.oneface} when $M$ has boundary. Combining
Theorem~\ref{thm.oneface} with the results from the last section, we
can now give:
\begin{proof}[Proof of Theorem~\ref{thm.main} when $M$ is closed]
  As $M$ is closed, a (necessarily closed) surface in $M$ is essential
  exactly when it is incompressible.  Therefore, by
  Theorem~\ref{thm:disjointdep}, each isotopy class of essential
  surface is carried by $\dep{C}$ for a unique complete lw-face $C$.
  As sums of short generating functions are also short,
  Theorem~\ref{thm.main} now follows from Theorem~\ref{thm.oneface}.
\end{proof}

\subsection{Counting surfaces via lattice points}

We now turn to the proof of Theorem~\ref{thm.oneface}, so let $M$ be a
closed irreducible atoroidal \3-manifold with triangulation $\cT$.
From now on, fix a complete lw-face $C$ of $\PT$ and consider the cone
$\Ctil \subset \ST$, that is
$\R_{\geq 0} \cdot C = \setdef{t \xx}{t \in \R_{\geq 0}, \xx \in C}$.
Recall that $\dep C \subset C$ is the complement of its
$C$-independent faces, and define
$\dep{\Ctil} = \left(\R_{\geq 0} \cdot \dep{C}\right) \setminus
\{\vec{0}\}$.  Now $b_C(n)$ in Theorem~\ref{thm.oneface} is the number
of isotopy classes of normal surfaces $F$ with $\vF \in \dep{\Ctil}$
and $\chi(F) = n$.

As motivation, let us start with the easy case when the subspace $W_C$
from Theorem~\ref{thm.compiso} is zero.  Then
$\dep{\Ctil} = \intCtil$, and two surfaces $F$ and $G$ carried by
$\dep{\Ctil}$ are isotopic if and only if $\vF = \vG$. In our
triangle-quad coordinates, there is a linear function
$\chi \maps \R^{7t} \to \R$ such that the Euler characteristic of a
normal surface $F$ is given by $\chi(\vF)$, see
\cite[Algorithm~9.1]{JacoTollefson1995}.  Every normal surface $F$
carried by $C$ is incompressible and hence $\chi(F) \leq -1$ as $M$ is
closed, irreducible, and atoroidal.  Thus $\chi < 0$ on every vertex
of $C$ which implies $\chi$ is proper on $\Ctil$ and so the set
$X = \setdef{\xx \in \Ctil}{\chi(x) = -1}$ is a compact polyhedron.
Now $b_C(-n)$ is simply the size of the set $(n \intX) \cap \Z^{7t}$,
and counting lattice points in dilations of a compact polyhedron has
been studied extensively starting with the work of Ehrhart in the
1960s.  In particular, Theorem~4.6.26 of \cite{Stanley}, whose proof
uses Ehrhart-Macdonald reciprocity, tells us that the generating
function $B_C(x)$ is short, proving Theorem~\ref{thm.oneface} when
$W_C = 0$.

When $W_C$ is nonzero, to count isotopy classes of surfaces we need to
identify lattice points in $\dep{\Ctil}$ that differ by an
element of $W_C$.  We do so in the following way.  Let $V$ be the
linear subspace of $\R^{7 t}$ spanned by all vectors in $C$, and let $W$ be
$W_C$.  Define $V(\Z) = V \cap \Z^{7t}$ and $W(\Z) = W \cap \Z^{7t}$.
Using Smith normal form, we can find a complementary rational subspace
$L \subset V$ to $W$ such that the lattice $V(\Z)$ is the direct sum
of $W(\Z)$ and $L(\Z) = L \cap \Z^{7t}$.  Let $T \maps V \to L$ be the
projection operator associated with the decomposition
$V = W \oplus L$.   We can now turn our question of counting isotopy
classes of surfaces into one about counting certain lattice points in
$L(\Z)$:
\begin{lemma}\label{lem:normbij}
  The set $T\left(\dep{\Ctil} \cap \Z^{7t}\right)$ is in bijection with
  isotopy classes of normal surfaces carried by $\dep{C}$.
\end{lemma}

\begin{proof}
  Normal surfaces carried by $\dep{C}$ correspond to lattice points in
  $\dep{\Ctil}$.  As $W$ is the kernel of $T$, the claim is equivalent
  to saying that if $F$ and $G$ are normal surfaces with $\vF$ and
  $\vG$ in $\dep{\Ctil}$, then $F$ is isotopic to $G$ if and only if
  $\vF - \vG \in W$.  As we are assuming that all incompressible
  surfaces in $M$ are orientable, this follows immediately from
  Theorem~\ref{thm.compiso}.
\end{proof}

To prove Theorem~\ref{thm.oneface}, we will need a tool for counting
points in sets such as $T\left(\dep{\Ctil} \cap \Z^{7t}\right)$.
Recently, Nguyen and Pak \cite{N-Pak}, building on
\cite{Barvinok-Woods}, established exactly the result we need here.
To apply \cite{N-Pak}, we need the linear map $T$ to be
\emph{integral} in the sense that its matrix with respect to any
$\Z$-bases of $V(\Z)$ and $L(\Z)$ has integer entries, but that is
clear from its definition.  Since we want to count by Euler
characteristic, we first study $\chi \maps V \to \R$:

\begin{lemma}\label{lem:eulerfn}
  The restriction $\chi \maps V(\Z) \to \R$ is integral and, since $M$
  is closed, irreducible, and atoroidal, the function $\chi$ is proper
  on $\Ctil$ and negative on $\Ctil \setminus \{ \vec{0} \}$.
\end{lemma}

\begin{proof}
  For $\chi|_V$, note that $\Ctil$ has nonempty interior as a subset
  of $V$ and contains open balls of arbitrary size.  Hence, given any
  $\vv \in V(\Z)$, we can find $\xx, \yy \in \Ctil(\Z)$ with
  $\vv = \xx - \yy$.  There are normal surfaces $F$ and $G$ with
  $\vF = \xx$ and $\vG = \yy$, and so $\chi(\vv) = \chi(F) - \chi(G)$
  is in $\Z$ as needed to show $\chi|_V$ is integral. 

  For $\chi|_\Ctil$, every normal surface $F$ carried by $C$ is
  incompressible and hence $\chi(F) \leq -1$ as $M$ is closed,
  irreducible, and atoroidal.  Thus $\chi < 0$ on every vertex of $C$
  which implies it is proper on $\Ctil$ and negative on
  $\Ctil \setminus \{ \vec{0} \}$.
\end{proof}

Now we combine $\chi$ and $T$ as follows.  Define
$\Tbar \maps V \to L \oplus \R$ by
$\Tbar(\xx) = \left(T(\xx), -\chi(\xx) \right)$, which is integral as
both its component functions are, and we have:
\begin{lemma}\label{lem:normbijbar}
  The set $\barT\left(\dep{\Ctil} \cap \Z^{7t}\right)$ is in bijection with
  isotopy classes of normal surfaces carried by $\dep{C}$.
\end{lemma}
\begin{proof}
  By Lemma~\ref{lem:normbij}, it suffices to show that projecting away
  the second factor of $L \oplus \R$ gives a bijection between
  $\Tbar\left(\dep{\Ctil} \cap \Z^{7t}\right)$ and
  $T\left(\dep{\Ctil} \cap \Z^{7t}\right)$.  This projection is
  clearly onto, so this reduces to showing that for normal surfaces
  $F$ and $G$ with $\vF$ and $\vG$ in $\dep{C}$ and $T(\vF) = T(\vG)$
  then $-\chi(\vF) = -\chi(\vG)$.  But the latter holds since
  $\vF - \vG \in W$ implies the surfaces $F$ and $G$ must be isotopic
  by Theorem~\ref{thm.compiso} and thus homeomorphic.
\end{proof}

To apply \cite{N-Pak}, we will need one more property about $T$ and $\Tbar$:

\begin{lemma}\label{lem:orthant}
  There is a lattice $L'(\Z)$ containing $L(\Z)$ which has a
  $\Z$-basis such that $T(\Ctil) \subset \Rp^d$ under the induced
  identification of $L$ with $\R^d$.  Moreover, the same holds for $\Tbar$.
\end{lemma}

\begin{proof}
  The claim for $\Tbar$ follows immediately from that for $T$ since
  $-\chi(\Ctil) = [0, \infty)$ by Lemma~\ref{lem:eulerfn}. So
  now we consider only $T$.

  Fourier-Motzkin elimination tells us that the image $T(\Ctil)$ is
  again a polyhedral cone.  We first show that $T(\Ctil)$ is a
  \emph{pointed cone}, that is, one that does not contain a line. By
  Theorem~\ref{thm.compiso}, we know $W \subset
  \ker(\weight)$. Therefore, the map $\weight \maps V \to \R$ factors
  through the projection $T \maps V \to L$.  Hence all of
  $T(\Ctil) \setminus \{\vec{0}\}$ is strictly to the positive side of
  the hyperplane $(\weight|_L)^{-1}(0)$ and so $T(\Ctil)$ is
  a pointed cone.

  We will now find a basis for $L(\Q)$ as a $\Q$-vector space with the
  property that $T(\Ctil)$ lies in the postive orthant; this suffices
  to prove the lemma as we can scale the basis elements by $a > 0$ in
  $\Q$ so that the lattice they generate contains $L(\Z)$.  Let
  $\{\vv_i\}$ denote the vertices of $C$.  Note that if we can find a
  basis $\{\ell_j\}$ of $L(\Q)^* = \Hom\big(L(\Q), \Q\big)$ where
  $\ell_j(\vv_i) > 0$ for all $i$ and $j$, then the algebraically dual
  basis $\evec_k$ of $L(\Q)$, that is, the one where
  $\ell_j(\evec_k) = \delta_{jk}$, is the basis we seek.  Fix any
  basis $\{\beta_j\}$ of $L(\Q)^*$ where $\beta_1 = \weight$ and for
  $\epsilon \in \Q^\times$ consider the new basis $\{\ell_j\}$ where
  $\ell_1 = \beta_1$ and all other
  $\ell_j = \beta_1 + \epsilon \beta_j$.  Since we showed above that
  $\beta_1(\vv_i) = \weight(\vv_i) > 0$ for each $i$, for small enough
  $\epsilon$ we have $\ell_j(\vv_i) > 0$ for all $i$ and $j$ as needed
  to prove the lemma.
\end{proof}

Next, we introduce the language needed to state the conclusion
of \cite{N-Pak}.  A set $A$ of points in $\N^n$ has an associated
generating function:
\[
  f_A(\bt) = \sum_{\va \in A} \bt^{\va}
  \mtext{in $\Z[[t_1,\ldots,t_n]]$ where $\bt^{\va} = t_1^{a_1} \cdots
    t_n^{a_n}$ for $\va = (a_1, \ldots, a_n)$.}
\]
We say that $A$ has a \emph{short generating function} when there are $c_i
\in \Q$ and $\va_i, \vb_{i j} \in \Z^{n}$ such that:
\begin{equation}\label{eq:multishort}
  f_A(\bt) = \sum_{i = 1}^{N} \frac{c_i \bt^{\va_i}}{%
    \big(1 - \bt^{\vb_{i 1}}\big) \cdots \big(1 - \bt^{\vb_{i k_i}}\big)} \,.
\end{equation}
These multivariable short generating series were introduced by
Barvinok and play a key role in polynomial time algorithms for
counting lattice points in convex
polyhedra~\cite{Barvinok:polynomial}.

Using the lattice $L'(\Z) \oplus \Z \subset L \oplus \R$, where
$L'(\Z)$ is from Lemma~\ref{lem:orthant}, we henceforth view
$\Tbar\big(\dep{\Ctil} \cap \Z^{7t}\big)$ as subset of $\N^{d +
  1}$. The key to Theorem~\ref{thm.oneface} is:

\begin{lemma}\label{lem:bigshort}
  The set $\Tbar\big(\dep{\Ctil} \cap \Z^{7t}\big)$ has a short
  generating function.
\end{lemma}

\begin{proof}
  We will construct a rational polyhedron $Q \subset \Ctil$ such
  that
  \[
    Q \cap \Z^{7t} = \dep{\Ctil} \cap \Z^{7t}
  \]
  By Lemma~\ref{lem:orthant}, we have $\Tbar(Q) \subset \Rp^{d+1}$
  using the basis of $L'(\Z)$ given there.  Additionally, the
  projection $\Tbar$ is integral with respect to the lattices $V(\Z)$
  and $L'(\Z)$ since $L'(\Z) \supset L(\Z)$.  Therefore, Theorem~1.1
  of Nguyen-Pak \cite{N-Pak} will apply and give that the generating
  function for $\Tbar\big(Q \cap \Z^{7t}\big)$ is short, proving the
  lemma.

  Now $\dep{\Ctil}$ is simply $\Ctil$ with some closed faces removed,
  and we can use the following standard trick to construct $Q$.  For a
  face $D$ of $C$ its \emph{active variables} are
  \[
    I_D = \setdefm{\big}{i \in [1, 2, \ldots, 7t]}%
          {\mbox{$x_i = 0$ on $D$ but $x_i > 0$ somewhere on $C$}}
  \]
  Thus $D$ is the subset of $C$ cut out by $x_i = 0$ for $i \in
  I_D$, or equivalently the locus where $\sum_{i \in I_D} x_i =
  0$ since each $x_i \geq 0$ on $C$.  Then
  $\dep{\Ctil}$ consists of those $\xx \in V$ where:
  \begin{enumerate}[noitemsep, topsep=0.4em]
  \item all $x_i \geq  0$,
  \item for each $C$-independent face $D$ one has $\sum_{i \in I_D} x_i
    > 0$,
  \item and finally $\sum_{i=1}^{7t} x_i >
     0$ as the origin is not in $\dep{\Ctil}$.
   \end{enumerate}
   If we define $Q$ to be those $\xx \in V$ where all $x_i \geq 0$,
   where for each $C$-independent face $D$ one has
   $\sum_{i \in I_D} x_i \geq 1$, and finally where
   $\sum_{i=1}^{7t} x_i \geq 1$, then we have
   $Q \cap \Z^{7t} = \dep{\Ctil} \cap \Z^{7t}$ as needed.
\end{proof}

\begin{proof}[Proof of Theorem~\ref{thm.oneface}]
  Let $f(\bt)$ be the generating function for
  $\Tbar\big(\dep{\Ctil} \cap \Z^{7t}\big)$.  The variable $t_{d + 1}$ in
  $f(\bt)$ corresponds to $-\chi$ and by Lemma~\ref{lem:eulerfn} the
  function $\chi$ is proper on $\Ctil$; thus, there are only finitely
  many terms of $f(\bt)$ with any given power of $t_{d+1}$.  Hence
  $g(t) = f(1, \ldots, 1, t)$ is a well-defined element of $\Z[[t]]$,
  and indeed by Lemma~\ref{lem:normbij} it is the generating function
  $B_C(x)$ we seek with $x$ replaced by $t^2$.

  Since $f(\bt)$ is short by Lemma~\ref{lem:bigshort}, it remains to
  use this to see that $g(t)$ is also short. Provided no denominator
  in (\ref{eq:multishort}) has a factor of
  $\left(1 - t_1^{a_1} t_2^{a_2} t_3^{a_3} \dots t_d^{a_{d}}\right)$,
  that is, has no $t_{d+1}$-term, then this is immediate. To handle
  the general case, we will use results from \cite{Woods}, noting that
  our notion of a generating function being short is equivalent to
  rationality in the sense of Definition 1.4 of \cite{Woods}.  First,
  set $S = \Tbar\big(\dep{\Ctil} \cap \Z^{7t}\big)$.  As the
  generating function $f(\bt)$ of $S$ is short, Theorem~1.5 of
  \cite{Woods} gives that $S$ is a Presburger set, that is, there is a
  Presburger formula $F$ which tests points in $\N^{d+1}$ for
  membership in $S$.  Writing points in $\N^{d + 1}$ as $(\vc, p)$
  with $\vc \in \N^d$ and $p \in \N$, we see from Definition 1.6 of
  \cite{Woods} that $g(t)$ is the generating function of the
  Presburger counting function
  $p \mapsto \# \setdef{\vc \in \N^d}{F(\vc, p)}$.  Therefore, by
  Theorem~1.10 of \cite{Woods}, specifically $A \Rightarrow C$, the
  generating function $g(t)$ is short.
\end{proof}

\subsection{Ideal triangulations}
\label{sec:ideal}

Our proof of Theorem~\ref{thm.main} when $\partial M$ is nonempty will
use ideal triangulations instead of finite ones.  The theory of
\emph{closed} normal surfaces in ideal triangulations is nearly
identical to that of normal surfaces in finite triangulations; after
all, normal surfaces stay away from the vertices of the ambient
triangulation, which is the only place where the topology differs
between the two cases.  Indeed, we claim that all the results of
\cite{To:isotopy} hold for closed normal surfaces in ideal
triangulations without any changes to the proofs. Manifolds with
boundary are allowed in \cite{To:isotopy}, so switching from finite to
ideal triangulations can be viewed as using a slightly different type
of finite cellulation as the background to do normal surface
theory, specifically a cellulation by \emph{truncated} tetrahedra.
The combinatorics of closed normal surfaces in truncated
tetrahedra is almost indistinguishable from standard normal
surface theory in a finite triangulation, and hence the proofs in
\cite{To:isotopy} work as written in our new context.
Consequently, the results of Section~\ref{sec.thm.main} also hold
for closed surfaces in ideal triangulations.

From now on, suppose $M$ is a compact irreducible \bdry irreducible
\3-manifold with nonempty boundary and $\cT$ an ideal triangulation of
$M$, which exists by
e.g.~\cite[Proposition~3]{JacoRubinsteinSpreerTillmann2018}.  The
\emph{vertex link} $H_v$ of a vertex $v \in \cT^0$ is the normal
surface consisting of one triangle in each tetrahedron corner where
the vertex in that corner corresponds to $v$. The vertex link $H_v$
should be viewed as a parallel copy of the corresponding boundary
component of $M$.  An ideal triangulation $\cT$ is \emph{\bdry
  efficient} when the only connected normal surfaces that are boundary
parallel are the vertex links.  For example, if $\cT$ has a positive
angle structure then it is \bdry efficient by \cite[Proposition
4.4]{Lackenby2000}.  Provided $M$ is acylindrical, then any minimal
ideal triangulation is \bdry efficient by \cite[Theorem
4]{JacoRubinsteinSpreerTillmann2018}, so such triangulations always
exist for the manifolds we consider in Theorem~\ref{thm.main}.

Our goal now is to weed out the inessential incompressible surfaces,
i.e.~those with a \bdry parallel component, from our counts. When
$\partial M$ includes a torus, this is not just an aesthetic
preference but a requirement since there are infinitely many
isotopy classes of (disconnected) incompressible surfaces with the
same Euler characteristic.

\begin{lemma}\label{lem:allornone}
  Suppose $\cT$ is a \bdry efficient ideal triangulation of a
  3-manifold $M$ that contains no nonorientable closed incompressible
  surfaces. Let $C$ be a lw-face of $\PT$.  If $C$ carries no vertex
  link then every normal surface carried by $C$ is essential.  If $C$
  carries some vertex link then no normal surface carried by $\dep{C}$
  is essential.
\end{lemma}
\begin{proof}
  Let $I_C \subset \{1, 2, \ldots, 7t\}$ be the indices of the
  coordinates on $\R^{7t}$ which vanish on all of $C$.  Then
  $
    C = \setdef{\xx \in \PT}{\mbox{$x_i = 0$ for all $i \in I_C$}}
    $ and $
    \intC = \setdef{\xx \in C}%
    {\mbox{$x_i > 0$ for all $i \notin I_C$}}.
  $

  First, suppose $C$ carries an inessential normal surface $F$. As
  $\cT$ is \bdry efficient, the surface $F$ is the disjoint union of a
  normal surface $G$ (possibly empty) and some vertex link $H_v$; in
  particular $F = G + H_v$.  From the above description of $C$, it is
  clear that as $G + H_v$ is carried by $C$, both $G$ and $H_v$ are also
  carried by $C$.  This proves the first claim.

  Second, suppose $C$ carries some vertex link $H_v$.  Since $H_v$ is
  carried by $C$, for each index $i$ corresponding to a triangle in
  $H_v$ we have $i \notin I_C$.  Hence, for any normal surface $F$
  carried by $\intC$, each triangle that appears in $H_v$ also has
  positive weight in $F$.  Consequently, the surface $F$ has a
  component which is normally isotopic to some $H_v$ and in particular
  is inessential.  More broadly, suppose $G$ is a normal surface
  carried by $\dep{C}$.  As $C$ is least-weight, by
  Theorem~\ref{thm.compiso}, the surface $G$ is projectively isotopic
  to some $F$ carried by $\intC$, and by the previous argument the
  latter has a component which is $H_v$.  By the hypotheses on $M$,
  both $G$ and $F$ are orientable since they are incompressible; as
  they are projectively isotopic, it follows that $G$ also has a
  component isotopic to $H_v$.  In particular, the surface $G$ is
  inessential. This proves the second claim.
\end{proof}

We call a lw-face $C$ \emph{essential} if every normal surface carried
by $C$ is essential.  By Lemma~\ref{lem:allornone}, a lw-face $C$ is
either essential or \emph{every} normal surface carried by $\dep{C}$
is inessential.  Hence Theorem~\ref{thm:disjointdep} gives:
\begin{theorem}
  \label{thm:disjointdepessential}
  Suppose $\cT$ is a \bdry efficient ideal triangulation of a
  3-manifold $M$ that contains no nonorientable closed incompressible
  surfaces.  For each orientable essential surface $F$ there exists a
  unique complete essential lw-face $C$ such that $\dep{C}$ carries a
  surface (non-projectively) isotopic to $F$.
\end{theorem}

We can now prove the analog of Theorem~\ref{thm.oneface} for manifolds
with boundary: 
\begin{theorem}\label{thm.idealoneface}
  Suppose $M$ is an irreducible \bdry irreducible atoroidal
  acylindrical 3-manifold with $\partial M \neq \emptyset$ that
  contains no nonorientable essential surfaces.  Suppose $\cT$ is a
  \bdry efficient ideal triangulation of $M$ and $C$ is a complete
  essential lw-face of $\PT$.  The generating function $B_C(x)$
  corresponding to the counts of isotopy classes of closed essential
  surfaces carried by $\dep C$ is short.
\end{theorem}
\begin{proof}
  As $C$ carries only essential surfaces, we have
  $\chi \maps \Ctil \to \R$ is proper since all essential surfaces
  have $\chi \leq -1$. This gives the analog of
  Lemma~\ref{lem:eulerfn} in our setting, and the proof of the theorem
  is now identical to that of Theorem~\ref{thm.oneface}.
\end{proof}

We now complete the proof of the first main theorem of this paper:

\begin{proof}[Proof of Theorem~\ref{thm.main} when $M$ has boundary]
  Take $\cT$ to be a minimal ideal triangulation of $M$, which is
  \bdry efficient by
  \cite[Theorem~4]{JacoRubinsteinSpreerTillmann2018} since $M$ is
  acylindrical.  By Theorem~\ref{thm:disjointdepessential}, every
  isotopy class of essential surface is carried by $\dep{C}$ for a
  unique complete essential lw-face $C$.  As sums of short generating
  functions are also short, Theorem~\ref{thm.main} now follows from
  Theorem~\ref{thm.idealoneface}.
\end{proof}

\subsection{Whither nonorientable and bounded surfaces}
\label{sec:nonorient}

We would of course like to remove the hypothesis in
Theorem~\ref{thm.main} that $M$ contains no closed nonorientable
essential surfaces, and also broaden the count to allow essential
surfaces with boundary.  We now outline some of the difficulties
inherent in such extensions.

For nonorientable closed essential surfaces, what prevents us from
just including them in the count is that while $\LW$ can carry
nonorientable essential surfaces, it need not carry \emph{all} of
them.  (This distinction is also present in the branched surface
perspective of \cite{Oertel:branched}.)  The issue is that you can
have a nonorientable essential normal surface $F$ which is least
weight in its isotopy class but where its double $2 F$, while normal
and essential, may not be least weight.  This does not happen for an
orientable $F$ since the double is just two parallel copies of $F$.
One could sidestep this issue by just counting orientable surfaces,
but picking those out of each lw-face seems tricky for the following
reason. Note that $F$ is orientable if and only if $2 F$ has twice the
number of connected components as $F$.  As we see from the $a_M$
versus $b_M$ discussion in Sections~\ref{sec.introduction} and
\ref{sec.genuspat}, counting components is subtle.  Consequently, we
suspect there are examples where the count of orientable surfaces does
not have a short generating function.

A key obstruction to counting surfaces with boundary is actually the
issue of orientability.  Unlike in the closed case with
Proposition~\ref{prop:nonorient}, there is no homological condition we
can impose that a priori eliminates the possibility of nonorientable
essential surfaces with boundary.  For example, the exterior of a knot
in $S^3$ can contain such nonorientable surfaces (e.g.~many
checkerboard surfaces for alternating knots).  As \cite{To:isotopy}
and our Section~\ref{sec.thm.main} do allow orientable surfaces with
boundary, we are hopeful that if nonorientable closed surfaces can be
dealt with, then counting bounded surfaces will also be possible.

\section{Proof of the decision theorem}
\label{sec:decision}

This section is devoted to proving Theorem~\ref{thm.algol}, which says
that there are algorithms for finding the generating function in
Theorem~\ref{thm.main} as well as enumerating representatives of the
isotopy class of essential surfaces and determining which isotopy
class a given surface belongs to.  In this section, we do not worry
about the efficiency of these algorithms, merely their existence; the
actual method used to compute the examples in
Section~\ref{sec.computations} uses some of the ideas here but in the
modified form of Section~\ref{sec.ideal} which is specific to when
$\partial M$ is a nonempty union of tori.

Throughout this section, let $M$ be a compact orientable irreducible
$\partial$-irreducible 3-manifold with a fixed triangulation $\calT$,
which is a finite triangulation when $M$ is closed or an ideal
triangulation otherwise.  With the exception of the proof of
Theorem~\ref{thm.algol} itself at the very end, in this section we do
\emph{not} require that $M$ is acylindrical or atoroidal, nor that
$\calT$ is $\partial$-efficient; moreover, the manifold $M$ may
contain nonorientable closed essential surfaces. The algorithms in
Theorem~\ref{thm.algol} follow the approach of the proof of
Theorem~\ref{thm.main} closely, so the first thing we will need is:
\begin{theorem}\label{thm.lwfaces}
  There exists an algorithm for computing the collection $\CLW$ of
  complete lw-faces of $\cP_\cT$.
\end{theorem}
Given a normal surface $F$ in $\cT$, one can algorithmically determine
whether or not it is incompressible and $\partial$-incompressible;
indeed, this is essentially Haken's original application of normal
surface theory, see e.g.~\cite[Algorithm 9.6]{JacoTollefson1995}.  The
tricky part of computing $\CLW$ is figuring out the isotopy
relationships between different such normal surfaces.

\subsection{Graph of incompressible surfaces}
\label{sec:graphincomp}

We use the following framework for understanding isotopies among
normal surfaces.  Let $\calG_\cT$ be the graph whose vertices are
connected incompressible closed normal surfaces, more precisely the
normal isotopy classes of such surfaces, and where there is an edge
between surfaces $F$ and $G$ exactly when $F$ and $G$ can be normally
isotoped to be disjoint and cobound a product region. Given an integer
$w$, we use $\cG_\cT^{\leq w}$ to denote the subgraph whose vertices
are all surfaces $F \in \cG_\cT$ of weight at most $w$.  Since there
are only finitely many normal surfaces of bounded weight, each
$\cG_\cT^{\leq w}$ is finite.  Moreover, given $w$, the graph
$\cG_\cT^{\leq w}$ can be algorithmically constructed as follows.
First, the vertices of $\cG_\cT^{\leq w}$ can be found by enumerating
all connected normal surfaces of weight at most $w$ and then testing
each for incompressibility.  Second, for each pair of surfaces $F$ and
$G$ in $\cG_\cT^{\leq w}$, one can test if they can be normally
isotoped to be disjoint using Algorithm 9.5 of
\cite{JacoTollefson1995}; specifically, this can be done if and only
if $F$ and $G$ are compatible and the normal sum $F + G$ consists of
two connected components where one is normally isotopic to $F$ and the
other to $G$.  When they can be made disjoint in this way, there is a
unique way to do so up to normal isotopy.  Finally, for each pair of
surfaces $F$ and $G$ that can be normally isotoped apart, we test all
components of $M$ cut along $F \cup G$ for being products using
Algorithm 9.7 of \cite{JacoTollefson1995}.  This completes the
algorithm for constructing $\cG_\cT^{\leq w}$.

Any two surfaces in the same connected component of $\cG_\cT$ are of
course isotopic.  It turns out the converse is true as well, in the
following strong form:
\begin{theorem}\label{thm.path}
  If $F$ and $G$ are two vertices of $\cG_\cT$ that are isotopic, then
  they are joined by a path in $\cG_\cT$ passing only through vertices
  $H$ with $\weight(H) \leq \max(\weight(F), \weight(G))$.  In
  particular, the isotopy classes of surfaces in any
  $\cG_\cT^{\leq w}$ correspond precisely to the connected components
  of $\cG_\cT^{\leq w}$.
\end{theorem}
We prove this theorem in Section~\ref{sec.normpairs} below, but we first
use it to derive Theorem~\ref{thm.lwfaces}.

\begin{proof}[Proof of Theorem~\ref{thm.lwfaces}]
  First, compute the polytope $\cP_\cT$ from the normal surface
  equations.  For each face $C$ of $\cP_\cT$, fix a normal surface
  $F_C$ that is carried by its interior.  Compute the graph
  $\cG_\cT^{\leq w}$ where $w$ is the maximum weight of any $F_C$.
  The incompressible $F_C$ are those that are vertices of
  $\cG_\cT^{\leq w}$, and, by Theorem~\ref{thm.path}, we know exactly
  which $F_C$ are least-weight.  When $F_C$ is least-weight, we also
  know every other least-weight surface isotopic to it.  Applying
  Theorem~\ref{thm.interior} now identifies exactly the faces $C$ that
  are in $\CLW$.
\end{proof}

\subsection{Isotopic normal pairs}
\label{sec.normpairs}

Throughout, let $M$ be a compact orientable irreducible 3-manifold
with a fixed triangulation $\calT$ as in the previous section.  An
\emph{isotopic normal pair} $(F, G)$ is an isotopic pair of closed
incompressible normal surfaces $F$ and $G$ that meet transversely in
the sense of \cite[Page 1091]{To:isotopy}.  Define the complexity of
such a pair by
\[
  c(F, G) = \big( \max(\weight(F), \weight(G)), \
                  \min(\weight(F), \weight(G)), \
                  \#(F \cap G) \big)
\]
where $\#(F \cap G)$ denotes the number of connected components of $F
\cap G$.  We will compare complexities lexicographically.

If $(F, G)$ is an isotopic normal pair where $F$ and $G$ are disjoint,
then by Lemma~5.3 of \cite{Waldhausen} the surfaces $F$ and $G$ are
parallel, i.e.~cobound a region homeomorphic to $F \times I$.  In
this case, the pair $(F, G)$ gives rise to an edge of
$\cG_\cT$. The key result of this subsection is:

\begin{theorem}\label{thm.pairs}
  If $(F, G)$ is an isotopic normal pair with $F \cap G \neq
  \emptyset$ then, after possibly interchanging $F$ and $G$, there
  exists a normal surface $F'$ that is isotopic to $F$ and disjoint
  from it that meets $G$ transversely with $c(F', G) < c(F, G)$.
\end{theorem}
Given any $(F, G)$ isotopic normal pair with
$w = \max(\weight(F), \weight(G))$, we can apply
Theorem~\ref{thm.pairs} repeatedly until we arrive at a pair
$(F'', G'')$ where $F''$ and $G''$ are disjoint.  This proves
Theorem~\ref{thm.path} above, since each application of
Theorem~\ref{thm.pairs} gives an edge in $\cG_\cT^{\leq w}$ and there is
also an edge from $F''$ to $G''$.

Suppose $\wt F$ and $\wt G$ are subsurfaces of $F$ and $G$
respectively, with $\partial \wt F = \partial \wt G$.  Here is one way
to make precise the notation that $\wt F$ and $\wt G$ are ``parallel
rel boundary''.  Given a compact surface $H$, define $P(H)$ as the
quotient of $H \times I$ where for each $h \in \partial H$ the set
$\{h\} \times I$ has been collapsed to a point.  A \emph{product
  region} between $\wt F$ and $\wt G$ is an embedding $f: P(H) \to M$
where $f(H \times \{0\}) = \wt F$ and $f(H \times \{1\}) = \wt G$;
here we do not insist that $P = f(P(H))$ meets $F \cup G$ only in
$\wt F \cup \wt G$. The first step in proving Theorem~\ref{thm.pairs}
is to show:

\begin{lemma}\label{lem.pairs}
  Suppose $(F, G)$ is an isotopic normal pair with
  $F \cap G \neq \emptyset$. After possibly interchanging $F$ and $G$,
  there exist subsurfaces $\wt F \subset F$ and $\wt G \subset G$
  where $\wt G \cap F = \partial \wt G$ and
  $\partial \wt F = \partial \wt G$ and
  $\weight( \wt G) \leq \weight( \wt F)$ with $\wt G$ and $\wt F$
  bounding a product region $P$ with $P \cap F = \wt F$.
\end{lemma}

\begin{proof}
  Suppose all components of $F \cap G$ are essential in both $F$ and
  $G$.  Then by Proposition~5.4 of \cite{Waldhausen}, there exist
  subsurfaces $\wt F \subset F$ and $\wt G \subset G$ with
  $\wt F \cap G = \partial \wt F = \partial \wt G = F \cap \wt G$
  where $\wt F \cup \wt G$ bounds a product region $P$ where
  $P \cap F = \wt F$ and $P \cap G = \wt G$.  Relabeling, we can
  arrange that $\weight( \wt G) \leq \weight( \wt F)$ to prove the
  lemma in this case.

  Suppose instead some component of $F \cap G$ is inessential in one
  of $F$ or $G$. Among all disks contained in one of $F$ or $G$
  bounded by a component of $F \cap G$, let $D$ be one of least
  weight, which exists since the weight (i.e., the number of
  intersection points with the 1-skeleton of $\calT$) of any
  subsurface is a nonnegative integer.  By passing to an innermost
  component, we can assume $D$ meets $F \cap G$ only along
  $\partial D$.  After relabeling, we can assume this $D$ is contained
  in $G$ and then set $\wt G = D$.  As $F$ is incompressible, the
  curve $\partial \wt G $ must bound a disk $\wt F$ in $F$. Together
  the disks $\wt F \cup \wt G$ form a sphere which must bound a ball
  as $M$ is irreducible, and hence $\wt F$ and $\wt G$ bound the
  required product region $P$.  By our initial choice of $\wt G$, we
  must have $\weight( \wt G) \leq \weight( \wt F)$ as desired.
\end{proof}

\begin{proof}[Proof of Theorem~\ref{thm.pairs}]
  Let $\wt F$ and $\wt G$ be given by Lemma~\ref{lem.pairs}.  Set
  $F_0 = (F \setminus \wt F) \cup \wt G$ which is isotopic to $F$ via
  the product region $P$.  Move $F_0$ slightly so that it is disjoint
  from $F$, and notice that
  $\weight(F_0) = \weight(F) - \weight(\wt F) + \weight(\wt G) \leq
  \weight(F)$.  If $F'$ is a normalization of the incompressible
  surface $F_0$, we have $\weight(F') \leq \weight(F_0)$; by the
  barrier theory \cite[Theorem 3.2(1)]{JacoRubinstein2003}, the
  surface $F'$ is disjoint from $F$, and we can perturb $F'$ slightly
  to be transverse to $G$.

  If $\weight(F') < \weight(F)$ we now have our desired $(F', G)$ as
  $c(F', G) < c(F, G)$ where the two complexities differ in one of the
  first two components.  If instead $\weight(F') = \weight(F)$, then
  we must have $\weight(F_0) = \weight(F)$.  This means that $F_0$ is
  \emph{normally} isotopic to $F'$ as all normalization moves that
  change the normal isotopy class strictly reduce the weight.  Then
  $\#(F' \cap G) = \#(F_0 \cap G) = \#(F \cap G) - \#\partial \wt F <
  \# (F \cap G)$ and hence $c(F', G) < c(F, G)$ with the complexities
  differing only in the last component.
\end{proof}

We turn now to the proof of Theorem~\ref{thm.algol}, so now the
manifold $M$ is atoroidal, acylindrical, and does not contain a
nonorientable essential surface.

\begin{proof}[Proof of Theorem~\ref{thm.algol}]
  As input, we are given a triangulation $\cT$ of $M$ which will be
  finite if $M$ is closed or could be finite or ideal if $M$ has
  boundary.  If $M$ has boundary and we are given a finite
  triangulation, convert it to an ideal one using the procedure
  described in the proof of \cite[Theorem 1.1.13]{Matveev2007}, which
  is relevant as per
  \cite[Proposition~3]{JacoRubinsteinSpreerTillmann2018}.  When $M$
  has boundary, apply the algorithm of \cite[Theorem
  4.7]{JacoRubinstein2011} so that the ideal triangulation $\cT$ we
  are working with is \bdry efficient.
  
  Start by computing $\CLW$ via Theorem~\ref{thm.lwfaces}. As in the
  proof of Theorem~\ref{thm.main}, the claim that we can compute the
  overall generating function algorithmically follows if we can
  implement Theorem~\ref{thm.oneface} or
  Theorem~\ref{thm.idealoneface} as appropriate for a particular face
  $C$ of $\CLW$.  (In the case when $M$ has boundary, by
  Lemma~\ref{lem:allornone} we can skip any $C$ which carries a vertex
  link, which is easy to test.)  From the proofs of those theorems, we
  need algorithms for two things: finding the subspace $W_C$ from
  Theorem~\ref{thm.compiso} and applying Theorem~1.1 of \cite{N-Pak}.
  The latter is provided by \cite{N-Pak} itself, so we focus on the
  former.

  To compute $W_C$, first let $F_0$ be any normal surface carried by
  the interior of $C$ and compute $\cG_\cT^{\leq  \weight(F_0)}$.  Now
  apply Theorem~\ref{thm.path} to find all least-weight surfaces
  $F_1, \ldots, F_k$ that are isotopic to $F$. As $C$ is complete, all
  the $F_i$ are carried by $C$ and hence are mutually compatible. By
  Lemma~\ref{thm.compiso}, we have that $W_C$ is spanned by the
  $\vF_0 - \vF_k$, giving us the needed description of $W_C$.

  The second claim of the theorem, that we can give unique normal
  representatives of the isotopy classes of incompressible surfaces
  with $\chi = -2n$, is easy by looking at the lattice points in the
  sublevel sets of $\chi$ on $\Ctil \subset \cS_\cT$ for each face $C$
  of $\CLW$ and modding out by $W_C$.  The final claim, that we can
  determine the isotopy class of a given incompressible surface $F$,
  can be done using this list and $\cG_\cT^{\leq \weight(F)}$ because
  of Theorem~\ref{thm.path}.
\end{proof}

\section{Almost normal surfaces in ideal triangulations}
\label{sec.ideal}

This section discusses the algorithm used to implement
Theorem~\ref{thm.algol} for the computations in
Section~\ref{sec.computations}.  The key difference compared to the
proof of Theorem~\ref{thm.algol} in Section~\ref{sec:decision} is that
we use almost normal surfaces, rather than normal ones, to determine
which normal surfaces are incompressible and to find isotopies between
them.

In this section, we study manifolds with boundary a union of tori
using ideal triangulations admitting a \emph{partially flat angle
  structure} in the sense of \cite{Lackenby:heeg}.  Since we are
restricting to $M$ with $\chi(\partial M) = 0$, we require the angles
at each ideal vertex to sum to exactly $\pi$ rather than at most $\pi$
as in (i) on page 916 of \cite{Lackenby:heeg}.  Such triangulations
impose restrictions on the topology of the underlying manifold
$M$. The only connected closed normal surfaces in $\cT$ with
$\chi \geq 0$ are vertex links, and $M$ is irreducible, \bdry
irreducible, atoroidal, and acylindrical
\cite[Theorem~2.2]{Lackenby:heeg}; in particular, the interior of $M$
admits a finite-volume complete hyperbolic metric, and $\cT$ is a
$\partial$-efficient triangulation of $M$.  A general algorithm for
finding such a $\cT$ in this setting is given in
\cite[\sec2]{Lackenby:heeg} and in practice one easily finds a $\cT$
admitting the stronger notion of a \emph{strict angle structure}
from~\cite{HRS:strict}.

\subsection{Tightening almost normal surfaces}
\label{sec:tight}

An \emph{almost normal surface} is a surface $S$ in $\cT$ built from
the same elementary discs as normal surfaces except for exactly one
piece, which is either an almost normal octagon or made by joining
two elementary discs in the same tetrahedron by an unknotted tube.
Given a transverse orientation of an orientable almost normal surface
$A$ in $\cT$, we can ``destabilize'' the exceptional piece in that
direction and then perform normalization moves.  This process, called
\emph{tightening} the surface $A$, moves it in only one direction
and terminates in a normal surface which we denote $\tight_{+}(A)$,
see \cite[Chapter~ 4]{Schleimer:thesis} for details.  While the
sequence of normalization moves is not unique, the tightened surface
$\tight_{+}(A)$ is well-defined: together $A$ and $\tight_{+}(A)$
bound the \emph{canonical compression body} of $A$ defined in
\cite[\sec 4.1]{Schleimer:thesis} which we denote $V_{+}(A)$.  Here,
the compression body $V_{+}(A)$ is built from $A \times I$ by adding
2- and \3-handles, so that $A$ is the minus boundary of $V_{+}(A)$ and
$\tight_{+}(A)$ is the plus boundary.  In particular, we have
$\chi(A) \leq \chi\left(\tight_{+}(A)\right)$ with equality if and
only if $V_{+}(A)$ is just $A \times I$.  Moreover, since $\cT$
contains no normal 2-spheres, every component of $\tight_{+}(A)$ has
genus at least 1.  We will use $\tight_{-}(A)$ to denote the
tightening of $A$ in the opposite transverse direction with $V_{-}(A)$
the corresponding canonical compression body.

As per \cite{Schleimer:thesis}, the tightening process can be followed
by tracking just the intersection of each surface with the \2-skeleton
of $\cT$.  Thus it amounts to looking at a family of arcs in $\cT^2$,
which need not all be normal, and then doing a sequence of bigon moves
across edges of $\cT^1$ until one is left only with normal arcs.  It
is thus straightforward to implement once one creates an appropriate
data structure to do the bookkeeping, though our code is the first
time this has been done.

\subsection{Sweepouts and thin position}

The notions of sweepouts \cite{Ru} and Gabai's thin position
\cite{Thompson} can independently be used to prove the existence of
almost normal surfaces in many situations.  We will need the following
two such results, which are quite standard.

\begin{lemma}\label{lem:prod}
  Suppose $N$ and $N'$ are normal surfaces cobounding a product
  region $V$.  Then there are disjoint surfaces
  $N = N_0, A_1, N_1, A_2, \ldots, N_{n - 1}, A_n, N_n = N'$ in $V$
  with the $N_k$ normal and the $A_k$ almost normals such that
  $\tight_{-}(A_k) = N_{k-1}$ and $\tight_{+}(A_k) = N_k$.
\end{lemma}
\begin{proof}
  Since $\cT$ is ideal and the surfaces $N$ and $N'$ are closed, the
  product region $V$ between them contains no vertices of $\cT$.  The
  usual sweepout or thin position argument for a product, see e.g.
  \cite[Theorem~6.2.2]{Schleimer:thesis}, gives the needed sequence of
  surfaces.
\end{proof}

\begin{lemma}\label{lem:comp}
  Suppose $V$ is a nonproduct compression body in $M$ where both
  $\partial_{-}V$ and $\partial_{+} V$ are normal surfaces in $\cT$.
  Then there exists an almost normal surface $A \subset V$ such that
  $\tight_{-}(A)$ and $A$ are parallel inside $V$ to $\partial_{-}V$
  and $\tight_{+}(A)$ is a (proper) compression of $\partial_{-}V$.
\end{lemma}
\begin{proof}
  As $\cT$ is ideal, there are no vertices of $\cT$ inside $V$, and a
  push off of $\partial_{-}V$ into $V$ is a strongly irreducible
  Heegaard surface for $V$.  By a slight strengthening of \cite{Ru,
    Stocking} in the same manner as \cite[Theorem 4.2]{Lackenby:heeg},
  we can find an almost normal surface $A$ in $V$ that is parallel to
  $\partial_{-}V$.  Transversely orient $A$ away from $\partial_{-}V$.
  Since $A$ is incompressible in the negative direction, we have
  $\tight_{-}(A)$ is parallel to $\partial_{-}V$.  We are done if
  $\tight_{+}(A)$ is a compression of $A$.  Otherwise, the region
  between $\partial_{-}V$ and $\tight_{+}(A)$ is a product, and we
  repeat the argument on the compression body bounded by the normal
  surfaces $\tight_{+}(A)$ and $\partial_{+} V$.  As there is a bound
  on the number of disjoint normal surfaces in $\cT$, none of which are
  normally isotopic, this will terminate and so produce the surface we
  seek.
\end{proof}

\subsection{Finiteness of (almost) normal surfaces}

For $g \geq 2$, we define $\NT^g$ to be the set of connected normal
surfaces of genus $g$ in $\cT$, up to normal isotopy. Correspondingly,
the set of such almost normal surfaces is $\AT^g$, again up to normal
isotopy.  A key result for us is Theorem 4.3 of \cite{Lackenby:heeg}
\begin{theorem}[\cite{Lackenby:heeg}]
  \label{thm:NTATfinite}
  When $\cT$ has a partially flat angle structure, both $\NT^g$ and
  $\AT^g$ are finite and algorithmically computable.
\end{theorem}
We will sketch the proof of Theorem~\ref{thm:NTATfinite} as it
outlines the algorithm for finding $\NT^g$ and $\AT^g$, which is an
important component of the overall algorithm given in
Section~\ref{sec:actualalgorithm}.  This discussion is most natural in
the setting of the \emph{quadrilateral coordinates} for normal
surfaces introduced in \cite{To} rather than the standard
triangle-quad coordinates we've used so far.  We now describe
the basics of quad coordinates, referring to \cite{Burton2009} for
details.  As the name suggests, in these coordinates a normal surface
$F$ is recorded by just the $3t$ weights on the quadrilateral discs,
where $t$ is the number of tetrahedra of $\cT$.  It turns out this
determines $F$ up to any vertex-linking components $H_v$ as in
Section~\ref{sec:ideal}. There are still linear equations, one for
each edge of $\cT$, characterizing the admissible vectors in $\N^{3t}$
that give normal surfaces; we use $\ST'$ and $\PT'$ to denote the
corresponding linear solution space and its intersection with the
positive orthant.  (The relationship between the vertices of $\PT$ and
$\PT'$ is described in detail in \cite{Burton2009}.)  Given an
admissible vector $\vv \in \N^{3t}$ carried by $\PT'$, we take the
associated normal surface to be the one with those quad weights and no
vertex-linking components; following \cite{Burton2009}, we call such
surfaces \emph{canonical}. Two things to keep in mind about quad
coordinates:
\begin{enumerate}
\item In standard coordinates, adding vector representatives
  corresponds to the geometric Haken sum. In quad coordinates, adding
  vector representatives corresponds to geometric Haken sum
  \emph{followed by removing all copies of the vertex links}.  Hence
  the total weight is additive in standard coordinates but only
  subadditive in quad coordinates.  Correspondingly, the total weight
  of a surface is only piecewise linear in quad coordinates.

\item\label{item:proper} Because we have an angle structure on $\cT$,
  the Euler characteristic function is linear in quad coordinates.
  (In contrast, it is only piecewise linear in quad coordinates for
  finite triangulations.)  Specifically, consider the linear function
  $\chi \maps \R^{3t} \to \R$ defined as follows.  Consider the
  basis vector $e_i$ corresponding to a quad $Q$ in a
  tetrahedron $\sigma$.  We set
  \[
  \chi(e_i) = -1 + \frac{\theta_1 + \theta_2 + \theta_3 + \theta_4}{2\pi}
  \]
  where the $\theta_k$ are the angles assigned to the four edges of
  $\sigma$ that $Q$ meets.  Then $\chi(F) = \chi(\vF)$ by \cite[Proposition
  4.3]{Lackenby2000}.  Moreover, for each face $C$ of $\PT'$ that
  carries only admissible vectors, the function
  $\chi \maps (\R_+ \cdot C) \to \R$ is in fact proper, nonpositive,
  and zero only at the origin \cite[Theorem~2.1]{Lackenby:heeg}; when
  the angle structure is strict, this is immediate for all of $\PT'$
  since each $\chi(e_i) < 0$.
\end{enumerate}

Turning to almost normal surfaces, those with octagons can be
described in terms of lattice points in certain polytopes, and we will
use the quad-octagon coordinates of \cite{Burton2010}, as opposed to
the standard quad-octagon-tri coordinates, to record them.  Almost
normal surfaces with tubes will be encoded by a \emph{normal} surface
together with the pair of adjacent normal discs that the tube runs
between.  With these preliminaries in hand, we can now give:
\begin{proof}[Proof of Theorem~\ref{thm:NTATfinite}]
  First, consider the case of $\NT^g$ which is contained in the
  preimage $\chi^{-1}(2 - 2g)$ for the map $\chi \maps \R^{3t} \to \R$
  defined in (\ref{item:proper}) above.  Since $\chi$ is proper on
  the cone over each admissible face of $\PT'$, there are only
  finitely many lattice points in $\chi^{-1}(2 - 2g)$ corresponding to
  surfaces.  These can be enumerated and tested for whether the
  surfaces are connected, giving us exactly $\NT^g$.

  For $\AT^g$, we consider the cases of octagons and tubes
  separately. For octagons, the map $\chi$ is again proper on the
  relevant polytope, and so this case works out the same as $\NT^g$.
  For tubes, one first enumerates all normal surfaces (not necessarily
  connected) with $\chi = 4 - 2g$.  For each such surface $N$, one
  considers all possible tubes and selects those that produce
  a connected surface, i.e.~an element of $\AT^g$.
\end{proof}

\begin{remark}
  A single normal surface $F$ can give rise to many different almost
  normal surfaces with tubes, where we are considering almost normal
  surfaces up to normal isotopy.  To keep the computation manageable,
  we considered non-normal isotopies of tubes for a fixed normal
  surface $F$.  That is, for an almost normal surface $A$ made by
  adding a tube to $F$, we can ``slide'' the attaching points of the
  tube through one of the faces of the tetrahedron that contains it to
  get another almost normal surface built on the same $F$. In our
  actual computations, we considered such surfaces up to this
  equivalence.  It is not hard to show that two tubed surfaces that
  are equivalent in this sense have the same canonical
  compression body and hence the same tightenings.
\end{remark}

\subsection{Another graph of normal surfaces}
\label{sec:NTgraph}

We now turn $\NT^g$ into a graph by adding edges as follows.  For each
$A \in \AT^g$, we pick a transverse orientation arbitrarily and
consider its two tightenings $\tight_{\pm}(A)$:
\begin{enumerate}
  \item If both $\tight_{\pm}(A)$ are homeomorphic to $A$, we
    add an (undirected) edge joining  $\tight_{-}(A)$ and
    $\tight_{+}(A)$.  In this situation, both $V_{\pm}(A)$ are
    products and so $\tight_{-}(A)$ and $\tight_{+}(A)$ are isotopic.

  \item If $\tight_{+}(A)$ is homeomorphic to $A$ but $\tight_{-}(A)$
    is not, mark the vertex $\tight_{+}(A)$ in $\NT^g$ as
    compressible.  Also do the same with the roles of $\tight_{+}(A)$
    and $\tight_{-}(A)$ reversed.
\end{enumerate}
The main result of this subsection is:
\begin{theorem}\label{thm.graphideal}
  Isotopy classes of closed essential surfaces in $M$ of genus $g$ are
  in bijection with the connected components of $\NT^g$ where no
  surface is marked as compressible.
\end{theorem}
Given how the edges in $\NT^g$ were defined, to prove
Theorem~\ref{thm.graphideal}, it suffices to show the following two
lemmas:
\begin{lemma}
  If $N \in \NT^g$ is essential and isotopic to $N' \in \NT^g$ then
  there is a path joining them in $\NT^g$.
\end{lemma}
\begin{proof}
  Let $w = \max\big(\weight(N),\, \weight(N')\big)$ and consider the
  graph $\cG_\cT^{\leq w}$ from Section~\ref{sec:graphincomp}.
  Temporarily viewing $N$ and $N'$ as vertices of $\cG_\cT^{\leq w}$,
  since they are isotopic surfaces, Theorem~\ref{thm.path} gives a
  sequence of normal surfaces $N = N_0, N_1, \ldots, N_n = N'$ where
  $N_k$ and $N_{k+1}$ can be normally isotoped to be disjoint and
  cobound a product region $P_k$.  Applying Lemma~\ref{lem:prod} to
  $P_k$ gives a path in $\NT^g$ joining $N_k$ to $N_{k+1}$.
  Concatenating these paths together gives a path in $\NT^g$ joining
  $N$ to $N'$ as needed.
\end{proof}

\begin{lemma}
  If $N \in \NT^g$ is compressible, it can be joined by a path in
  $\NT^g$ to a surface $N'$ that is marked as compressible.
\end{lemma}
\begin{proof}
  We first show there exists a nontrivial compression body $V$ in $M$
  with $\partial_-V = N$ and $\partial_+V$ a normal surface.
  Splitting $M$ open along $N$ and using the characteristic
  compression body of \cite[\sec 2]{Bonahon1983}, we can find a
  nontrivial compression body $V \subset M$ with $\partial_-V = N$ and
  $\partial_+V$ is incompressible in the complement of $N$.  Since $N$
  is normal, barrier theory \cite[Theorem~3.2(1)]{JacoRubinstein2003}
  tells us that we can normalize $\partial_+V$ in the complement of
  $N$ via an isotopy, giving us the desired compression body.

  By Lemma~\ref{lem:comp}, there is an almost normal surface $A$ in
  $V$ so that $\tight_{-}(A)$ and $A$ are parallel to $N$ and
  $\tight_{+}(A)$ is a compression of $A$.  Set $N' = \tight_{-}(A)$,
  which is marked as compressible because of the surface $A$.  As $N$
  and $N'$ are parallel, by Lemma~\ref{lem:prod} they are joined by a
  path in $\NT^g$, completing the proof of the lemma.
\end{proof}

\subsection{Algorithm}\label{sec:actualalgorithm}

The input for the algorithm is an ideal triangulation $\cT$ with a
partially flat angle structure of a manifold $M$ where
$H_2(\partial M; \F_2) \to H_2(M; \F_2)$ is onto so $M$ contains no
nonorientable surfaces by Proposition~\ref{prop:nonorient}. (Here, the
angles are given as rational multiples of $\pi$; the set of partially
flat angle structures form a convex polytope with rational vertices,
so this is not a real restriction.)  The output
is the list $\{(C,W_C)\}$ of complete essential lw-faces $C$ of
$\LW \subset \PT$ together with the corresponding subspaces $W_C$.
Before starting, recall that a \emph{vertex surface} is a normal
surface $F$ where $\vF$ is a primitive lattice point on the ray
corresponding to an admissible vertex of $\PT$.  Note that every
vertex surface is connected, and that each vertex-linking torus $H_v$
for $v \in \cT^0$ is a vertex surface \cite[Corollary
4.4]{Burton2009}.

\begin{enumerate}[label={(\arabic*)}, ref={\arabic*}]
\item \label{step:start} Enumerate all vertex surfaces for the normal
  surface equations for $\cT$ in standard triangle-quad coordinates
  via \cite[Algorithm 5.17]{Burton2009}. Then use Algorithm~3.2 of
  \cite{Burton2014} to find all admissible faces of $\PT$.  Set $g_0$
  to be the maximum genus of any vertex surface.

\item \label{step:NTAT}
  For each $g$ with $2 \leq g \leq g_0$, enumerate the finite sets
  $\NT^g$ and $\AT^g$ as described in the proof of
  Theorem~\ref{thm:NTATfinite}. Apply the tightening procedure of
  Section~\ref{sec:tight} to each surface $A$ in $\AT^g$ to compute
  $T_{-}(A)$ and $T_{+}(A)$.  As detailed in
  Section~\ref{sec:NTgraph}, this information makes $\NT^g$ into a
  graph where certain vertices are labeled compressible.

  We then compute a complete list of lw-surfaces of genus $g$ from the
  graph $\NT^g$ using Theorem~\ref{thm.graphideal} as follows: for
  each connected component of $\NT^g$ where no surface was marked as
  compressible, compute the weight of each surface and then take all
  those of minimal weight for that component.  This also computes all
  isotopy relations among the lw-surfaces of genus $g$.  Because of
  the angle structure on $M$, there are no essential tori in $M$ and
  the only nonessential lw-surfaces are the vertex links.  Thus we now
  have a complete list of all essential lw-surfaces of genus at most
  $g_0$.

\item \label{step:lwfaces} From the list of admissible faces of $\PT$
  which were computed in Step~\ref{step:start}, select those where all
  vertices are among the essential lw-surfaces enumerated in
  Step~\ref{step:NTAT}.  For each such $C$, select a surface $F_C$
  carried by its interior, e.g.~take $F_C$ to be the sum of the vertex
  surfaces of $C$. Use Algorithm~9.4 of \cite{JacoTollefson1995} to
  decompose $F_C$ into its connected components which are again normal
  surfaces.  If any component of $F_C$ has genus greater than $g_0$,
  replace $g_0$ with the maximum genus of any component of $F_C$
  and re-run Step~\ref{step:NTAT}.  By Theorem~\ref{thm.interior}, the
  face $C$ is least-weight if and only if every connected component of
  $F_{C}$ is a lw-surface.  So we now determine whether or not $C$ is
  least-weight by using the list of lw-surfaces of genus at most $g_0$
  from Step \ref{step:NTAT}.  Finally, if $C$ is least-weight, it is
  essential by Lemma~\ref{lem:allornone} and the observation that if
  $C$ carried a vertex link $H_v$, then $H_v$ would have to be one of
  the vertex surfaces of $C$.  We now have a complete list of all
  essential lw-faces of $\PT$.

\item Now we determine which essential lw-faces are complete.  Given
  such a face $C$, let $F_C$ be the preferred surface in its interior.
  By Step \ref{step:NTAT}, we know all lw-surfaces isotopic to a
  connected component of $F_C$.  By Theorem~\ref{thm.interior}, the
  face $C$ is complete if and only if all these other surfaces are
  carried by $C$.

\item \label{step:WC}
  It remains to determine the subspace $W_C$ for each complete
  essential lw-face $C$.  Let $G_1, \ldots, G_k$ be all lw-surfaces
  isotopic to a connected component of $F_{C}$, all of which will be
  carried by $C$ as it is complete.  By the last part of
  Corollary~\ref{cor:WC}, the vectors $\vG_i - \vG_j$ span $W_C$.
  This concludes the algorithm.
\end{enumerate}

\begin{remark}
  It is clearly to our advantage to keep $g_0$ as small as possible,
  which suggests several performance improvements.  For example, in
  Step~\ref{step:lwfaces} it pays to search the interior of $C$ for an
  $F_C$ whose components have the least genus.  More elaborately,
  say that a normal surface $N$ has an \emph{obvious compression}
  when there is a chain of quads forming an annulus around a thin
  edge.  In our setting, such surfaces cannot be essential, and we can
  discard in Step~\ref{step:start} any vertex surfaces with obvious
  compressions before setting $g_0$.  Because the notion of obvious
  compression can be framed as an admissibility criteria on the faces
  of $\PT$ that is compatible with \cite [Algorithm 3.2]{Burton2014},
  it is not hard to check that this does not affect the correctness of
  the answer.
\end{remark}

\begin{remark}
  The above algorithm in particular determines whether or not $M$
  contains a closed essential surface.  It would be very interesting
  to study the practical efficiency of this algorithm as compared to
  the more traditional approach of \cite{BurtonCowardTillmann2013,
  BurtonTillmann2018} involving testing for incompressibility by
  cutting $M$ open along candidate surfaces.  While we did implement
  Steps~\ref{step:start}--\ref{step:lwfaces} of our algorithm for the
  computations in Section~\ref{sec.computations}, we used a high-level
  but slow programming language and did not optimize the code
  extensively.  Consequently, we did not have a good basis for making
  this comparison.
\end{remark}
\section{Computations, examples, and patterns}
\label{sec.computations}

\begin{table}
\begin{center}
  \tablefont
  \newcommand{\phzero}{\phantom{0}}
  \begin{tabular}{lrrrrcc}
    \toprule
    \thead{sample} & \thead{count} & \thead{small} & \thead{barely \\ large}
                   & \thead{very \\ large}
                   & \thead{isotopy \\ of lw} & \thead{failed} \\
    \midrule
    Cusped census & 44,692 & 38,358 & 6,046 & 288 & \phzero 0 & \phzero 0 \\
    Knot exteriors & 14,656 & 10,554 & 3 & 4,049 & 14  & 36 \\
    \midrule
    Combined & 59,132 & 48,703 & 6,049 & 4,330 & 14  & 36  \\
    \bottomrule
  \end{tabular}
\end{center}
\caption{Summary of the manifolds where we tried to compute $\LW$.
  There are 216 manifolds common to both samples, which is why the
  last row is not the sum of the previous two.  The 14 knot
  exteriors where there was a non-normal isotopy of lw-surfaces are all
  large; they are likely very large, but we did not check this.
}
\label{tab:overall}
\end{table}

In this section, we describe the results of computing $\LW$ for some
59,096 manifolds with torus boundary.  These manifolds were drawn from
two censuses.  The first was the 44,692 orientable hyperbolic
3-manifolds that have ideal triangulations with at most 9 tetrahedra
where $\partial M$ is a single torus and $H_1(M; \F_2) = 0$
\cite{Burton2014census}.  The second was the 14,656 hyperbolic knots
in $S^3$ with at most 15 crossings whose exteriors have ideal
triangulations with at most 17 ideal tetrahedra
\cite{HosteThistlethwaiteWeeks1998}.  These two censuses have little
overlap, with only 216 manifolds common to both.

Using the default triangulation for each manifold provided by SnapPy
\cite{SnapPy}, we used Algorithm~\ref{sec:actualalgorithm} to try to
compute $\LW$.  We succeeded except for 36 triangulations where the
computation ran out of time or memory.  Combined, the
computations took about 8 CPU-months, with the running time for a
single manifold having a mean of 5.85 minutes and a maximum of 3
days.  The median time was 0.8 seconds for the cusped census and 1.3
minutes for the knot exteriors.

With the initial triangulations, some 182 manifolds had distinct
essential lw-surfaces that were isotopic.  To avoid computing the
generating function $B_M(x)$ in the general case where one is taking
the quotient by the subspaces $W_C$, we replaced 168 of these
triangulations with others where there were no such isotopies.  Except
for Section~\ref{sec:isoexs}, we will unfairly lump the remaining 14
manifolds with distinct isotopic essential lw-surfaces in with the 36
whose computations timed out, and restrict our analysis to the other
59,082.  A summary of these manifolds is given in
Table~\ref{tab:overall}, where we use the following terminology.
Recall a 3-manifold $M$ is \emph{large} when it contains a closed
essential surface and \emph{small} otherwise.  We call a large
manifold $M$ \emph{barely large} when every closed essential surface
is a multiple of a finite collection of such surfaces; otherwise $M$
will be \emph{very large}.  In the language of Section~\ref{sec:ML},
the terms small, barely large, and very large correspond,
respectively, to $\ML_0(M) = \emptyset$, $\dim( \ML_0(M)) = 1 $, and
$\dim(\ML_0(M)) > 1$.

It is natural to ask what is the smallest volume of a hyperbolic
manifold that is barely or very large.  In our sample, the smallest
manifold that is barely large is $m137$, which has volume
$V_{\mathit{oct}} \approx 3.663862376$, and the smallest knot exterior
that is barely large is that of $K15n153789$, which has volume about
$9.077985047$.  Similarly, the smallest manifold we found that is very
large is $s783$ which has volume about $5.333489566$, and the smallest
such knot exterior is that of $K10n10 = 10_{153}$, which has volume about
$7.374343889$.

\begin{table}
\begin{center}
  \tablefont
  \newcolumntype{i}{@{\hspace{0.8em}}r}
  \newcolumntype{s}{p{0em}}
  \newcommand{\muc}{\multicolumn{1}{c}{$\mu$}}	
  \begin{tabular}{rr*{4}{sri}}
    \toprule
    \multirowthead{2}{dim}
    & \multirowthead{2}{count} 
    & & \multicolumn{2}{c}{comps}
    & & \multicolumn{2}{c}{verts}
    & & \multicolumn{2}{c}{max faces}
    & & \multicolumn{2}{c}{face size} \\
    & & & \muc & range & & \muc & range & & \muc & range & & \muc & range\\
    \midrule 
    1 & 1,697 & &  1.1 & $[1, 4]$ & &  2.8 & $[2,  14]$ & &  1.8 & $[1, 10]$ & &    2 & $[2,  2]$\\
    2 & 1,810 & &  1.1 & $[1, 2]$ & &  5.6 & $[3,  16]$ & &  3.3 & $[1, 13]$ & &  3.1 & $[3,  4]$\\
    3 &   606 & &  1.2 & $[1, 3]$ & & 10.6 & $[6,  21]$ & &  7.8 & $[1, 26]$ & &  4.8 & $[4,  7]$\\
    4 &   205 & &  1.0 & $[1, 2]$ & & 16.4 & $[8,  44]$ & & 11.3 & $[2, 48]$ & &  7.7 & $[5, 12]$\\
    5 &    12 & &  1.0 & $[1, 1]$ & & 19.3 & $[16, 21]$ & & 13.6 & $[7, 18]$ & & 10.8 & $[10, 12]$\\
    \bottomrule
  \end{tabular}
\end{center}

\caption{
Statistics about the complexes $\LW$ for the 4,330 very large
manifolds, broken down by $\dim \LW$. The properties recorded are: the
number of connected components (comps), the number of vertices
(verts), the number of maximal faces (max faces), and the largest
number of vertices in any face (face size). For each numerical
property, we give the mean in the $\mu$ column as well as
the min-max interval in the range column.
}
\label{tab:stats}
\end{table}

\subsection{Very large manifolds}

For the 4,330 manifolds where $\dim \LW \geq 1$, the complexes $\LW$
run the gamut from a single edge (for 760 manifolds) up to monsters
like $\LW$ for $K13n3838$ which is connected with 44 vertices and
48 maximal faces all of dimension 4, where each maximal face has
between 5 and 9 vertices.  Basic statistics about the topology and
combinatorics of the $\LW$ are given in Table~\ref{tab:stats}.  All
but 178 of these complexes are \emph{pure}, that is, every maximal
face has the same dimension; the exceptions are 140 cases where each
component of $\LW$ is pure but there are components of differing
dimensions, and 38 cases where $\LW$ is connected and impure.

While the combinatorics of some of these complexes is quite elaborate,
the underlying topology of all $\LW$ in our sample is simple in that
every connected component is actually contractible.  Moreover, for a
component $Y$ of dimension $d$, each $(d - 1)$--face is glued to at
most two $d$-faces; consequently, all components of dimension 1 are
homeomorphic to intervals rather than more general trees.  Here,
contractibility was checked as follows.  First, each $\LW$ was
converted to a simplicial complex (some 3,603 of the $\LW$ are in fact
simplicial, for the rest a barycentric subdivision of the polyhedral
complex was used).  We then checked that every component had vanishing
reduced homology and trivial fundamental group using \cite{SageMath},
which implies contractibility.

\begin{table}
\begin{center}
  \tablefont
  \newcommand{\phz}{\phantom{0}}
  \begin{tabular}{ccllr}
    \toprule
    dim & count & \multicolumn{1}{c}{degree} & periods & $\ell_1$-norm \\
    \midrule 
    1 &     35  &  2, 3, 4, 6, 8  &  1, 2, 3, 6  &  [7, 53] \\
    2 &     18  &  3, 6, 7        &  1, 2, 3     & [16, 38] \\
    3 &     24  &  4, 5, 6, 7, 8  &  1, 2        & [26, 71] \\
    4 & \phz 9  &  5, 6, 8        &  1, 2        & [38, 94] \\
    5 & \phz 2  &  7, 8           &  2           & [78, 88] \\
    \bottomrule
  \end{tabular}
\end{center}

\caption{Statistics about the 88 distinct generating functions
  $B_M(x)$ for the 4,330 very large manifolds, broken down by
  $\dim \LW$. Here, each $B_M(x)$ has rational form $P(x)/Q(x)$ for
  some $P, Q \in \Z[x]$ with $\deg P = \deg Q$. The properties
  recorded are: the number of distinct $B_M(x)$ (count), the values of $\deg P$
  (degree), the observed periods of $B_M(x)$ (periods), and the range of the
  $\ell_1$-norm of the combined coefficients of the polynomials $P$
  and $Q$ ($\ell_1$-norm). }
\label{tab:gen}
\end{table}

\begin{table}
  \centering
  \tablefont
  \newcommand{\showBM}[2]{$\displaystyle \frac{#1}{#2}$}
  \newcommand{\BMlinesp}{\addlinespace[0.6em]}
  \begin{tabular}{clcccr}
    \toprule
    dim & \multicolumn{1}{c}{$B_M(x)$} & per & $\ell_1$ & sample $M$ & count\\
    \midrule 
    1 & \showBM{-x^2 + 2x}{(x - 1)^2}
      &  1  &  7 & $K10n10$ & 1,009 \\ \BMlinesp
    1 & \showBM{-x^4 + 2x^2}{(x - 1)^2(x + 1)^2}
      & 2 & 7 & $K14n11913$ & 259 \\ \BMlinesp
    2 & \showBM{-x^3 + 3 x^2 - 4x}{(x - 1)^3}
      & 1 & 16 & $t12766$ & 1,459 \\ \BMlinesp 
    2 & \showBM{-x^6 + 3 x^4 - 6 x^2}{(x - 1)^3(x + 1)^3}
      & 2 & 18 & $K15n93515$ & 82 \\ \BMlinesp
    3 & \showBM{-x^4 + 4x^3 - 5x^2 + 6x}{(x - 1)^4}
      & 1 & 32 & $K12n605$ & 219 \\ \BMlinesp 
    3 & \showBM{-x^5 + 3x^4 - 2x^3 + 2x^2 + 6x}{(x - 1)^4 (x + 1)}
      & 2 & 26 & $K11n34$ & 139 \\ \BMlinesp
    4 & \showBM{-x^6 + 4 x^5 - 5 x^4 - 2 x^3 - 2 x^2 - 8x}{(x - 1)^5 (x + 1)}
      & 2 & 42 & $K14n1808$ & 62 \\ \BMlinesp
    4 & \showBM{-x^5 + 5 x^4 - 10 x^3 + 10 x^2 - 8 x}{(x-1)^5}
      & 1 & 66 & $K12n214$ & 44 \\ \BMlinesp
    5 & \showBM{-x^8 + 4x^7 - 3x^6 - 4x^5 + 14x^4 + 2x^3 + 14x^2 + 10x}{(x - 1)^6 (x + 1)^2}
      & 2 & 88 & $K15n15582$ & 11\\ \BMlinesp
    \bottomrule
  \end{tabular}

\caption{ Nine of the most common $B_M(x)$, which together account for
  3,284 (75.8\%) of the 4,330 very large manifolds. The properties
  recorded are the dimension of $\LW$ (dim), the rational form
  $P(x)/Q(x)$ of $B_M(x)$, the period of $B_M(x)$ (per), the
  $\ell_1$-norm of the combined coefficients of $P$ and $Q$, an example 
  manifold with this $B_M(x)$ (sample $M$), and the number of
  manifolds with this $B_M(x)$ (count). }
\label{tab:gencommon}
\end{table}

\begin{table}
  \centering
  \tablefont
  \newcommand{\showBM}[2]{$\displaystyle \frac{#1}{#2}$}
  \newcommand{\BMlinesp}{\addlinespace[0.6em]}

  \begin{tabular}{clccc}
    \toprule
    dim & \multicolumn{1}{c}{$B_M(x)$} & per & $\ell_1$ & $M$  \\
    \midrule 
    1 & \showBM{-2x^8 - 4x^7 - 2x^6 + 6x^5 + 13x^4 + 8x^3 + 2x^2}{(x - 1)^2  (x + 1)^2 (x^2 + x + 1)^2}
      & 6 & 53 & $K15n138922$ \\ \BMlinesp
    2 & \showBM{-2 x^6 + 5 x^4 - 4 x^3 - 15 x^2 - 4 x}{(x - 1)^3 (x + 1)^3}
      & 2 & 38 & $K15n27228$ \\ \BMlinesp
    2 & \showBM{-2 x^7 + 2 x^6 - x^5 + x^4 - 9 x^3 - 5 x^2 - 4 x}{(x - 1)^3 (x^2 + x + 1)^2}
      & 3 & 32 & $K15n86383$ \\ \BMlinesp
    3 & \showBM{-3 x^8 + 13 x^6 + 2 x^5 - 14 x^4 - 4 x^3 + 17 x^2 + 2 x}{(x - 1)^4  (x + 1)^4}
      & 2 & 71 & $K15n139871$ \\ \BMlinesp
    4 & \showBM{-2 x^8 + 4 x^7 + 4 x^6 - 14 x^5 - 12 x^4 - 6 x^3 - 22 x^2 - 8 x}{(x - 1)^5 (x + 1)^3}
      & 2 & 94 & $K13n1795$ \\ \BMlinesp
    5 & \showBM{-x^7 + 5 x^6 - 9 x^5 + 5 x^4 + 8 x^3 + 10 x}{(x - 1)^6 (x + 1)}
      & 2 & 78 & $K13n2458$ \\ \BMlinesp
    \bottomrule
  \end{tabular}

\caption{Six of the most complicated $B_M(x)$ in our sample. The
  properties recorded are the dimension of $\LW$ (dim), the rational
  form $P(x)/Q(x)$ of $B_M(x)$, the period of $B_M(x)$ (per), the
  $\ell_1$-norm of the combined coefficients of $P$ and $Q$
  ($\ell_1$), and a manifold with this generating function ($M$). }
\label{tab:genmessy}
\end{table}

\subsection{Surface counts by Euler characteristic}
\label{sec:expsurfcounts}

For each of the 4,330 very large
manifolds, we computed the generating function $B_M(x)$ from
Theorem~\ref{thm.main} starting from $\LW$ by using Normaliz
\cite{Normaliz}.  This resulted in only 88 distinct generating
functions whose properties are summarized in Table~\ref{tab:gen} and
examples of which are given in Tables~\ref{tab:gencommon} and
\ref{tab:genmessy}.

\subsection{Sample LW complexes}

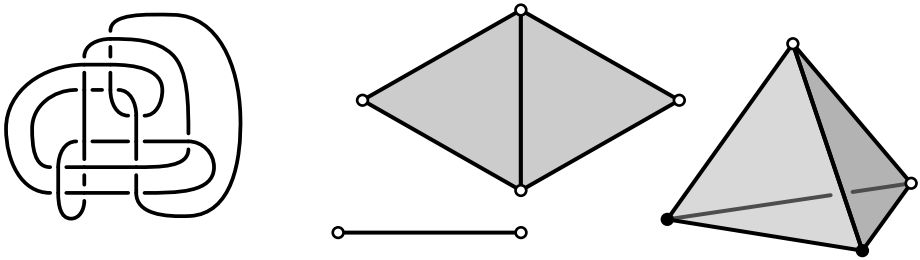
\begin{figure}[tb]
  \centering
  \begin{tikzpicture}[scale=0.8, line cap=round, line join=round]
    \begin{scope}[scale=0.4, line width=1.4pt, color=black, shift={(-13.5, -4.3)}]
    \draw (4.41, 5.48) .. controls (4.41, 4.94) and (4.66, 4.41) .. (5.14, 4.41);
    \draw (5.81, 4.41) .. controls (6.28, 4.41) and (6.54, 4.94) .. 
          (6.54, 5.47) .. controls (6.54, 6.35) and (5.43, 6.53) .. (4.41, 6.53);
    \draw (4.41, 6.53) .. controls (4.05, 6.53) and (3.69, 6.53) .. (3.34, 6.53);
    \draw (3.34, 6.53) .. controls (1.66, 6.53) and (0.13, 5.45) .. 
          (0.13, 3.87) .. controls (0.13, 2.53) and (0.75, 1.20) .. (1.94, 1.20);
    \draw (2.60, 1.20) .. controls (2.85, 1.20) and (3.09, 1.20) .. (3.34, 1.20);
    \draw (3.34, 1.20) .. controls (3.94, 1.20) and (4.54, 1.20) .. (5.14, 1.20);
    \draw (5.81, 1.20) .. controls (7.12, 1.20) and (8.66, 1.20) .. 
          (8.66, 2.27) .. controls (8.66, 2.86) and (8.19, 3.34) .. (7.61, 3.34);
    \draw (7.61, 3.34) .. controls (7.01, 3.34) and (6.41, 3.34) .. (5.81, 3.34);
    \draw (5.14, 3.34) .. controls (4.65, 3.34) and (4.16, 3.34) .. (3.67, 3.34);
    \draw (3.01, 3.34) .. controls (2.52, 3.34) and (2.27, 2.81) .. (2.27, 2.27);
    \draw (2.27, 2.27) .. controls (2.27, 1.92) and (2.27, 1.56) .. (2.27, 1.20);
    \draw (2.27, 1.20) .. controls (2.27, 0.69) and (2.36, 0.13) .. 
          (2.80, 0.13) .. controls (3.14, 0.13) and (3.34, 0.50) .. (3.34, 0.87);
    \draw (3.34, 1.54) .. controls (3.34, 1.67) and (3.34, 1.81) .. (3.34, 1.94);
    \draw (3.34, 2.60) .. controls (3.34, 2.85) and (3.34, 3.10) .. (3.34, 3.34);
    \draw (3.34, 3.34) .. controls (3.34, 4.06) and (3.34, 4.77) .. (3.34, 5.48);
    \draw (3.34, 5.48) .. controls (3.34, 5.72) and (3.34, 5.96) .. (3.34, 6.20);
    \draw (3.34, 6.86) .. controls (3.34, 7.34) and (3.87, 7.60) .. (4.41, 7.60);
    \draw (4.41, 7.60) .. controls (5.32, 7.60) and (6.29, 7.59) .. 
          (6.95, 6.95) .. controls (7.61, 6.31) and (7.61, 4.85) .. (7.61, 3.67);
    \draw (7.61, 3.01) .. controls (7.61, 2.27) and (6.45, 2.27) .. (5.47, 2.27);
    \draw (5.47, 2.27) .. controls (4.76, 2.27) and (4.05, 2.27) .. (3.34, 2.27);
    \draw (3.34, 2.27) .. controls (3.09, 2.27) and (2.85, 2.27) .. (2.60, 2.27);
    \draw (1.94, 2.27) .. controls (1.28, 2.27) and (1.20, 3.12) .. 
          (1.20, 3.88) .. controls (1.20, 4.80) and (2.04, 5.48) .. (3.01, 5.48);
    \draw (3.67, 5.48) .. controls (3.80, 5.48) and (3.94, 5.48) .. (4.07, 5.48);
    \draw (4.74, 5.48) .. controls (5.22, 5.48) and (5.47, 4.95) .. (5.47, 4.41);
    \draw (5.47, 4.41) .. controls (5.47, 4.06) and (5.47, 3.70) .. (5.47, 3.34);
    \draw (5.47, 3.34) .. controls (5.47, 3.10) and (5.47, 2.85) .. (5.47, 2.60);
    \draw (5.47, 1.94) .. controls (5.47, 1.70) and (5.47, 1.45) .. (5.47, 1.20);
    \draw (5.47, 1.20) .. controls (5.47, 0.32) and (6.62, 0.19) .. 
          (7.67, 0.24) .. controls (9.38, 0.33) and (9.80, 2.47) .. 
          (9.74, 4.44) .. controls (9.68, 6.49) and (8.84, 8.56) .. 
          (7.01, 8.60) .. controls (5.84, 8.63) and (4.41, 8.67) .. (4.41, 7.93);
    \draw (4.41, 7.27) .. controls (4.41, 7.13) and (4.41, 7.00) .. (4.41, 6.86);
    \draw (4.41, 6.20) .. controls (4.41, 5.96) and (4.41, 5.72) .. (4.41, 5.48);
  \end{scope}

  \begin{scope}[line width=1.5pt, shift={(-0.5, -1.2)}]
    \begin{scope}[shift={(0.6, -0.7)}]
      \coordinate (s) at (0, 0);
      \coordinate (t) at (3, 0);
      \draw (s) -- (t);
    \end{scope}

    \begin{scope}[shift={(1, 1.5)}]
      \coordinate (a) at (0, 0);
      \coordinate (b) at (-30:3);
      \coordinate (c) at (30:3);
      \coordinate (d) at ($(b) + (c)$);
      \filldraw[fill=black!20] (a) -- (b) -- (c) -- cycle;
      \filldraw[fill=black!20] (d) -- (b) -- (c) -- cycle;
    \end{scope}

    \begin{scope}[shift={(6, -1.0)}, scale=0.4]
      \coordinate (u) at (0,1.3);
      \coordinate (v) at (8.0, 0);
      \coordinate (w) at (10.0, 2.8);
      \coordinate (z) at (5.15, 8.6);

      \fill[fill=black!15] (u) -- (v) -- (z) -- cycle;
      \fill[fill=black!30] (v) -- (w) -- (z) -- cycle;

      \begin{scope}[color=black!70]
        \draw (u) -- ($(u)!0.67!(w)$); 
        \draw ($(u)!0.76!(w)$) -- (w); 
      \end{scope}

      \draw (u) -- (v) -- (z) -- cycle;
      \draw (v) -- (w) -- (z) -- cycle;

    \end{scope}
    
    \foreach \n in {s, t, a, b, c, d, w, z}{
      \filldraw[fill=white, line width=1pt] (\n) circle (2.5pt);
    }

    \foreach \n in {u, v}{
      \filldraw[fill=black, line width=1pt] (\n) circle (2.5pt);
    }
    
  \end{scope}
\end{tikzpicture}

  \caption{For the knot $K15n51747$ shown at left, at right is the
    complex $\LW$ for a triangulation of its exterior with 17 ideal
    tetrahedra.  This example is unusual in that there are components
    of different dimensions.  The vertex surfaces are either
    genus 2 (solid vertices) or genus 3 (open vertices).  Here
    $B_M(x) = (-3x^7 + 3x^6 + 9x^5 - 9x^4 - 9x^3 + 9x^2 +
    2x)/\big((x-1)^4(x+1)^3\big)$.}
  \label{fig:K15n51747}
\end{figure}

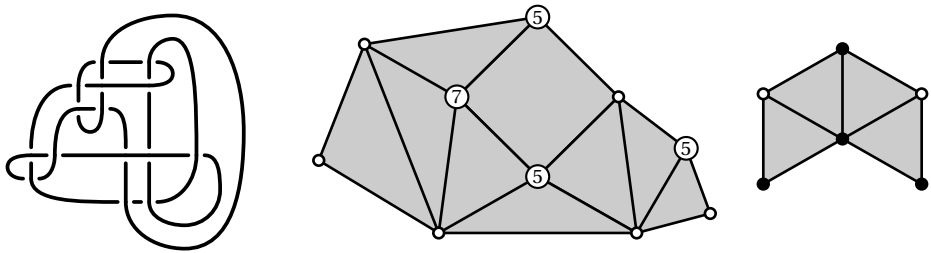
\begin{figure}[thb]
  \centering
\begin{tikzpicture}[line cap=round, line join=round, font=\scriptsize]

  \begin{scope}[line width=1pt]
    \begin{scope}[shift={(4.0, 0.5)}]
      \coordinate (N8) at (0, 0);
      \coordinate (N2) at (210:1.2);
      \coordinate (N22) at (150:1.2);
      \coordinate (N4) at (90:1.2);
      \coordinate (N21) at (30:1.2);
      \coordinate (N5) at (-30:1.2);
      
      \draw[fill=black!20] (N8) -- (N2) -- (N22) -- (N4) -- (N21) -- (N5) -- cycle;
      \draw (N8) -- (N22);
      \draw (N8) -- (N4);
      \draw (N8) -- (N21);
    \end{scope}
    
    \begin{scope}
      \coordinate (N224) at (0, 0);
      \coordinate (N601) at (135:1.5);
      \coordinate (N36) at (45:1.5);
      \coordinate (N226) at ($(N601) + (N36)$);

      \coordinate (N23) at (-30:1.5);
      \coordinate (N227) at ($(N23) + (60:1.3)$);
      \coordinate (N38) at ($(N23) + (15:1)$);
      
      \coordinate (N32) at (-150:1.5);
      \coordinate (N25) at ($(N601) + (150:1.4)$);
      \coordinate (N33) at ($(N601) + (205:2.0)$);

      \begin{scope}[fill=black!20]
        \filldraw (N224) -- (N36)  -- (N226) -- (N601) -- cycle;
        \filldraw (N224) -- (N23)  -- (N38)  -- (N227) -- (N36) -- cycle;
        \filldraw (N226) -- (N25)  -- (N33)  -- (N32) -- (N23)
                         -- (N224) -- (N601) -- cycle;
      \end{scope}
    
      \draw (N23) -- (N36); 
      \draw (N23) -- (N227); 
      \draw (N32) -- (N224);
      \draw (N32) -- (N601);
      \draw (N32) -- (N25) -- (N601);
    
    \end{scope}

    \foreach \n in {2, 4, 5, 8}{
      \filldraw (N\n) circle (2pt);
    };
    
    \foreach \n in {21, 22, 23, 25, 32, 33, 36, 38}{
      \filldraw[fill=white, line width=1pt] (N\n) circle (2pt);
    };

    \foreach \n in {224, 226, 227}{
      \node[fill=white, circle, draw, line width=0.8pt, inner sep=0.8pt] at (N\n) {5};
    };

    \node[fill=white, circle, draw, line width=0.8pt, inner sep=0.8pt] at (N601) {7};
  \end{scope}

  \begin{scope}[shift={(-7, -1)}, scale=0.32, line width=1.4pt]
    \draw (8.21, 4.01) .. controls (8.80, 4.01) and (8.85, 3.24) .. 
          (8.85, 2.56) .. controls (8.85, 1.75) and (8.20, 1.10) .. 
          (7.40, 1.10) .. controls (6.67, 1.10) and (5.94, 1.43) .. (5.94, 2.07);
    \draw (5.94, 2.07) .. controls (5.94, 2.61) and (5.94, 3.14) .. (5.94, 3.68);
    \draw (5.94, 4.34) .. controls (5.94, 5.09) and (5.94, 5.83) .. (5.94, 6.58);
    \draw (5.94, 7.24) .. controls (5.94, 7.46) and (5.94, 7.67) .. (5.94, 7.88);
    \draw (5.94, 7.88) .. controls (5.94, 8.42) and (6.38, 8.85) .. 
          (6.91, 8.85) .. controls (7.88, 8.85) and (7.88, 6.11) .. (7.88, 4.01);
    \draw (7.88, 4.01) .. controls (7.88, 2.99) and (7.23, 2.07) .. (6.27, 2.07);
    \draw (5.61, 2.07) .. controls (5.51, 2.07) and (5.41, 2.07) .. (5.31, 2.07);
    \draw (4.65, 2.07) .. controls (3.07, 2.07) and (1.10, 2.07) .. (1.10, 3.04);
    \draw (1.10, 3.04) .. controls (1.10, 3.25) and (1.10, 3.47) .. (1.10, 3.68);
    \draw (1.10, 4.34) .. controls (1.10, 5.61) and (1.59, 6.91) .. (2.71, 6.91);
    \draw (3.37, 6.91) .. controls (3.58, 6.91) and (3.80, 6.91) .. (4.01, 6.91);
    \draw (4.01, 6.91) .. controls (4.65, 6.91) and (5.30, 6.91) .. (5.94, 6.91);
    \draw (5.94, 6.91) .. controls (6.41, 6.91) and (6.91, 7.00) .. 
          (6.91, 7.40) .. controls (6.91, 7.70) and (6.60, 7.88) .. (6.27, 7.88);
    \draw (5.61, 7.88) .. controls (5.19, 7.88) and (4.76, 7.88) .. (4.34, 7.88);
    \draw (3.68, 7.88) .. controls (3.25, 7.88) and (3.04, 7.40) .. (3.04, 6.91);
    \draw (3.04, 6.91) .. controls (3.04, 6.70) and (3.04, 6.49) .. (3.04, 6.27);
    \draw (3.04, 5.61) .. controls (3.04, 5.29) and (3.22, 4.98) .. 
          (3.52, 4.98) .. controls (3.92, 4.98) and (4.01, 5.48) .. (4.01, 5.94);
    \draw (4.01, 5.94) .. controls (4.01, 6.16) and (4.01, 6.37) .. (4.01, 6.58);
    \draw (4.01, 7.24) .. controls (4.01, 7.46) and (4.01, 7.67) .. (4.01, 7.88);
    \draw (4.01, 7.88) .. controls (4.01, 9.17) and (5.46, 9.82) .. 
          (6.91, 9.82) .. controls (8.99, 9.82) and (9.82, 7.36) .. 
          (9.82, 4.98) .. controls (9.82, 2.66) and (9.39, 0.13) .. 
          (7.40, 0.13) .. controls (6.14, 0.13) and (4.98, 0.90) .. (4.98, 2.07);
    \draw (4.98, 2.07) .. controls (4.98, 2.61) and (4.98, 3.14) .. (4.98, 3.68);
    \draw (4.98, 4.34) .. controls (4.98, 5.09) and (4.98, 5.94) .. (4.34, 5.94);
    \draw (3.68, 5.94) .. controls (3.47, 5.94) and (3.25, 5.94) .. (3.04, 5.94);
    \draw (3.04, 5.94) .. controls (2.24, 5.94) and (2.07, 4.94) .. (2.07, 4.01);
    \draw (2.07, 4.01) .. controls (2.07, 3.52) and (1.86, 3.04) .. (1.43, 3.04);
    \draw (0.77, 3.04) .. controls (0.44, 3.04) and (0.13, 3.22) .. 
          (0.13, 3.52) .. controls (0.13, 3.92) and (0.64, 4.01) .. (1.10, 4.01);
    \draw (1.10, 4.01) .. controls (1.32, 4.01) and (1.53, 4.01) .. (1.74, 4.01);
    \draw (2.40, 4.01) .. controls (3.26, 4.01) and (4.12, 4.01) .. (4.98, 4.01);
    \draw (4.98, 4.01) .. controls (5.30, 4.01) and (5.62, 4.01) .. (5.94, 4.01);
    \draw (5.94, 4.01) .. controls (6.48, 4.01) and (7.02, 4.01) .. (7.55, 4.01);
  \end{scope}
  
\end{tikzpicture}

  \caption{For the knot $K15n18579$ shown at left, at right is the
    complex $\LW$ for a triangulation of its exterior with 17 ideal
    tetrahedra.  It is one of the most complicated examples in our
    sample with $\dim \LW = 2$; note that one face is a square rather
    than a triangle.  The vertex surfaces are genus 2 (solid
    vertices), genus 3 (open vertices), or genus 5 or 7 as labeled.
    Here
    $B_M(x) = (-2x^6 + 5x^4 - 4x^3 - 15x^2 - 4x)/\big((x-1)^3(x+1)^3\big)$.}
 \label{fig:K15n18579}
\end{figure}

\begin{figure}[thb]
  \centering
  \input figures/K11n34

  \caption{
    The complex $\LW$ for a triangulation of the exterior of
    the Conway knot $K11n34$.  The 11 vertex surfaces are
    $N_{5}, N_{6}, N_8, N_9, N_{11}, N_{12}, N_{18}, N_{19}, N_{37},
    N_{38}, N_{40}$, where the notation follows
    \protect\cite{PaperData}.  The first six surfaces have genus 2
    (dark vertices above) and the rest genus 3 (white vertices above).
    There are seven \3-dimensional faces: four tetrahedra, a pyramid
    with quadrilateral base ($N_5, N_6, N_{40}, N_{38}, N_{11}$), a
    triangular prism ($N_{11}, N_{12}, N_{38}, N_8, N_9, N_{37}$), and
    the one in the lower right whose faces are four triangles and two
    quadrilaterals ($N_{11}, N_{8}, N_{37}, N_{38}$ and
    $N_{11}, N_{8}, N_{18}, N_{5}$).  Here $B_M(x) = (-x^5 + 3x^4 -
    2x^3 + 2x^2 + 6x)/((x + 1)(x - 1)^4)$.
  }
\label{fig:LWconway}
\end{figure}

\begin{figure}[thb]
  \centering
  \newcommand{\LWinFourD}[2]{
\begin{tikzoverlay}[scale=0.3]{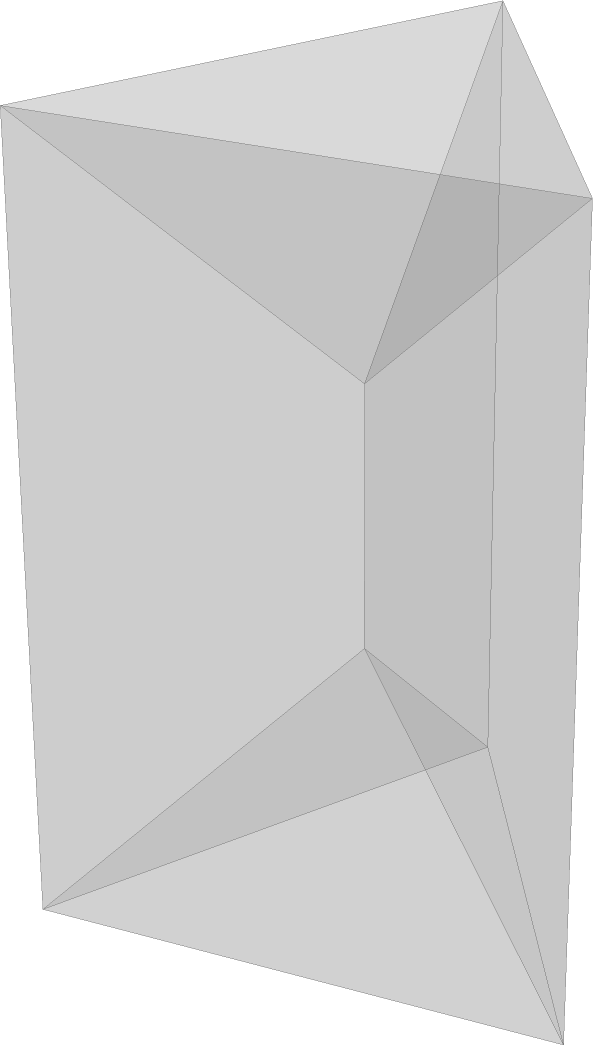}[font=\footnotesize]
  \begin{scope}[line width=1.2pt, line cap=round, line join=round]
    \coordinate (a) at (7.4,23.1);
    \coordinate (b) at (94.8,0.5);
    \coordinate (c) at (82.1,50.2);
    \coordinate (d) at (61.4,66.4);
    \coordinate (e) at (61.4,111.5);
    \coordinate (f) at (0.4,157.9);
    \coordinate (g) at (99.4,142.4);
    \coordinate (h) at (84.5,175.3);


    \begin{scope}[color=black!70]
      \draw (a) -- (d) -- (b) -- (c) -- (d) -- (e) -- (f);
      \draw (e) -- (g);
      \draw (a) -- ($(a)!0.80!(c)$);
      \draw ($(a)!0.91!(c)$) -- (c);
      \draw (e) -- ($(e)!0.48!(h)$);
      \draw ($(e)!0.61!(h)$) -- (h);
      \draw (c) -- ($(c)!0.59!(h)$);
      \draw ($(c)!0.675!(h)$) -- ($(c)!0.72!(h)$);
      \draw ($(c)!0.785!(h)$) -- (h);
    \end{scope}

    \draw (a) -- (b) -- (g) -- (f) -- cycle;
    \draw (f) -- (h) -- (g);

    \begin{scope}[line width=1.5pt]
      \foreach \n in {a, f, h}{
        \filldraw (\n) circle (2.0);
      }
      
      \foreach \n in {c, d, e}{
        \filldraw[color=black!70] (\n) circle (2.0);
      }

      \foreach \n in {b, g}{
        \filldraw[fill=white, line width=1.0pt] (\n) circle (3.0);
      };
    \end{scope}

    \node[left=3] at (a) {10};
    \node[left=3] at (f) {15};
    \node[right=3] at (b) {25};
    \node[right=3] at (g) {29};
    \node[right=3] at (h) {14};
    \node[left] at ($(d) + (-2, 3)$) {#1};
    \node[left] at ($(e) + (-2, -3)$) {#2};
    \node[right] at ($(c) + (-0.8, 0)$) {9};
  \end{scope}
\end{tikzoverlay}
}

\begin{tikzpicture}[scale=0.4, line width=1.4, line cap=round, line join=round]
  \begin{scope}[color=black]
    \draw (4.97, 4.10) .. controls (4.97, 4.55) and (5.18, 4.99) .. 
          (5.58, 4.99) .. controls (6.08, 4.99) and (6.18, 4.35) .. (6.18, 3.77);
    \draw (6.18, 3.77) .. controls (6.18, 3.48) and (6.18, 3.19) .. (6.18, 2.89);
    \draw (6.18, 2.23) .. controls (6.18, 1.35) and (4.88, 1.35) .. 
          (3.76, 1.35) .. controls (2.60, 1.35) and (1.34, 1.57) .. 
          (1.34, 2.56) .. controls (1.34, 3.23) and (1.88, 3.77) .. (2.55, 3.77);
    \draw (2.55, 3.77) .. controls (2.85, 3.77) and (3.14, 3.77) .. (3.43, 3.77);
    \draw (4.09, 3.77) .. controls (4.39, 3.77) and (4.68, 3.77) .. (4.97, 3.77);
    \draw (4.97, 3.77) .. controls (5.27, 3.77) and (5.56, 3.77) .. (5.85, 3.77);
    \draw (6.51, 3.77) .. controls (6.70, 3.77) and (6.88, 3.77) .. (7.07, 3.77);
    \draw (7.73, 3.77) .. controls (8.02, 3.77) and (8.31, 3.77) .. (8.61, 3.77);
    \draw (8.61, 3.77) .. controls (9.41, 3.77) and (9.82, 4.68) .. 
          (9.82, 5.59) .. controls (9.82, 6.46) and (9.68, 7.41) .. (8.94, 7.41);
    \draw (8.28, 7.41) .. controls (6.77, 7.41) and (5.27, 7.41) .. (3.76, 7.41);
    \draw (3.76, 7.41) .. controls (3.15, 7.41) and (2.55, 7.08) .. (2.55, 6.53);
    \draw (2.55, 5.87) .. controls (2.55, 5.28) and (2.55, 4.69) .. (2.55, 4.10);
    \draw (2.55, 3.44) .. controls (2.55, 3.00) and (2.76, 2.56) .. 
          (3.16, 2.56) .. controls (3.66, 2.56) and (3.76, 3.20) .. (3.76, 3.77);
    \draw (3.76, 3.77) .. controls (3.76, 4.58) and (3.76, 5.39) .. (3.76, 6.20);
    \draw (3.76, 6.20) .. controls (3.76, 6.49) and (3.76, 6.78) .. (3.76, 7.08);
    \draw (3.76, 7.74) .. controls (3.76, 8.62) and (5.07, 8.62) .. 
          (6.18, 8.62) .. controls (7.34, 8.62) and (8.61, 8.40) .. (8.61, 7.41);
    \draw (8.61, 7.41) .. controls (8.61, 6.31) and (8.61, 5.21) .. (8.61, 4.10);
    \draw (8.61, 3.44) .. controls (8.61, 3.15) and (8.61, 2.86) .. (8.61, 2.56);
    \draw (8.61, 2.56) .. controls (8.61, 1.98) and (8.50, 1.35) .. 
          (8.00, 1.35) .. controls (7.60, 1.35) and (7.40, 1.79) .. (7.40, 2.23);
    \draw (7.40, 2.89) .. controls (7.40, 3.19) and (7.40, 3.48) .. (7.40, 3.77);
    \draw (7.40, 3.77) .. controls (7.40, 5.30) and (5.78, 6.20) .. (4.09, 6.20);
    \draw (3.43, 6.20) .. controls (3.14, 6.20) and (2.85, 6.20) .. (2.55, 6.20);
    \draw (2.55, 6.20) .. controls (1.09, 6.20) and (0.13, 4.74) .. 
          (0.13, 3.17) .. controls (0.13, 1.07) and (2.58, 0.14) .. 
          (4.97, 0.14) .. controls (7.12, 0.14) and (9.82, 0.14) .. 
          (9.82, 1.35) .. controls (9.82, 1.97) and (9.49, 2.56) .. (8.94, 2.56);
    \draw (8.28, 2.56) .. controls (7.98, 2.56) and (7.69, 2.56) .. (7.40, 2.56);
    \draw (7.40, 2.56) .. controls (6.99, 2.56) and (6.59, 2.56) .. (6.18, 2.56);
    \draw (6.18, 2.56) .. controls (5.57, 2.56) and (4.97, 2.89) .. (4.97, 3.44);
  \end{scope}

  \node at (16, 4.5) {\LWinFourD{12}{16}};
  \node at (16, -1.5) {$\partial C_1$};

  \draw[dashed, color=black!40] (23.1, -0.5) -- +(0, 10);

  \node at (30, 4.5) {\LWinFourD{13}{17}};
  \node at (30, -1.5) {$\partial C_2$};
\end{tikzpicture}

  \caption{For a 15-tetrahedra triangulation of the exterior of the
    knot $K13n1019$ shown at left, the complex $\LW$ consists of two
    4-dimensional faces $C_1$ and $C_2$ with the same combinatorics
    that are glued together along a single \3-dimensional face.  The
    boundary of each $C_i$ is depicted above via an identification of
    $\partial C_i$ with $S^3$; hence, in each case there is an
    additional face of $\partial C_i$ on the outside, namely a
    triangular prism whose vertices are
    $N_{10}, N_{25}, N_{9}, N_{15}, N_{29}, N_{14}$.  It is these two
    outside faces that are identified to form $\LW$.  As usual, solid
    and open vertices correspond to surfaces of genus 2 and 3
    respectively, and the numbering follows \protect\cite{PaperData}.
    Here $B_M(x) = (-x^5 + 5x^4 - 10x^3 + 10x^2 - 8x)/(x - 1)^5$.  }
  \label{fig:K13n1019}
\end{figure}

We next give several examples of $\LW$ for specific triangulations
of knot exteriors.  To start off, Figure~\ref{fig:K15n51747} gives an
example of a simple $\LW$ which is unusual in having components of
different dimensions.  Then Figure~\ref{fig:K15n18579} describes one of the
most complicated \2-dimensional examples we found.
Figure~\ref{fig:LWconway} shows the fairly complicated \3-dimensional
complex coming from the Conway knot $K11n34$.  In dimension four, we were
only able to visualize one of the very simplest examples in
Figure~\ref{fig:K13n1019}, and for dimension five we simply gave up.

All 38 examples where $\LW$ is connected and impure have dimension 4
or 5; those of dimension 4 also have maximal faces of dimension 2 and
those of dimension 5 also have maximal faces of dimension 3.  One of
the simplest such is $K13n857$ where $\LW$ consists of seven
4-simplices plus two triangles, where the triangles are glued together
to form a square, and then one edge of that square is glued to the
main mass of 4-simplices.

\subsection{Isotopies of lw-surfaces}
\label{sec:isoexs}

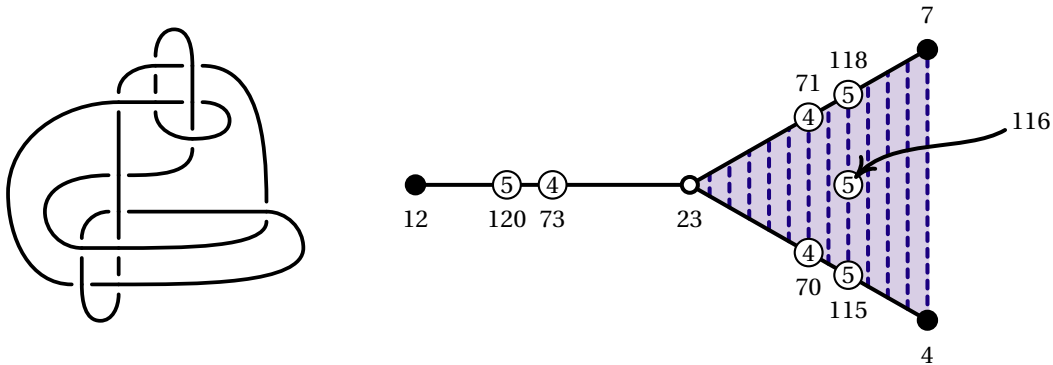
\begin{figure}
  \centering
  \begin{tikzpicture}[nmdstd]
  \begin{scope}[scale=0.4, line width=1.4]
    \draw (6.51, 7.40) .. controls (6.96, 7.40) and (7.40, 7.19) .. 
          (7.40, 6.79) .. controls (7.40, 6.29) and (6.76, 6.18) .. (6.18, 6.18);
    \draw (6.18, 6.18) .. controls (5.57, 6.18) and (4.97, 6.51) .. (4.97, 7.07);
    \draw (4.97, 7.73) .. controls (4.97, 7.91) and (4.97, 8.09) .. (4.97, 8.28);
    \draw (4.97, 8.94) .. controls (4.97, 9.38) and (5.18, 9.82) .. 
          (5.58, 9.82) .. controls (6.08, 9.82) and (6.18, 9.19) .. (6.18, 8.61);
    \draw (6.18, 8.61) .. controls (6.18, 8.20) and (6.18, 7.80) .. (6.18, 7.40);
    \draw (6.18, 7.40) .. controls (6.18, 7.10) and (6.18, 6.81) .. (6.18, 6.51);
    \draw (6.18, 5.85) .. controls (6.18, 5.01) and (5.07, 4.97) .. (4.09, 4.97);
    \draw (3.43, 4.97) .. controls (2.41, 4.97) and (1.34, 4.65) .. 
          (1.34, 3.76) .. controls (1.34, 3.10) and (1.88, 2.55) .. (2.55, 2.55);
    \draw (2.55, 2.55) .. controls (2.96, 2.55) and (3.36, 2.55) .. (3.76, 2.55);
    \draw (3.76, 2.55) .. controls (5.85, 2.55) and (8.61, 2.55) .. (8.61, 3.43);
    \draw (8.61, 4.09) .. controls (8.61, 6.23) and (8.35, 8.61) .. (6.51, 8.61);
    \draw (5.85, 8.61) .. controls (5.56, 8.61) and (5.27, 8.61) .. (4.97, 8.61);
    \draw (4.97, 8.61) .. controls (4.36, 8.61) and (3.76, 8.28) .. (3.76, 7.73);
    \draw (3.76, 7.07) .. controls (3.76, 6.37) and (3.76, 5.67) .. (3.76, 4.97);
    \draw (3.76, 4.97) .. controls (3.76, 4.57) and (3.76, 4.17) .. (3.76, 3.76);
    \draw (3.76, 3.76) .. controls (3.76, 3.47) and (3.76, 3.18) .. (3.76, 2.88);
    \draw (3.76, 2.22) .. controls (3.76, 2.04) and (3.76, 1.86) .. (3.76, 1.67);
    \draw (3.76, 1.01) .. controls (3.76, 0.57) and (3.55, 0.13) .. 
          (3.16, 0.13) .. controls (2.66, 0.13) and (2.55, 0.76) .. (2.55, 1.34);
    \draw (2.55, 1.34) .. controls (2.55, 1.64) and (2.55, 1.93) .. (2.55, 2.22);
    \draw (2.55, 2.88) .. controls (2.55, 3.37) and (2.95, 3.76) .. (3.43, 3.76);
    \draw (4.09, 3.76) .. controls (5.60, 3.76) and (7.10, 3.76) .. (8.61, 3.76);
    \draw (8.61, 3.76) .. controls (9.27, 3.76) and (9.82, 3.22) .. 
          (9.82, 2.55) .. controls (9.82, 1.34) and (6.39, 1.34) .. (3.76, 1.34);
    \draw (3.76, 1.34) .. controls (3.47, 1.34) and (3.18, 1.34) .. (2.88, 1.34);
    \draw (2.22, 1.34) .. controls (0.86, 1.34) and (0.13, 2.84) .. 
          (0.13, 4.37) .. controls (0.13, 6.16) and (1.86, 7.40) .. (3.76, 7.40);
    \draw (3.76, 7.40) .. controls (4.17, 7.40) and (4.57, 7.40) .. (4.97, 7.40);
    \draw (4.97, 7.40) .. controls (5.27, 7.40) and (5.56, 7.40) .. (5.85, 7.40);
  \end{scope}

  \begin{scope}[scale=1.2, shift={(7.5, 1.55)}, line width=1.5]
    \coordinate (o) at (0, 0);
    \coordinate (a) at (-30:3);
    \coordinate (b) at (30:3);
    \coordinate (c) at (180:3);

    \coordinate (u) at ($(o)!0.5!(a)$);
    \coordinate (v) at ($(o)!0.5!(b)$);

    \coordinate (x) at ($(o)!0.66667!(a)$);
    \coordinate (z) at ($(o)!0.66667!(b)$);
    \coordinate (y) at ($(x)!0.5!(z)$);

    \coordinate (co) at ($(c)!0.5!(o)$);
    \coordinate (cco) at ($(c)!0.33334!(o)$);

    \coordinate (e) at (10:3.5);
    

    \filldraw[nmdlight] (o) -- (a) -- (b) -- cycle;
    \foreach \t in {0.0833333, 0.1666667,...,1.02}{
      \draw[color=nmddark, dashed] ($(o)!\t!(a)$) --  ($(o)!\t!(b)$);
    };

    \draw (o) -- (c);
    \draw (o) -- (a);
    \draw (o) -- (b);

    \foreach \n in {c, a, b}{
      \filldraw (\n) circle (2.7pt);
    };

    \draw[fill=white] (o) circle (2.5pt);
    
    \foreach \n in {co, u, v}{
      \node[fill=white, circle, draw, line width=0.8pt, inner sep=0.9pt] at (\n) {4};
    };

    \foreach \n in {cco, x, y, z}{
      \node[fill=white, circle, draw, line width=0.8pt, inner sep=0.9pt] at (\n) {5};
    };

    \begin{scope}[every node/.style={below=0.2}]
      \node at (c) {12};
      \node at (cco) {120};
      \node at (co) {73};
      \node at (o) {23};
      \node at (u) {70};
      \node at (x) {115};
      \node at (a) {4};      
    \end{scope}

    \begin{scope}[every node/.style={above=0.2}]  
      \node at (v) {71};
      \node at (z) {118};
      \node at (b) {7};
    \end{scope}

    \draw[->] (e) ..controls +(210:0.5) and +(50:0.6).. ($(y) + (0.07, 0.07)$);

    \node at ($(e) + (20:0.3)$) {116};
  \end{scope}  
\end{tikzpicture}

  \caption {For this triangulation of the exterior of $K13n585$, there
    are 11 connected lw-surfaces down to $\chi = -8$: three of genus 2
    ($N_4$, $N_7$, $N_{12}$), one of genus 3 ($N_{23}$), three of
    genus 4 ($N_{70} = N_4 + N_{23}$ and $N_{71} = N_7 + N_{23}$ and
    $N_{73} = N_{12} + N_{23}$), and four of genus 5
    ($N_{115} = 2 N_4 + N_{23}$ and $N_{116} = N_4 + N_7 + N_{23}$ and
    $N_{117} = 2 N_7 + N_{23}$ and $N_{120} = 2 N_{12} + N_{23}$).
    The complete list of isotopies between them is: $N_4 \sim N_7$,
    $N_{70} \sim N_{71}$, and $N_{115} \sim N_{116} \sim N_{118}$.
    The complex $\LW$ consists of the edge $B = [N_{12}, N_{23}]$ and
    the triangle $C = [N_{23}, N_{4}, N_{7}]$ shown above. Here
    $W_B = 0$, but for the faces $C$ and $E = [N_{4}, N_{7}]$, the
    subspaces $W_C$ and $W_E$ are both 1-dimensional; indeed,
    $W_C = W_E$ with the induced decomposition of $C$ into projective
    isotopy classes indicated by the dashed vertical lines.  }
  \label{fig:K13n585}
\end{figure}

An example of a non-normal isotopy of lw-surfaces occurs in the
13-tetrahedra triangulation:
\[
  \cT = \mathit{nvLAAvAPQkcdfgfhkmjlmklmwcadtfaaoaedrg}
\]
of the exterior of $K13n585$.  To determine $\LW$, we enumerated
normal and almost normal surfaces down to $\chi = -8$.  In this range,
there are 138 connected normal surfaces, 261 connected almost normal
surfaces with octagons, and 603 almost normal surfaces with tubes.  By
tightening the almost normal surfaces, we found there are 11 connected
essential lw-surfaces with $\chi \geq -8$ with four non-normal
isotopies among them.  Figure~\ref{fig:K13n585} shows the complex
$\LW$ which consists of an edge $B = [N_{12}, N_{23}]$ and a triangle
$C = [N_{23}, N_{4}, N_{7}]$.  For the face $C$ it is the surface
$N_{116}$ that plays the role of $F_C$ in steps
(\ref{step:lwfaces}--\ref{step:WC}) of
Algorithm~\ref{sec:actualalgorithm}, and the isotopies
$N_{115} \sim N_{116} \sim N_{118}$ are what determine the subspace
$W_C$.  Here, the subspace $W_E$ for $E = [N_4, N_7]$ is the same as
$W_C$.  In general, even if a face $E$ of $C$ is parallel to $W_C$, it
could be that $W_E$ is a proper subspace of $W_C$; see the example at
the start of Section~5 of \cite{To:isotopy} for more on this important
phenomenon.  Notice also that $\dep C$ is the complement of
$\{N_{23}\} \cup E$ and that the surfaces $N_{70}$ and $N_{71}$ are
\emph{projectively isotopic} to a surface carried by the interior of
$C$, but not isotopic to such a surface.

Because $W_E = W_C$, every isotopy class of essential lw-surface is
uniquely represented by a surface carried by $B \cup [N_{23}, N_4]$.
This allows us to easily compute that $B_M(x) = (-x^2 + 2x)/(x -
1)^2$. This is also what one gets from the triangulation:
\[
  \cS= \mathit{nvLALAwAQkedffgiijkmlmlmfvaeetcaangcbn}
\]
where $\LW[\cS]$ is two edges sharing a common vertex and there are no
isotopies between essential lw-surfaces.

\subsection{Barely large knots and those without meridional essential
  surfaces}

A striking contrast in Table~\ref{tab:overall} is that there are more
than 6,000 barely large manifolds in the cusped census yet only three
such knot exteriors.  Many constructions of closed essential surfaces
in knot exteriors come from tubing essential surfaces with meridional
boundary, and there are classes of knots where all closed essential
surfaces are of this form, including Montesinos knots
\cite{Oertel1984}, alternating knots \cite{Menasco}, and their
generalizations \cite{AdamsBrockBugbeeFaiginHustonPesikoff1992,
  Adams1994}.  A connected meridional surface $F$ in $M$ can be tubed
along $\partial M$ in two distinct ways, resulting in a pair of
disjoint surfaces; hence if both tubings are essential then $M$ is
very large.  However, barely large knot exteriors do exist: Baker
identified an infinite family of barely large knots with a single
incompressible genus $2$ surface in \cite[\sec 4.7.1]{Baker:thesis}.
Additionally, Adams-Reid \cite{AdamsReid1993} and Eudave-Mu\~noz
\cite{Eudave-Munoz2006} gave examples of closed essential surfaces
that cannot come from a meridional tubing construction.  Still, the
following appears to be new:

\begin{theorem}
  There exists a knot in $S^3$, namely $K15n153789$, whose exterior is
  large (indeed, barely large) and where the meridian is not the
  boundary slope of any essential surface.
\end{theorem}
Here, the knot $K15n153789$ is one of the three examples of barely
large knots we found; its exterior contains a unique essential
surface, which has genus 2.  We checked the boundary slope
condition by noting that there are no spunnormal surfaces with
meridional boundary slope in the triangulation:
\[
  \mathit{kLLLzPQkccfegjihijjlnahwdavhqk\_bBaB}
\]
of its exterior.

\subsection{Code and data} Complete data and the code used to compute
it are available at \cite{PaperData}.  Regina \cite{Regina} was used as
the underlying engine for triangulations and normal surface
computations, including enumeration of vertex and fundamental (almost)
normal surfaces, and Normaliz \cite{Normaliz} was used for computing
$B_M(x)$, with the whole computation taking place inside SageMath
\cite{SageMath} using the Python wrappings of these libraries.  The
code for dealing with almost normal surfaces with tubes, tightening
almost normal surfaces (with either tubes or octagons), and
implementing Algorithm~\ref{sec:actualalgorithm} was all completely
new.

To help validate our code, we started with a sample of 6,510 of the
manifolds from Table~\ref{tab:overall} and generated 5 random
triangulations of each. Then the complete algorithm was run on all
32{,}550 triangulations and the output compared to ensure that each
triangulation gave the same surface counts and other associated data.
This technique proved extremely effective at finding bugs in the code
(and, if we are being honest, our thinking), including subtle ones
that only manifest themselves in corner cases.  Additionally, we
compared our data to the lists of which knots are small/large from
\cite{BurtonCowardTillmann2013}; on the common set of 1,764 knots, our
data matched theirs exactly.

\FloatBarrier

\section{Patterns of surface counts by genus}
\label{sec.genuspat}

\begin{table}
\centering
  \tablefont
  \newcommand{\showBM}[2]{$\displaystyle \frac{#1}{#2}$}
  \newcommand{\aglinesp}{\addlinespace[0.4em]}
  \begin{tabular}{m{0.7\textwidth}cr}
    \toprule
    \multicolumn{1}{c}{$a_M(g)$} & $M$ & count \\
    \midrule
    4, 2, 4, 4, 8, 4, 12, 8, 12, 8, 20, 8, 24, 12, 16, 16, 32, 12, 36, 16     & $t09753$    & 1,473 \\ \aglinesp
    2, 1, 2, 2, 4, 2, 6, 4, 6, 4, 10, 4, 12, 6, 8, 8, 16, 6, 18, 8            & $t12198$    & 918   \\ \aglinesp
    0, 2, 0, 1, 0, 2, 0, 2, 0, 4, 0, 2, 0, 6, 0, 4, 0, 6, 0, 4                & $K14n11913$ & 259   \\ \aglinesp 
    6, 4, 8, 8, 16, 8, 24, 16, 24, 16, 40, 16, 48, 24, 32, 32, 64, 24, 72, 32 & $K12n605$   & 219   \\ \aglinesp
    8, 4, 8, 8, 16, 8, 24, 16, 24, 16, 40, 16, 48, 24, 32, 32, 64, 24, 72, 32 & $K11n73$    & 169   \\ \aglinesp
    0, 4, 0, 0, 0, 0, 0, 0, 0, 0, 0, 0, 0, 0, 0, 0, 0, 0, 0, 0                & $K14n13645$ & 148   \\ \aglinesp
    6, 9, 24, 37, 86, 87, 208, 220, 366, 386, 722, 602, 1168, 1039, 1498, %
        1564, 2514, 1993, 3484, 2924                                          & $K11n34$    & 139   \\ \aglinesp
    6, 7, 18, 29, 64, 73, 156, 177, 290, 321, 550, 521, 896, 865, 1236, %
        1297, 1950, 1731, 2714, 2499                                          & $K11n42$    & 131   \\ \aglinesp
    2, 0, 0, 0, 0, 0, 0, 0, 0, 0, 0, 0, 0, 0, 0, 0, 0, 0, 0, 0                & $o9_{37085}$ & 91    \\ \aglinesp
    0, 6, 0, 5, 0, 12, 0, 16, 0, 31, 0, 28, 0, 58, 0, 53, 0, 82, 0, 79        & $K15n93515$ & 82    \\ \aglinesp
    \bottomrule
  \end{tabular}

  \caption{The ten most common patterns of $a_M(g)$ for $2 \leq g \leq 21$,
    which together account for 3,629 (83.8\%) of the 4,330 very large
    manifolds. A sample manifold for each pattern is given in the second column, and
    the final column is the number of times the pattern appears.
  }
\label{tab:ag_common}
\end{table}

\begin{table}
\centering
  \tablefont
  \newcommand{\aglinesp}{\addlinespace[0.4em]}
  \begin{tabular}{%
    m{0.7\textwidth}
    >{$}c<{$}
    }
    \toprule
    \multicolumn{1}{c}{$a_M(g)$} & M \\
    \midrule
    8, 14, 46, 89, 224, 305, 674, 905, 1536, 1955, 3326, 3771,
    6150, 7019, 9850, 11611, 16714, 17767, 25490, 27415 & K12n214 \\ \aglinesp
    8, 16, 54, 98, 264, 318, 806, 984, 1794, 2098, 3994, 4074, 7368,
    7632, 11552, 12976, 20114, 19396, 30670, 30550 & K12n210 \\ \aglinesp
    12, 21, 61, 109, 261, 320, 721, 880, 1480, 1762, 3094, 3115,
    5429, 5666, 8019, 9086, 13596, 13059, 20062, 19841 & K13n3763 \\ \aglinesp
    10, 25, 71, 140, 352, 473, 1058, 1386, 2389, 2939, 5152, 5585,
    9422, 10311, 14887, 17057, 25304, 25573, 38238, 39603 & K15n15582 \\ \aglinesp
    12, 16, 51, 99, 235, 345, 711, 999, 1649, 2209, 3551, 4319,
    6593, 7919, 10971, 13231, 18275, 20555, 28063, 31485 & K15n15220 \\ \aglinesp
    8, 18, 57, 110, 270, 356, 785, 1013, 1737, 2092, 3667, 3942,
    6614, 7134, 10397, 11710, 17426, 17422, 26131, 26891 & K15n23198 \\ \aglinesp
    12, 34, 110, 216, 532, 708, 1558, 2018, 3462, 4176, 7314, 7876,
    13204, 14256, 20778, 23404, 34820, 34832, 52226, 53766 & K13n3838 \\ \aglinesp
    12, 30, 109, 231, 549, 861, 1737, 2511, 4059, 5643, 8859, 10941,
    16623, 20229, 27303, 33729, 46215, 52455, 71079, 80271 & K15n33595 \\ \aglinesp
    10, 21, 73, 143, 385, 513, 1224, 1605, 2870, 3542, 6409, 7010, 12051,
    13231, 19463, 22436, 33614, 34307, 51700, 53862 & K13n2458 \\ \aglinesp
    \bottomrule
  \end{tabular}

  \caption{The eight of the most complicated patterns of $a_M(g)$ for $2 \leq g \leq 21$.
    These all come from examples where $\dim \LW \geq 4$.}
\label{tab:ag_messy}
\end{table}

We now return to the question of counting \emph{connected} essential
surfaces in a given \3-manifold in terms of their genus.  As we know
how to count all essential surfaces by Euler characteristic given the
complex $\LW$, we approach this by identifying the connected surfaces
in that larger count.  This problem has an arithmetic flavor, and is
related to counting primitive lattice points, as well as to the work
of Mirzakhani discussed in Section~\ref{sec:motivation}.

Let $a_M(g)$ denote the number of isotopy classes of connected
essential surfaces of genus $g$.  For each of the 4,330 very large
examples in Table~\ref{tab:overall}, we computed the first 20 values
of $a_M(g)$ starting from $\LW$ as follows.  Let $g$ be fixed. For
each face $C$ of $\LW$, let $\Ctil = \R_{\geq 0} \cdot C \subset \ST$
be the corresponding cone.  We used Normaliz \cite{Normaliz} to find
all lattice points carried by the interior of the rational polytope
$\setdef{\xx \in \Ctil}{\chi(\xx)= 2 - 2g}$.  For each corresponding
normal surface, we checked connectivity using Algorithm 9.4 of
\cite{JacoTollefson1995}.  As $W_C = 0$ for all these examples, the
number of such lattice points corresponding to connected surfaces is
the contribution of the interior $\intC$ to $a_M(g)$.

We found 94 distinct patterns for $\big(a_M(2), \ldots, a_M(21)\big)$.
Table~\ref{tab:ag_common} lists the most common patterns and
Table~\ref{tab:ag_messy} gives the most complicated. For one manifold
exhibiting each pattern, we computed additional values of $a_M(g)$,
nearly always up to at least $g = 50$ and in more than 40 cases up to
$g = 200$. This data is available at \cite{PaperData}, where the
largest single value is $a_M(51) = 3,072,351$ for the
exterior $M$ of $K15n33595$.

From now on, we work with $n = g - 1$ rather than $g$
as the index for the count of connected surfaces, 
and so define $\atil_M(n) = a_M(n + 1)$; this simplifies the arithmetic,
for example giving $\atil_M(n) \leq b_M(-2n)$ rather than $a_M(g) \leq
b_M(2 - 2g)$.

\subsection{Independence of $a_M(g)$ and $B_M(x)$}
\label{sec:indep}

We next give four examples showing that neither $\atil_M(n)$ nor
$B_M(x)$ determines the other.  We start with two manifolds with the
same $B_M(x)$ but different $\atil_M(n)$.  Let $A$ and $B$ be the
exteriors of the knots $K14n22185$ and $K13n586$ respectively, and we
use $\cT$ and $\cS$ to denote their standard triangulations. Both
$\LW$ and $\LW[\cS]$ consist of a single edge $C$ whose vertices
correspond to genus 2 surfaces $F$ and $G$.  Moreover, the lattice
points in the cone over $C$ are simply $u\vF + v\vG$ for
$u, v \in \N$. Thus, the surfaces with $\chi = -2n$ are the lattice
points in $\N^2$ on the line $x + y = n$, which gives $b_M(-2n) = n +
1$ and hence $B_M(x) = (-x^2 + 2x)/(x - 1)^2$.

For $\cT$, the surfaces $F$ and $G$ can be made disjoint after a normal
isotopy, and hence every normal surface carried by $C$ is a disjoint
union of parallel copies of $F$ and $G$.  Thus the only connected
essential surfaces in $A$ are $F$ and $G$, giving $\atil_A(1) = 2$ and
$\atil_A(n) = 0$ for $n > 1$. In contrast, we find that the first 30
values of $\atil_B(n)$ are:
\[
  2, 1, 2, 2, 4, 2, 6, 4, 6, 4, 10, 4, 12, 6, 8, 8, 16, 6, 18, 8, 12,
  10, 22, 8, 20, 12, 18, 12, 28, 8
\]
Now, if $u\vF + v\vG$ is connected then $\gcd(u, v) = 1$.  The above
data is consistent with the converse being true, or equivalently
$\atil_B(n)$ is exactly the number of primitive lattice points in
$\N^2$ on the line $x + y = n$, which is the Euler totent function
$\phi(n)$ when $n > 1$.  This pattern continues for all $n \leq 500$,
so we may safely posit:
\begin{conjecture}
  \label{conj:K13n586}
  For the exterior $B$ of the knot $K13n586$, one has $\atil_B(n) = \phi(n)$ for
  all $n > 1$.
\end{conjecture}
Since $\atil_B(1) = 2$, equivalently the conjecture is that
$\atil_B(n) = \epsilon(n) + \phi(n)$ for all $n \geq 1$ where
$\epsilon(n)$ is $1$ when $n = 1$ and $0$ otherwise.  This count of
primitive lattice points can be related to the corresponding count of
all lattice points via the M\"obius inversion formula, making
Conjecture~\ref{conj:K13n586} equivalent to:
\begin{equation}
  \label{eq.aBg}
  \atil_B(n) = \sum_{d|n} \mu\left(\frac{n}{d}\right) (d + 1) 
\end{equation}
for all $n \geq 1$, where $\mu$ is the M\"obius function.
\emph{Note added in proof:} Lee \cite{Lee2021} has recently
proved Conjecture~\ref{conj:K13n586} by analyzing the number of
connected components of $u\vF + v\vG$ using the method of
\cite{AgolHassThurston2006}.

A pair with the same $\atil_M(n)$ but different $B_M(x)$ are the
census manifolds $X = v3394$ and $Y = o9_{43058}$.  Both have exactly
four connected essential surfaces, all of genus two, but
\[
  B_X(x) = \frac{-2x^2 + 4x}{(x - 1)^2} \mtext{and}   B_Y(x) = \frac{-x^2 + 4x}{(x - 1)^2}
\]
For the standard triangulations $\cX$ and $\cY$ of $X$ and $Y$, the
complexes $\LW[\cX]$ and $\LW[\cY]$ are quite different: the first
consists of an edge and a disjoint vertex, but $\LW[\cY]$ consists of
two edges sharing a common vertex.  All vertex surfaces are connected,
and so the vertices of $\LW[\cX]$ and $\LW[\cY]$ correspond to three
of the four essential genus 2 surfaces; in both cases, the fourth is
hiding as a fundamental surface in the interior of an edge.

\subsection{Regular genus counts and the Lambert series}
\label{sec:reg}

For the manifold $B = K13n586$ in Section~\ref{sec:indep}, while the
count $\atil_B(n)$ does not have a short generating function, from
(\ref{eq.aBg}) we see its M\"obius transform
$p(n) = \sum_{d|n} \atil_B(d)$ is a polynomial, specifically
$p(n) = n + 1$.  This motivates our next definition.  Recall that
Dirichlet convolution on arithmetic functions
$f, g \maps \Z_{\geq 1} \to \C$ is defined by
$(f \ast g)(n) = \sum_{d|n} f(n/d) g(d)$.  We say that $\atil_M(n)$ is
\emph{regular} if $1 \ast \atil_M$ has a short generating series.
Equivalently, if we set $p_M = 1 \ast \atil_M$, regularity is
equivalent to $p_M(n)$ being a quasi-polynomial for all large $n$.
Thus, when $\atil_M(n)$ is regular, we can use M\"obius inversion
$ \atil_M(n) = (\mu \ast p_M)(n) = \sum_{d|n}
\mu\left(\frac{n}{d}\right) p_M(d) $ to compute $\atil_M$ from the
simpler $p_M$, as we did in (\ref{eq.aBg}).  In the language of
generating functions, the count $\atil_M$ is regular if and only if
its Lambert series \be
\label{eq.lambert}
\LA_M(x) = \sum_{n=1}^\infty \atil_M(n) \frac{x^n}{1-x^n} \,.
\ee
is short, since the coefficients of this series are precisely $1 \ast \atil_M$.

\begin{table}
  \centering
  \tablefont
  \newcommand{\showLA}[2]{$\displaystyle \frac{#1}{#2}$}
  \newcommand{\LAlinesp}{\addlinespace[0.6em]}
  \begin{tabular}{lcc>{$}c<{$}}
    \toprule
    \multicolumn{1}{c}{$LA_M(x)$} & per & $\ell_1$ & M \\
    \midrule

    \showLA{-2x^2 + 4x}{(x - 1)^2} & 1 & 10 & t09753 \\ \LAlinesp
    \showLA{-x^2 + 2x}{(x - 1)^2} & 1 & 7 & t12198 \\ \LAlinesp
    \showLA{-x^4 + 2x^2}{(x - 1)^2  (x + 1)^2} & 2 & 7 & K14n11913 \\ \LAlinesp
    \showLA{-2x^2 + 6x}{(x - 1)^2} & 1 & 12 & K12n605 \\ \LAlinesp
    \showLA{-4x^2 + 8x}{(x - 1)^2} & 1 & 16 & K11n73 \\ \LAlinesp
    \showLA{-4x^6 + 2x^5 + 16x^4 + 4x^3 - 14x^2 - 6x}{(x - 1)^3  (x + 1)^3}
        & 2 & 54 & K15n67261 \\ \LAlinesp
    \showLA{-2x^8 - 4x^7 - 2x^6 + 4x^5 + 9x^4 + 6x^3 + 2x^2}{(x - 1)^2  (x + 1)^2  (x^2 + x + 1)^2}
        & 6 & 45 & K15n129923 \\ \LAlinesp
    \showLA{-2x^8 - 4x^7 - 2x^6 + 6x^5 + 13x^4 + 8x^3 + 2x^2}{(x - 1)^2  (x + 1)^2  (x^2 + x + 1)^2}
        & 6 & 53 & K15n138922 \\ \LAlinesp
    \bottomrule
  \end{tabular}

  \caption{Eight examples of our conjectured Lambert series $LA_M(x)$
    for manifolds where $\atil_M(n)$ appears regular.  The first five
    are from Table~\ref{tab:ag_common} and the last three are
    among the most complicated we found.}
\label{tab:lambert}
\end{table}

\subsection{Examples}

Of the 94 observed patterns for $\atil_M(n)$ in our sample, we conjecture
that exactly 54 of them are regular, including 7 of the 10 manifolds
in Table~\ref{tab:ag_common},  with the exceptions being $K11n34$,
$K11n42$, and $K15n93515$.  Examples of our conjectured formulae for
$LA_M(x)$ are given in Table~\ref{tab:lambert}.  In contrast, we believe
all the examples in Table~\ref{tab:ag_messy} are irregular.

One example where the count appears irregular, though still highly
structured, is $M = o9_{41176 }$.  Specifically we conjecture that
$\atil_M(n)$ is equal to $f(n) = \phi(n) + 1$ for $n \geq 2$ (here
$\atil_M(1) = 5$).  While $f$ is quite simple, we have
$1 \ast f = n + \sigma_0$ where $\sigma_0(n)$ is the number of
divisors of $n$, and $\sigma_0(n)$ does not have a short generating
function.  For the standard triangulation $\cT$ of $M$, the complex
$\LW$ is:
\begin{center}
\begin{tikzpicture}[scale=0.75, nmdstd, line width=1.5pt, font=\small]
  \draw (0, 0) -- (10, 0);
  \begin{scope}[radius=3pt, every node/.style={above=4pt}]
    \filldraw (0, 0) circle node {$N_{9}$};
    \filldraw (2, 0) circle node {$N_{7}$};
    \filldraw[fill=white] (4, 0) circle node {$N_{16}$};
    \filldraw (6, 0) circle node {$N_{5}$};
    \filldraw[fill=white] (8, 0) circle node {$N_{14}$};
    \filldraw (10, 0) circle node {$N_{10}$};
  \end{scope}
  \foreach \i in {1,...,5}{
    \node[below=2pt] at ($2*(\i, 0)-(1,0)$) {$C_{\i}$};
  }
\end{tikzpicture}
\end{center}
using the conventions of \cite{PaperData}.  Here the vertex surfaces
$N_{14}$ and $N_{16}$ have genus 3 and the others have genus 2.  We
conjecture that the faces contribute to $\atil_M$ as follows:
\begin{enumerate}
\item The interior of $C_1$ carries a single connected surface $N_8$
  which has genus 2.  Here $N_7 + N_9 = 2 N_8$.

\item The interior of $C_2$ carries no connected surfaces as $N_7$ and
  $N_{16}$ are disjoint.

\item The interior of $C_3$ carries a unique surface of genus $g$ for
  each $g \geq 4$, namely $(g - 3) N_5 + N_{16}$.  It is this face
  that contributes the $+1$ to $\atil_M(n)$.

\item The connected surfaces carried by $C_4$ are exactly $u N_5 + v
  N_{14}$ for $u, v > 0$ and $\gcd(u, v) = 1$.  The situation is the
  same for $C_5$, with $N_5$ replaced with $N_{10}$.  Together, these
  faces contribute the $\phi(n)$ to $\atil_M(n)$.
\end{enumerate}
The manifold $N = t12071$ is similar in that $\atil_N(n)$ is irregular
but $\atil_N(n) - 4$ appears regular.

In our sample, a simple example where $\atil_M$ appears irregular and where
we cannot glean any other structure is $W = o9_{42517}$.  The first
50 values of $\atil_W(n)$ are:
\begin{align*}
  & 6, 4, 10, 14, 26, 26, 52, 46, 76, 76, 118, 96, 172, 136, 194, 196,
    286, 212, 354, 274, 388, \\
  & 360, 506, 378, 604, 490, 634, 574, 820, 568, 948, 728, 946, 846,
    1122, 864, 1356, 1040, \\
  & 1316, 1146, 1644, 1140, 1800, 1392, 1716, 1570, 2136, 1506, 2332, 1752
\end{align*}
and \cite{PaperData} has all values to $n=200$.  With its
usual triangulation $\cW$, the complex $\LW[\cW]$ is a triangle whose
vertex surfaces $N_3, N_9$, and $N_{11}$ all have genus
2.  Here the edge $[N_3, N_9]$ appears to carry a single connected
surface in its interior, which has genus 2, and the same for
$[N_9, N_{11}]$.  The remaining edge $[N_3, N_{11}]$ appears to
contribute $2 \phi(n)$ to $\atil_W(n)$ for $n > 1$, and the interior
of the triangle contributes the mysterious:
\[
  0, 2, 6, 10, 18, 22, 40, 38, 64, 68, 98, 88, 148, 124, 178, 180,
  254, 200, 318, 258, 364, 340
\]


\subsection{Asymptotics of genus counts}

\begin{figure}
  \centering
 \pgfkeys{/matplotlibfigure, default, width=0.8\textwidth}%
 \begin{tikzoverlay}[width=\matplotlibfigurewidth]{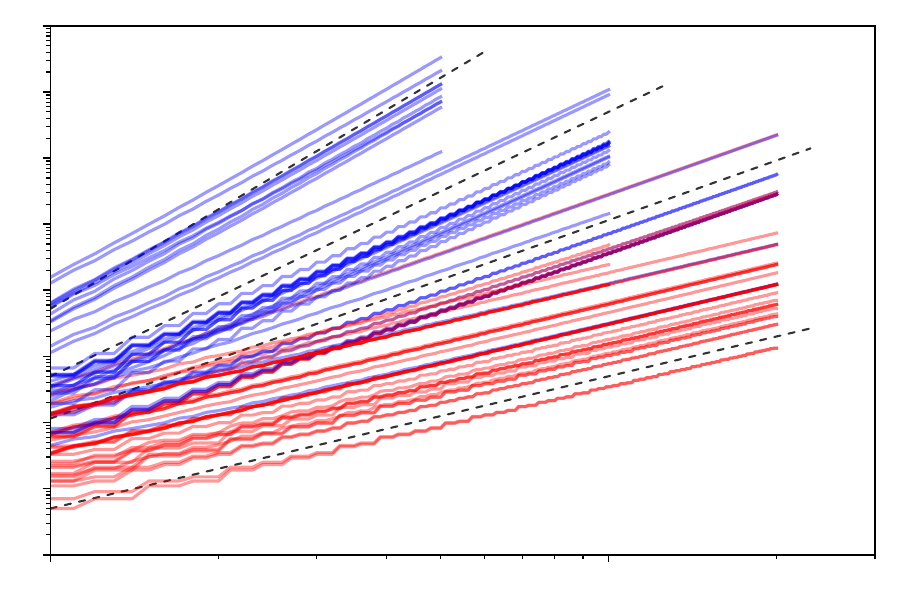}[\matplotlibfigurefont]
  \draw (5.620370, 3.416667) node[below] {${10^{1}}$};
  \draw (67.624734, 3.416667) node[below] {${10^{2}}$};
  \draw (4.000000, 5.037037) node[left] {${10^{0}}$};
  \draw (4.000000, 12.381773) node[left] {${10^{1}}$};
  \draw (4.000000, 19.726508) node[left] {${10^{2}}$};
  \draw (4.000000, 27.071244) node[left] {${10^{3}}$};
  \draw (4.000000, 34.415979) node[left] {${10^{4}}$};
  \draw (4.000000, 41.760715) node[left] {${10^{5}}$};
  \draw (4.000000, 49.105451) node[left] {${10^{6}}$};
  \draw (4.000000, 56.450186) node[left] {${10^{7}}$};
  \draw (4.000000, 63.794922) node[left] {${10^{8}}$};
  \node at (92.0, 32.75) {$s = 2$};
  \node at (92.0, 52.75) {$s = 3$};
  \node at (78.0, 58.75) {$s = 4$};
  \node at (58.0, 61.5) {$s = 5$};
  \node at (50, -1) {$n$};
  \node[left] at (-3,  34.415979) {$\abar_M$};
\end{tikzoverlay}

  \caption{This log-log plot shows the sequence $\{\abar_M(n)\}$ for
    76 manifolds, up to $n = 50, 100$, or $200$ depending.  Those
    coming from conjecturally regular $\atil_M(n)$ are in red whereas
    the likely irregular ones are in blue.  The dotted lines plot
    $c_s n^s$ for the indicated $s$ and some choice of $c_s$.  Each of
    $\abar_M(n)$ appears nearly parallel to one of these
    lines, consistent with $\abar_M(n)$ being asymptotic to $c n^s$ as
    $n \to \infty$ for some integer $s$ and $c_s > 0$.  }
  \label{fig:asymp}
\end{figure}

\begin{figure}
  \centering
 \pgfkeys{/matplotlibfigure, default, width=0.7\textwidth}%
 \begin{tikzoverlay}[width=\matplotlibfigurewidth]{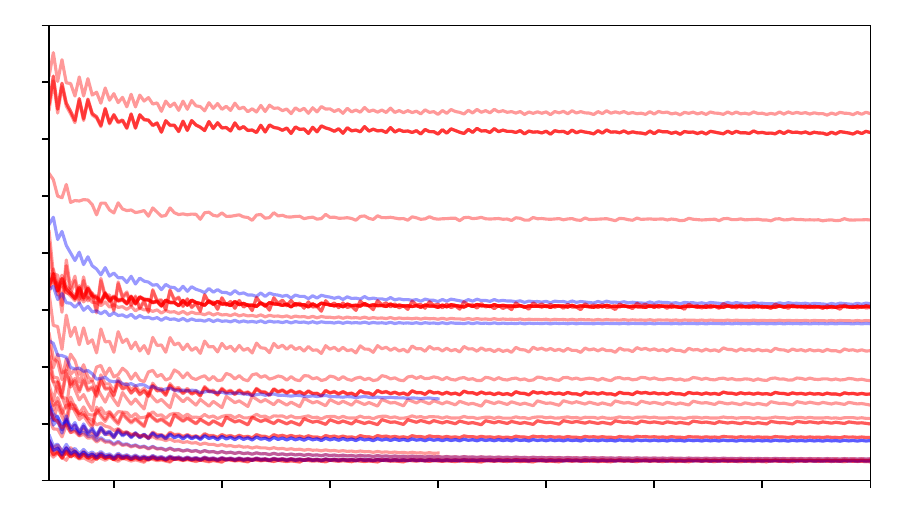}[\matplotlibfigurefont]
  \newcommand{\zpad}{\phantom{0}};
  \draw (12.650564, 3.305990) node[below] {$25$};
  \draw (24.658352, 3.305990) node[below] {$50$};
  \draw (36.666139, 3.305990) node[below] {$75$};
  \draw (48.673927, 3.305990) node[below] {$100$};
  \draw (60.681715, 3.305990) node[below] {$125$};
  \draw (72.689503, 3.305990) node[below] {$150$};
  \draw (84.697290, 3.305990) node[below] {$175$};
  \draw (96.705078, 3.305990) node[below] {$200$};
  \draw (3.825521, 4.926360) node[left]  {$\zpad 0.0$};
  \draw (3.825521, 11.250181) node[left] {$\zpad 0.1$};
  \draw (3.825521, 17.574002) node[left] {$\zpad 0.2$};
  \draw (3.825521, 23.897823) node[left] {$\zpad 0.3$};
  \draw (3.825521, 30.221644) node[left] {$\zpad 0.4$};
  \draw (3.825521, 36.545464) node[left] {$\zpad 0.5$};
  \draw (3.825521, 42.869285) node[left] {$\zpad 0.6$};
  \draw (3.825521, 49.193106) node[left] {$\zpad 0.7$};
  \draw (3.825521, 55.516927) node[left] {$\zpad 0.8$};
  \node at (48.673927, -4) {$n$};
  \node[left] at (-3, 30.221644) {$\displaystyle \frac{\abar_M(n)}{n^s}$};
\end{tikzoverlay}

  %
 \pgfkeys{/matplotlibfigure, default, width=0.7\textwidth}%
 \begin{tikzoverlay}[width=\matplotlibfigurewidth]{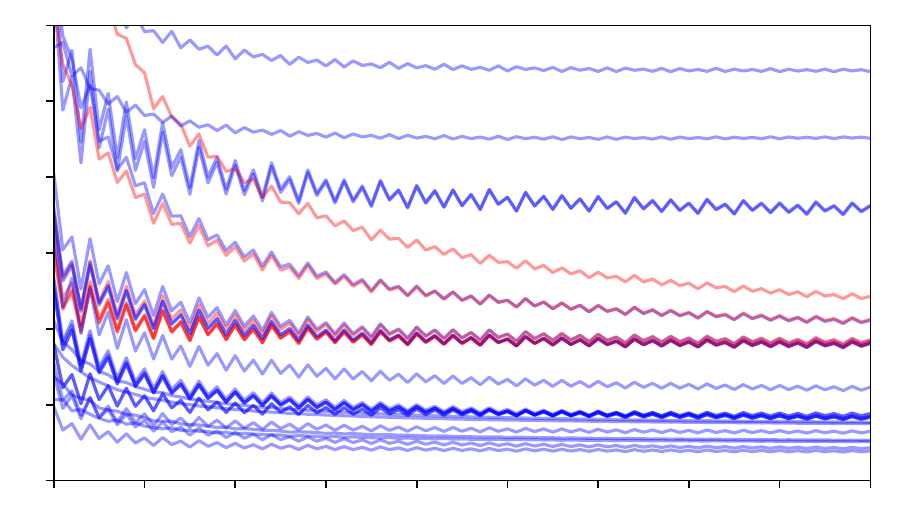}[\matplotlibfigurefont]
  \draw (5.975839, 3.305990) node[below] {$10$};
  \draw (16.056866, 3.305990) node[below] {$20$};
  \draw (26.137892, 3.305990) node[below] {$30$};
  \draw (36.218919, 3.305990) node[below] {$40$};
  \draw (46.299945, 3.305990) node[below] {$50$};
  \draw (56.380972, 3.305990) node[below] {$60$};
  \draw (66.461998, 3.305990) node[below] {$70$};
  \draw (76.543025, 3.305990) node[below] {$80$};
  \draw (86.624052, 3.305990) node[below] {$90$};
  \draw (96.705078, 3.305990) node[below] {$100$};
  \draw (4.355469, 4.926360) node[left] {$0.00$};
  \draw (4.355469, 13.358121) node[left] {$0.02$};
  \draw (4.355469, 21.789882) node[left] {$0.04$};
  \draw (4.355469, 30.221644) node[left] {$0.06$};
  \draw (4.355469, 38.653405) node[left] {$0.08$};
  \draw (4.355469, 47.085166) node[left] {$0.10$};
  \draw (4.355469, 55.516927) node[left] {$0.12$};
  \node at (48.673927, -4) {$n$};
  \node[left] at (-3, 30.221644) {$\displaystyle \frac{\abar_M(n)}{n^s}$};
\end{tikzoverlay}

  \caption{Using the predicted asymptotic exponent $s$ for each
    sequence $\abar_M(n)$ from Figure~\ref{fig:asymp}, we plot
    $\abar_M(n)/n^s$ to test Conjecture~\ref{conj:asymp}.  Here, the
    top plot shows those where $s=2,3$ and the bottom where $s=4,5$;
    again, red and blue correspond to (conjecturally) regular versus
    irregular sequences.  For better readability, 10 sequences that
    lie above the given vertical scales are omitted, 9 from the top
    plot (all but one regular) and 1 from the bottom; these look very
    similar to the 66 sequences shown.}
    \label{fig:asymp2345}
\end{figure}

We now explore the asymptotics of the sequences $\atil_M(n)$.  Since
$\atil_M(n) \leq b_M(-2n)$ and the latter grows polynomially, it is
natural to ask whether $\atil_M(n)$ does so as well.  Even in the
regular case, the sequence $\atil_M(n)$ depends arithmetically on the
divisors of $n$, so it is better to study the smoothed sequence:
\begin{equation}
\label{eq.asmooth}
\abar_M(n) = \sum_{k \leq n} \atil_M(k) 
\end{equation}
Of the 94 observed patterns for $\atil_M(n)$, there are 14 where
$\atil_M(n) = 0$ for all large $n$ and 4 where we were only able to
compute up to $\atil_M(20)$; we consider only the remaining 76. The
plots in Figures~\ref{fig:asymp} and~\ref{fig:asymp2345} together
suggest:
\begin{conjecture}\label{conj:asymp}
  Suppose $M$ as in Theorem~\ref{thm.main}.  Then either $\atil_M(n) =
  0$ for all large $n$ or there exists $s \in \N$ such
  that $\lim_{n \to \infty} \abar_M(n)/n^s$ exists and is positive.
\end{conjecture}

In fact, Conjecture~\ref{conj:asymp} holds whenever $\atil_M(n)$
is regular and the corresponding quasi-polynomial has constant leading
term as we now show; this includes all the conjecturally regular
examples in our sample.

\begin{lemma}
  \label{lem.areg}
  Suppose $\atil_M(n)$ is regular and the corresponding $p_M(n)$ has
  constant leading term, with
  $p_M(n)= c_r n^r + O\left(n^{r-1}\right)$ for some $r \geq 1$ and
  positive $c_r \in \Q$.  Then
  \begin{equation}
    \label{eq.aslim}
    \lim_{n \to \infty} \frac{1}{n^{r+1}}\abar _M(n)
    = \frac{c_r}{r+1} \frac{1}{\zeta(r+1)} 
  \end{equation}
  where $\zeta(s)$ is the Riemann $\zeta$-function. 
\end{lemma}
\begin{proof}
  As $\sum_{k = 1}^{m} k^s = \frac{1}{s + 1} m^{s +
    1} + O(m^s)$, we see $\sum_{k \leq m} p_M(k) = \frac{c_r}{r + 1} m^{r + 1} + O(m^r)$.
  Now
  \newcommand{\floornd}{\left\lfloor \frac{n}{d} \right\rfloor}
  \begin{align*}
    \abar_M(n) 
    \ &= \sum_{\ell \leq n} \ \sum_{d | \ell} \mu(d) p_M(\ell/d)
    \ = \sum_{d \cdot k \leq n} \mu(d) p_M(k)
    \ = \sum_{d \leq n} \Big( \mu(d) \sum_{k \leq \lfloor n/d \rfloor}
        p_M(k) \Big) \\
    &= \sum_{d \leq n} \mu(d) \left( \frac{c_r}{r + 1}
        \floornd^{r + 1} + O\left(\floornd^{r}\right) \right)
      = \sum_{d \leq n} \mu(d) \left( \frac{c_r}{r + 1}
        \frac{n^{r + 1}}{d^{r + 1}} + O\left(\frac{n^{r}}{d^{r}}\right) \right)
  \end{align*}
  where we have used that $\lfloor n/d\rfloor = n/d + O(1)$ and hence by
  the binomial theorem $\lfloor n/d\rfloor^{r + 1} = n^{r+1}/d^{r+1} +
  O(n^r/d^r)$.  Thus
  \[
    \frac{\abar_M(n)}{n^{r + 1}} \ = \sum_{d \leq n} \mu(d) \left( \frac{c_r}{r + 1}
        \frac{1}{d^{r + 1}} + O\left(\frac{1}{n d^{r}}\right) \right)
      = \left( \frac{c_r}{r+1} \sum_{d \leq n} \frac{\mu(d)}{d^{r+1}} \right)
      + O\left(\frac{\log(n)}{n}\right)
  \]
  where we have used in the last step that
  $\sum_{d \leq n} 1/d^r \leq \sum_{d \leq n} 1/d \approx \log(n)$.
  Now as $r \geq 1$, we have
  $\sum_{d \leq n} \mu(d)/d^{r+1} = 1/\zeta(r + 1) + o(n)$ by
  \cite[Section 11.4]{Apostol}, and so the lemma follows.
\end{proof}

\begin{remark}
  If one allows $r = 0$, then (\ref{eq.aslim}) still holds if you
  interpret the righthand side as 0, seeing that $\zeta$ has a pole at
  1; we leave the details in this case to the reader.
\end{remark}

{\RaggedRight
  \bibliographystyle{nmd/math}
  \small
  \bibliography{biblio}
}

\pagebreak

\end{document}